\documentclass[preprint,11pt,3p,authoryear]{elsarticle}

\usepackage[bookmarks=false,colorlinks]{hyperref}
\AtBeginDocument{%
	\hypersetup{
		linkcolor=blue,   
		citecolor=red,
		urlcolor=magenta}}
\usepackage{lmodern}
\usepackage[english]{babel}
\usepackage[utf8]{inputenc}
\usepackage{amsthm}
\usepackage{amsmath,amsfonts,amssymb}
\usepackage[normalem]{ulem}
\usepackage{subfig}
\usepackage[ruled,vlined,linesnumbered]{algorithm2e}
\usepackage{graphicx}
\usepackage{xcolor}
\usepackage{float}
\usepackage{bm}
\usepackage{natbib}
\usepackage{xspace}
\usepackage[tableposition=top]{caption}
\usepackage{textcomp}
\usepackage{booktabs}
\usepackage{gensymb}
\usepackage[toc,page]{appendix}

\newtheorem{lemma}{Lemma}

\newtheorem{theorem}{Theorem}
\newtheorem{prop}{Proposition}
\newtheorem{definition}{Definition}
\newtheorem{modl}{Model}
\newenvironment{model}{\begin{samepage}\begin{modl}}{\end{modl}\end{samepage}}

\newcommand{\sig}{\sigma_{ij}}

\newcommand{\viz}{\widehat{v}_{i}}
\newcommand{\vjz}{\widehat{v}_{j}}
\newcommand{\qi}{q_{i}}
\newcommand{\qj}{q_{j}}
\newcommand{\ti}{\theta_{i}}
\newcommand{\tj}{\theta_{j}}
\newcommand{\qmax}{\overline{q}_i}
\newcommand{\qmin}{\underline{q}_i}
\newcommand{\qub}{\overline{\bm{q}}}
\newcommand{\qlb}{\underline{\bm{q}}}

\newcommand{\pzx}{\widehat{x}_{ij}}
\newcommand{\pzy}{\widehat{y}_{ij}}
\newcommand{\tiz}{\widehat{\theta}_{i}}
\newcommand{\tjz}{\widehat{\theta}_{j}}
\newcommand{\tlb}{\underline{\bm{\theta}}}
\newcommand{\tub}{\overline{\bm{\theta}}}
\newcommand{\tmin}{\underline{\theta}_i}
\newcommand{\tmax}{\overline{\theta}_i}
\newcommand{\dx}{\delta_{i,x}}
\newcommand{\dy}{\delta_{i,y}}
\newcommand{\tdx}{\tilde{\delta}_{i,x}}
\newcommand{\tdy}{\tilde{\delta}_{i,y}}

\newcommand{\z}{z_{ij}}

\newcommand{\ii}{\mathrm{i}}
\newcommand{\tm}{t^{\text{min}}_{ij}}

\newcommand{\vijx}{v_{ij,x}}
\newcommand{\vijy}{v_{ij,y}}
\newcommand{\vpij}{\varphi_{ij}}
\newcommand{\rik}{\rho_{ik}}
\newcommand{\rjk}{\rho_{jk}}
\newcommand{\R}{\mathbb{R}}
\newcommand{\B}{\mathcal{B}_{ij}}
\newcommand{\C}{\mathcal{C}_{ij}}
\newcommand{\A}{\mathcal{A}}

\newcommand{\Z}{\mathcal{Z}}
\newcommand{\Zi}{\mathcal{Z}_i}
\newcommand{\Zj}{\mathcal{Z}_j}
\renewcommand{\P}{\mathcal{P}}
\newcommand{\PF}{\P_{\text{F}}}
\newcommand{\PS}{\P_{\text{S}}}
\newcommand{\PI}{\P_{\text{I}}}
\newcommand{\PP}{\mathcal{P}^3}

\newcommand{\PPS}{\PP_{\text{S}}}

\newcommand{\gx}{\Gamma_{i,x}}
\newcommand{\gy}{\Gamma_{i,y}}
\newcommand{\axk}{\alpha_{i,x}^k}
\newcommand{\ayk}{\alpha_{i,y}^k}
\newcommand{\bxk}{\beta_{i,x}^k}
\newcommand{\byk}{\beta_{i,y}^k}
\newcommand{\sxk}{s_{i,x}^k}
\newcommand{\syk}{s_{i,y}^k}
\newcommand{\dxlk}{\underline{\delta}_{i,x}^k}
\newcommand{\dxuk}{\overline{\delta}_{i,x}^k}
\newcommand{\dylk}{\underline{\delta}_{i,y}^k}
\newcommand{\dyuk}{\overline{\delta}_{i,y}^k}

\newcommand{\Python}{\textsc{Python}\xspace}
\newcommand{\Cplex}{\textsc{Cplex}\xspace}

\journal{European Journal of Operational Research}

\begin{document}

\begin{frontmatter}
\textcolor{red}{\title{Disjunctive linear separation conditions and mixed-integer formulations for aircraft conflict resolution}}

\author[UNSW]{Fernando H. C. Dias}
\author[LANL]{Hassan Hijazi}
\author[UNSW]{David Rey\corref{mycorrespondingauthor}}
\ead{drey@unsw.edu.au}
\cortext[mycorrespondingauthor]{Corresponding author}
\address[UNSW]{School of Civil and Environmental Engineering, UNSW Sydney, Sydney, NSW, 2052, Australia}
\address[LANL]{Los Alamos National Laboratory, P.O. Box 1663, Los Alamos, NM 87545, USA}

\begin{abstract}
We address the aircraft conflict resolution problem in air traffic control. We introduce new mixed-integer programming formulations for aircraft conflict resolution with speed, heading and altitude control which are based on disjunctive linear separation conditions. We first examine the two-dimensional aircraft conflict resolution problem with speed and heading control represented as continuous decision variables. We show that the proposed disjunctive linear separation conditions are equivalent to the traditional nonlinear conditions for aircraft separation. Further, we characterize conflict-free pairwise aircraft trajectories and propose a simple pre-processing algorithm to identify aircraft pairs which are either always conflict-free, or which cannot be separated using speed and heading control only. We then incorporate altitude control and propose a lexicographic optimization formulation that aims to minimize the number of flight level changes before resolving outstanding conflicts via two-dimensional velocity control. The proposed mixed-integer programming formulations are nonconvex, and we propose convex relaxations, decomposition methods and constraint generation algorithms to solve the two-dimensional and lexicographic optimization formulations to guaranteed optimality. Numerical experiments on four types of conflict resolution benchmarking instances are conducted to test the performance of the proposed mixed-integer formulations. Further, the proposed disjunctive  formulations are compared against state-of-the-art formulations based on the so-called shadow separation condition. Our numerical results show that the proposed disjunctive linear separation conditions outperform existing formulations in the literature and can solve significantly more instances to global optimality. For reproducibility purposes, all formulations and instances are made available on a public repository. 
\end{abstract}


\end{frontmatter}

\newpage

\section{Introduction}
Air traffic control (ATC) is an extremely dynamic and constrained environment where many decisions need to be taken in a short amount of time. Adopting automation within such environment can be vital to reduce controller workload and improve airspace capacity \citep{durand1997optimal,barnier20094d,rey2016subliminal}. Considering that traditional methods for air traffic control have been exhaustively used and are reaching their limits, automated approaches are receiving a significant and growing attention in the field \citep{vela2009two}. The International Civil Aviation Organization (ICAO) determines all regulations related to civil aviation \citep{icao2010convention}. One of its main roles is to set separation standards for commercial aviation. We focus on aircraft separation for en-route traffic. During cruise stage, separation conditions require a minimum of 5 NM horizontally and 1000 ft vertically between any pair of aircraft. A conflict between two or more aircraft is a loss of separation among these aircraft. Air traffic networks are organized in flights levels which are separated by at least 1000 ft, hence during cruise stage most conflicts occur among aircraft flying at the same flight level. Congested air traffic networks can lead to loss of separation between aircraft which impairs flight safety and may result in collisions. The aircraft conflict resolution problem (ACRP) is typically formulated as an optimization model in which the objective is to find conflict-free trajectories for a set of aircraft with intersecting trajectories. Different strategies have been used to address this problem based on the type of deconfliction manoeuvres available, namely: speed control (acceleration or deceleration), heading control, vertical control (flight level reassignment) or a combination of these manoeuvres. 

In this paper, we adopt and extend the complex number formulation for aircraft conflict resolution proposed by \citet{rey2017complex} to derive nonconvex mixed-integer programming (MIP) formulations for ACRPs and develop exact solution methods. Specifically, we propose exact optimization approaches for the two-dimensional (2D) conflict resolution problem with speed and heading control represented as continous variables; as well as for an extended problem wherein flight level (FL) changes are also available. We make the following methodological contributions to the field: i) we show that the proposed disjunctive linear pairwise aircraft separation conditions introduced by \citet{rey2017complex} are equivalent to the traditional nonlinear separation conditions; ii) we fully characterize the set of 2D conflict-free trajectories and propose propose a simple pre-processing algorithm to identify aircraft pairs which are either always conflict-free, or which cannot be separated using speed and heading control only; iii) we refine existing convex relaxations and present a novel exact constraint generation algorithm for the 2D ACRP; and, iv) we incorporate altitude control in the proposed formulations, and introduce a lexicographic optimization formulation which primary objective is to minimize the number of flight level changes before resolving outstanding conflicts via 2D velocity control. 

We conduct numerical experiments on four types of conflict resolution benchmarking instances consisting of a total of 2072 instances to test the performance of the proposed mixed-integer formulations and algorithms. We show that the proposed mixed-integer formulations and exact solution methods improve on existing exact approaches in the literature in terms of scalability. Further, we quantify the benefits of using the proposed disjunctive linear separation conditions instead of the classical so-called shadow separation conditions by implementing the latter in the proposed exact solution methods. Our experiments reveal that the disjunctive linear separation conditions significantly outperform the existing shadow separation conditions. Specifically, we report an average decrease of 86.6\% in computation time for instances solved within the available time limit, and we find that the proposed disjunctive formulation is able to solve 87.8\% of the instances compared to 69.7\% using an equivalent formulation based on the shadow separation conditions. For reproducibility purposes all formulations and instances are made available online at the public repository \small \url{https://github.com/acrp-lib/acrp-lib}.\normalsize

The paper is organized as follows. We review the state-of-the-art on aircraft conflict resolution, highlight existing research gaps and position our contributions to the field in Section \ref{lit}. We next characterize pairwise aircraft separation conditions and present nonconvex aircraft conflict resolution formulations for the 2D problem, as well as for an extension of the 2D problem with FL changes, hereby referred to as 2D+FL in Section \ref{model}. We present exact solution methods for the 2D and 2D+FL problems that build on existing convex relaxations and propose novel decomposition and constraint generation algorithms in Section \ref{methods}. Numerical results are provided in Section \ref{num} and concluding remarks are discussed in Section \ref{con}.


\section{Literature Review}\label{lit}

The first exact global optimization approach for conflict resolution problems was proposed by \cite{pallottino2002conflict} which introduced two MIP formulations: a first model was based on speed control only, and a second model was based on heading control only and assumed that all aircraft fly at the same speed. In the proposed MIP formulation for conflict resolution with speed control, the authors derived linear pairwise aircraft separation constraints based on the geometric construction introduced by \cite{bilimoria2000geometric}. These separation conditions are obtained by projecting the shadow of an intruder aircraft onto the trajectory of a reference aircraft, thus we henceforth refer to these separation conditions as the shadow separation conditions. \citet{frazzoli2001resolution} were the first to observe that this geometric construction provided a basis to characterize the set of aircraft pairwise conflict-free trajectories via linear half-planes in the relative velocity (speed and heading) plane. The authors introduced a nonconvex formulation for the conflict resolution problem with speed and heading control, and proposed a convex relaxation based on semi-definite programming as well as a heuristic algorithm to find feasible solutions on problems with up to 10 aircraft.
 
The shadow separation conditions were subsequently used in a number of formulations. \citet{alonso2011collision} proposed a mixed-integer linear programming (MILP) formulation for conflict resolution by speed and altitude control, and reported solving instances with up to 50 aircraft in competitive time. \citet{alonso2014exact} proposed a two-step formulation in which only heading control is available for deconfliction and the available angle changes are discretized. The same group of authors also proposed a nonconvex formulation involving trigonometric functions based on the shadow separation conditions that enables speed, heading and altitude control \citep{alonso2016exact}. The authors used a mixed-integer nonlinear programming (MINLP) solver to solve the resulting nonconvex formulations and reported results for the 2D ACRP with up to 7 aircraft on structured instances and up to 20 aircraft on un-structured instances. An alternative representation of pairwise aircraft separation based on conflict points has been proposed by several authors. In \cite{vela2009fuel,vela2009two,vela2009mixed}, the authors proposed several MILP formulations which aim to minimize fuel consumption, incorporate air traffic controller workload in the objective function, and account the impact of uncertainty in trajectory prediction due to wind effects. \citet{omer2015space} proposed a space-discretized MILP formulation involving a finite set of turning angles. In contrast to most other approaches, the heading control manoeuvres consist of two actions: a first heading change for collision avoidance and subsequent turn to recover the initial heading. \citet{rey2012minimization,rey2016subliminal} proposed linear upper bounds for the ACRP with speed control only and the resulting MILP formulations are able to solve realistic large-scale instances to optimality within a few seconds. 

More recently, nonlinear global optimization approaches received an increasing attention in the literature. \citet{omer2013hybridization} proposed a hybrid algorithm which uses the optimal solution of a MILP as the starting point for solving a nonlinear formulation of the same problem. \citet{cafieri2017maximizing} proposed a MINLP approach for conflict resolution with speed control only which highlights that subliminal speed control alone may not be sufficient to resolve all conflicts in dense traffic scenarios. Using a similar framework, \citet{cafieri2017mixed} presented a two-step approach where a maximum number of conflicts are first solved using speed control only and outstanding conflicts are solved by heading control. \citet{cerulli2020detecting} proposed a formulation based on bilevel optimization with multiple follower problems, each of which representing a two-aircraft separation problem. The authors presented two formulations one using speed control only and another using heading control only. A cut generation algorithm is proposed to solve the corresponding bilevel optimization problems. Recently, \citet{pelegrin2020airspace} conducted a review of the literature on exact conflict resolution approaches, and have shown that the disjunctive linear separation conditions introduced by \citet{rey2017complex} and the shadow separation conditions are equivalent. This result combined with the formal proof provided in the present manuscript linking the disjunctive separation conditions with the definitional nonlinear separation conditions (see Theorem \ref{theo:conditions}) further motivates the comparison of both (disjunctive and shadow) separation conditions in terms of computational performance which we address in our numerical experiments.


This review of the literature highlights that despite recent improvements in the development of optimization approaches for aircraft conflict resolution, there remain significant open challenges in the design of scalable and exact global optimization approaches. We next present new mixed-integer formulations for aircraft conflict resolution by speed, heading and altitude control. 

\section{Aircraft Conflict Resolution Formulations}\label{model}

In this section, we present mixed-integer formulations for aircraft conflict resolution. We first focus on the 2D problem under velocity, i.e. speed and heading, control which aims to represent a single flight level during cruise stage air traffic conditions. We first characterize 2D separation conditions for a pair of aircraft in Section \ref{2dsepa}, before introducing a compact nonconvex formulation in Section \ref{2dnonconvex}. We then extend this nonconvex formulation to the case of multiple separated flight levels (FL), denoted 2D+FL in Section \ref{3dnonconvex}.

\subsection{Characterization of 2D Separation Conditions}
\label{2dsepa}


Let $\A$ be the set of aircraft. For each $i \in \A$, $[\hat{x_i},\hat{y_i}]$ is the aircraft initial position in the 2D plane, $\hat{v}_i$ is its nominal speed (in NM/h) and $\tiz$ is its heading angle. Assuming uniform motion laws, aircraft motion can be described as: $p_i(t) = [x_i(t),y_i(t)]$, where $x_i(t) = \hat{x_i} + \qi \hat{v}_i\cos(\tiz + \ti) t$ and $y_i(t) = \hat{y_i} + \qi \hat{v}_i\sin(\tiz + \ti) t$. In this model, the decision variables are  $\qi$, which is the speed control variable that determines the acceleration or deceleration with regards to the nominal speed $\hat{v_i}$ ($\qi$ equals to 1 means no speed variation) and $\ti$, which is the heading control variable that determines the deviation with regards to the nominal trajectory ($\ti$ equal to 0 means no deviation in heading angle). 

Let $\P = \{i,j \in A : i < j\}$ be the set of aircraft pairs, the relative motion of $(i,j) \in \P$ is denoted $p_{ij}(t) = p_i(t) - p_j(t)$. Let $d$ be the minimum separation distance (e.g. 5 NM). We next define the notion of 2D separation for a pair of aircraft.

\begin{definition}[2D separated trajectories]
The trajectories of a pair of aircraft $(i,j) \in \P$ are said to be 2D separated, i.e. conflict-free, if and only if:
\begin{equation}\label{eq:sepcond}
||p_{ij}(t)|| \geq d, \qquad \forall t \geq 0.
\end{equation}
\end{definition}

Let $v_{ij} = v_i - v_j$ be the 2D relative velocity vector of $(i,j) \in \P$, i.e. $v_{ij} = [\vijx,\vijy]^\top$ with:
\begin{subequations}\label{eq:v}
\begin{align}
& \vijx = \qi \viz\cos(\tiz + \ti) - \qj \vjz\cos(\tjz + \tj), \\
& \vijy = \qi \viz\sin(\tiz + \ti) - \qj \vjz\sin(\tjz + \tj).
\end{align}
\end{subequations}

Aircraft relative velocity equations are linear with regards to speed control variables $\qi$ and $\qj$, but  nonlinear with regards to heading control variables $\ti$ and $\tj$. Expanding the expression in Eq. \eqref{eq:sepcond} and denoting the relative initial position of each pair $i,j \in \A: i < j$, $\hat{p}_{ij}$, we obtain a second-order polynomial function:
\begin{equation}\label{eq:sep1}
f_{ij}(t) \equiv ||v_{ij}||^2 t^2 + 2||\hat{p}_{ij}||\cdot v_{ij}t +||\hat{p}_{ij}||^2 - d^2 \geq 0,
\end{equation}

which is a minimum at:
\begin{equation}\label{eq:tmequation}
\tm(\vijx,\vijy) = \frac{-\hat{p}_{ij} \cdot v_{ij}}{||v_{ij}||^2}.
\end{equation}

Evaluating Eq. \eqref{eq:sep1} at $\tm(\vijx,\vijy)$  yields a time-independent separation condition \citep{rey2017complex,cafieri2017maximizing,cafieri2017mixed}. After multiplying by $||v_{ij}||^2$ on both sides, we obtain:
\begin{equation}\label{gfunction}
g_{ij}(\vijx,\vijy) = \vijx^2(\hat{y}_{ij}^2 - d^2) + \vijy^2(\hat{x}_{ij}^2  - d^2) -  \vijx\vijy(2\hat{x}_{ij}\hat{y}_{ij})\geq 0.
\end{equation}

Assuming aircraft are initially separated, if $\tm(\vijx,\vijy) \leq 0$, then they are diverging and do not incur any risk of future conflict. If $\tm(\vijx,\vijy) \geq 0$ and $g(\vijx,\vijy) \geq 0$, aircraft are converging but separation is ensured. Otherwise, if $\tm(\vijx,\vijy) \geq 0$ and $g(\vijx,\vijy) \leq 0$, there is a loss of separation and aircraft trajectories should be adjusted to avoid it. Hence, pairwise aircraft separation conditions for $(i,j) \in \P$ can be written as:
\begin{equation}\label{eq:sepconditions}
g_{ij}(\vijx,\vijy) \geq 0 \vee \tm(\vijx,\vijy) \leq 0.
\end{equation}

To linearize the separation condition \eqref{gfunction} with regards to variables $\vijx$ and $\vijy$, we adopt the approach proposed by \cite{rey2017complex} and recall the main steps hereafter. Observe that the solutions of the equation $g(\vijx,\vijy) = 0$ can be identified by alternatively fixing variable $\vijx$ and $\vijy$, and calculating the roots of the resulting single-variable quadratic equations. Isolating each variable, we obtain the discriminants:

\begin{equation}\label{discriminant}
\begin{cases}
\Delta_{\vijx} = 4d^2\vijy^2(\hat{x}^2_{ij} + \hat{y}^2_{ij} - d^2),\\
\Delta_{\vijy} = 4d^2\vijx^2(\hat{x}^2_{ij} + \hat{y}^2_{ij} - d^2).
\end{cases}
\end{equation}

Assuming aircraft are initially separated, then $\hat{x}^2_{ij} + \hat{y}^2_{ij} - d^2 \geq 0$ and the discriminants are positive, and the roots of equation $g(\vijx,\vijy) = 0$ are the lines defined by the system of equations:
\begin{subequations}
\begin{align}
(\hat{y}_{ij}^2 -d^2)\vijx - (\hat{x}_{ij}\hat{y}_{ij} + d\sqrt{\hat{x}_{ij}^2 + \hat{y}_{ij}^2 - d^2})\vijy = 0, \label{l1}\\
(\hat{y}_{ij}^2 -d^2)\vijx - (\hat{x}_{ij}\hat{y}_{ij} - d\sqrt{\hat{x}_{ij}^2 + \hat{y}_{ij}^2 - d^2})\vijy = 0, \label{l2}\\ 
(\hat{x}_{ij}^2 -d^2)\vijy - (\hat{x}_{ij}\hat{y}_{ij} + d\sqrt{\hat{x}_{ij}^2 + \hat{y}_{ij}^2 - d^2})\vijx = 0, \label{l3}\\
(\hat{x}_{ij}^2 -d^2)\vijy - (\hat{x}_{ij}\hat{y}_{ij} - d\sqrt{\hat{x}_{ij}^2 + \hat{y}_{ij}^2 - d^2})\vijx = 0. \label{l4}
\end{align}
\label{planesequations}
\end{subequations}

Let us emphasize that if all coefficients in \eqref{l1}-\eqref{l4} are non-zero, then \eqref{l1} is identical to \eqref{l3} and \eqref{l2} is identical to \eqref{l4}. Observe that	
\begin{align*}
&&\pzx \pzy \pm d \sqrt{\pzx^2 + \pzy^2 - d^2} &= 0, \\
\Rightarrow  && d^2 (\pzx^2 + \pzy^2 - d^2) &= \pzx^2 \pzy^2,\\
\Leftrightarrow  && \pzx^2 \pzy^2 - d^2 (\pzx^2 + \pzy^2 - d^2) &= 0,\\
\Leftrightarrow  && (\pzx^2 - d^2) (\pzy^2 - d^2) &= 0.
\end{align*}

Eqs. \eqref{l1}, \eqref{l2}, \eqref{l3} and \eqref{l4} define two lines, denoted $R_1$ and $R_2$, in the plane $\{(\vijx,\vijy) \in \mathbb{R}^2\}$ and the sign of $g(\vijx,\vijy)$ can be characterized based on the position of $(\vijx,\vijy)$ relative to these lines. Recall that according to Eq. \eqref{eq:tmequation}, the sign of the dot product~$\hat{p}_{ij} \cdot v_{ij}$ indicates aircraft convergence or divergence. Let \eqref{plane} be the equation of the line corresponding to $\hat{p}_{ij} \cdot v_{ij}$:
\begin{equation}\label{plane}\tag{$P$}
\vijx\hat{x}_{ij} + \vijy\hat{y}_{ij} = 0.
\end{equation}

\begin{figure}[t]
\centering
\subfloat[The hashed pink region represents $g(\vijx,\vijy) \geq 0$. The hashed green half-plane represents diverging trajectories, i.e. $\tm(\vijx,\vijy) \leq 0$.]{\includegraphics[width=0.47\linewidth]{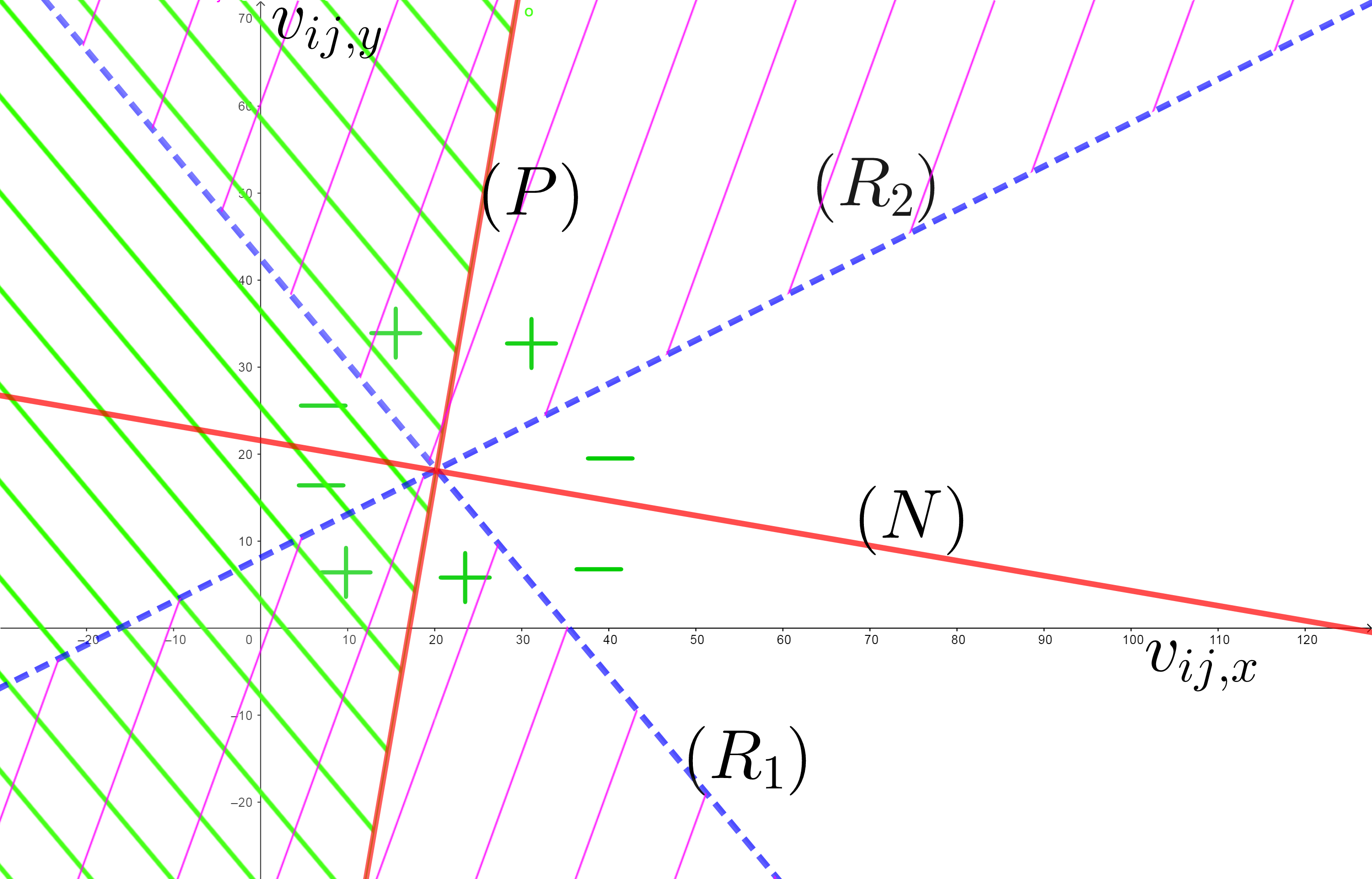}\label{fig:gznp}}
\hfil
\subfloat[Illustration of the disjunctive convex regions: $\z = 1$ correspond to the region hashed in yellow and $\z = 0$ corresponds to the region hashed in blue.]{\includegraphics[width=0.45\linewidth]{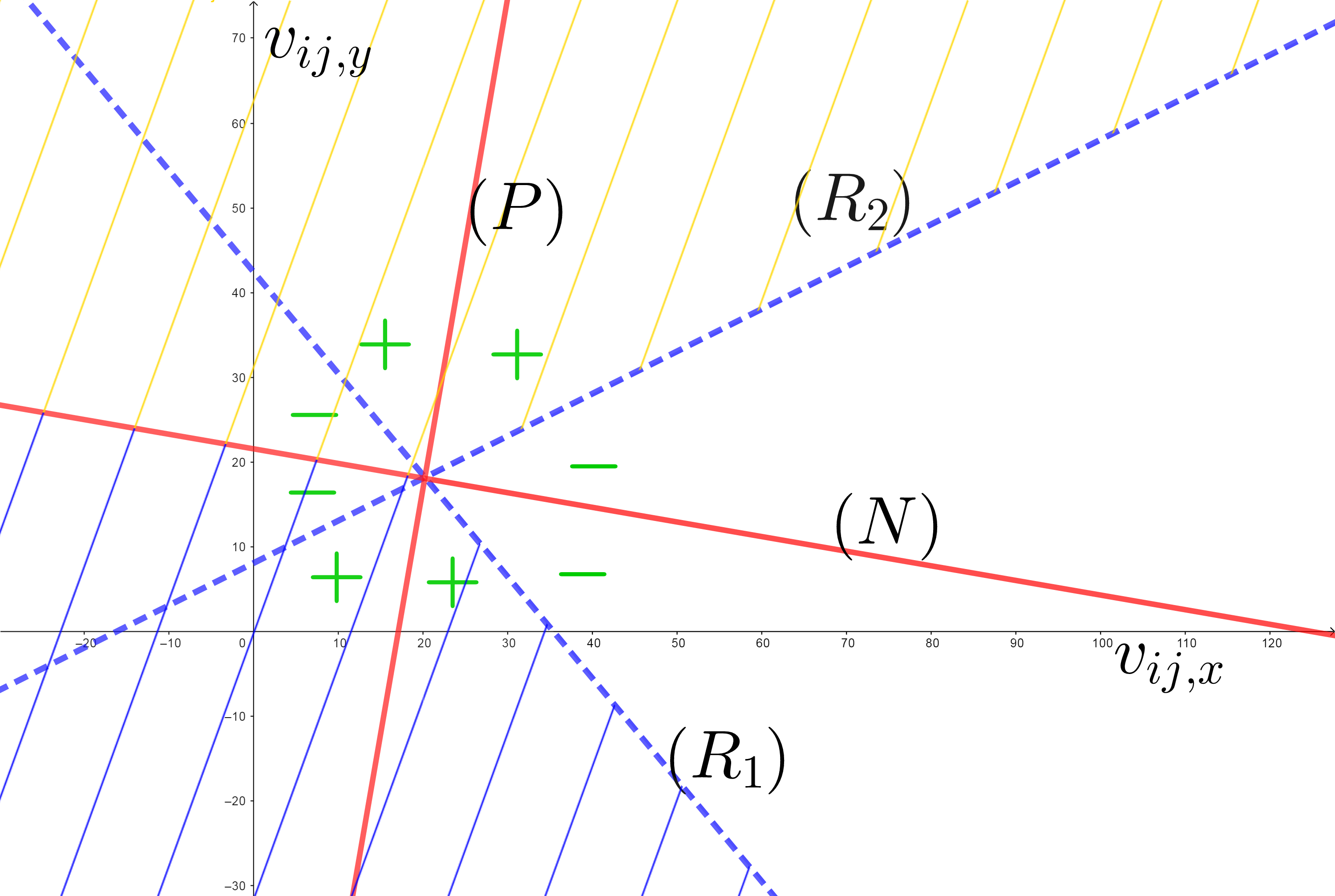}\label{fig:qznp}}
\caption{Illustration of a two-aircraft conflict in the plane $\{(\vijx,\vijy) \in \mathbb{R}^2\}$. The red lines represent the lines $P$ and $N$. The dashed blue lines correspond to the linear equations $R_1$ and $R_2$ that are the roots of $g(\vijx,\vijy) = 0$. The sign of $g(\vijx,\vijy)$ is shown by the + and - green symbols.}
\end{figure}

The line defined by \eqref{plane} splits the plane $\{(\vijx,\vijy) \in \mathbb{R}^2\}$ in two half-planes, each of which representing converging and diverging trajectories, respectively. This is illustrated in Figure \ref{fig:gznp} which depicts a two-aircraft conflict in the plane $\{(\vijx,\vijy) \in \mathbb{R}^2\}$. The sign of $g(\vijx,\vijy)$ is shown by the + and - green symbols and the hashed pink region corresponds to $ g(\vijx,\vijy) \geq 0$. The hashed green half-plane delimited by \eqref{plane} represents diverging trajectories, i.e. $\tm(\vijx,\vijy) \leq 0$.

Consider the line normal to \eqref{plane}, denoted \eqref{normalplane}:
\begin{equation}\label{normalplane}\tag{$N$}
\vijy\hat{x}_{ij} - \vijx\hat{y}_{ij} = 0.
\end{equation}

Recall that any point $(\vijx,\vijy)$ such that $\tm \leq 0$ or $g(\vijx,\vijy) \geq 0$ corresponds to a pair of conflict-free trajectories. Hence, the conflict-free region is nonconvex and is represented by the union of the green and pink hashed regions in Figure \ref{fig:gznp}. Reciprocally, the conflict region, corresponding to conflicting trajectories is convex and represented by the non-hashed region in Figure \ref{fig:gznp}. An equivalent expression of Eq. \eqref{planesequations} was proposed by \citet{frazzoli2001resolution} which observed that the set of conflict-free trajectories could be characterized by the union of two half-planes. We next show through Lemmas \ref{lem:normal} and \ref{lem:negative} that \eqref{normalplane} is a bisector of the angle formed by lines $R_1$ and $R_2$ in the conflict zone ($g(\vijx,\vijy) \leq 0$) and can be used to generate two disjunctive but convex conflict-free regions.


\begin{lemma}
\label{lem:normal}
The lines \eqref{plane} and \eqref{normalplane} are bisectors of the angles formed by the two lines $R_1$ and $R_2$ representing the solutions of $g(\vijx,\vijy)= 0$.
\end{lemma}

\begin{lemma}
\label{lem:negative}
$g(\vijx,\vijy) \leq 0$ for all points $(\vijx,\vijy)$ of the normal line \eqref{normalplane}.
\end{lemma}

\newpage

Lemmas \ref{lem:normal} and \ref{lem:negative} (the proofs are provided in \ref{pl1} and \ref{pl2}) assert that \eqref{normalplane} can be used to split the conflict-free region into two convex but disjunctive regions.  We model this disjunction using variable $z_{ij} \in \{0,1\}$ defined as:
\begin{subequations}
\label{eq:zd}
\begin{align}
\vijy\hat{x}_{ij} - \vijx\hat{y}_{ij} \leq 0, &\quad \text{ if } \z=1, \quad \forall (i,j) \in \P,\\
\vijy\hat{x}_{ij} - \vijx\hat{y}_{ij} \geq 0, &\quad \text{ if } \z=0, \quad \forall (i,j) \in \P.
\end{align}
\end{subequations}

In each convex sub-region, the lines defined by \eqref{l1}-\eqref{l4} delineate the conflict-free region. The expressions of these lines depends on aircraft initial positions, i.e. $\hat{x}_{ij}$, $\hat{y}_{ij}$. Recall that we denote $R_1$ and $R_2$ the equation of these lines. Integer-linear separation conditions with regards to aircraft velocity components can be derived as follows:
\begin{subequations}
\label{eq:rs}
\begin{align}
\vijy \gamma_{ij}^l - \vijx \phi_{ij}^l \leq 0, &\quad \text{ if } \z=1, \quad \forall (i,j) \in \P, \label{eq:r1}\\
\vijy \gamma_{ij}^u - \vijx \phi_{ij}^u \geq 0, &\quad \text{ if } \z=0, \quad \forall (i,j) \in \P, \label{eq:r2}
\end{align}
\end{subequations}

where $\gamma_{ij}^l$, $\phi_{ij}^l$ and $\gamma_{ij}^u$, $\phi_{ij}^u$ are coefficients of the lines \eqref{l1}-\eqref{l4} corresponding to the roots of $g(\vijx,\vijy) = 0$. The proposed linear disjunction is illustrated in Figure \ref{fig:qznp} which depicts the resulting convex sub-regions corresponding the disjunction $\z \in \{0,1\}$ for a two-aircraft conflict. This leads to the following result.

\begin{theorem}\label{theo:conditions}
The disjunctive linear separation conditions \eqref{eq:zd}-\eqref{eq:rs} fully characterize the set of aircraft pairwise conflict-free trajectories as given by Eq. \eqref{eq:sepcond}
\end{theorem}

Theorem \ref{theo:conditions} (the proof is provided in \ref{pt}), asserts that the  disjunctive separation conditions \eqref{eq:zd}-\eqref{eq:rs} equivalent to the definitional nonlinear separation conditions. Further, these disjunctive separation conditions are linear with regards to aircraft velocity variables $\vijx$ and $\vijy$, and only require a single binary variable per pair of aircraft. This is expected to improve on the so-called shadow separation conditions which are also linear with regards to aircraft velocity variables, but require four binary variable per pair of aircraft \citep{pallottino2002conflict,alonso2011collision,alonso2016exact}. 

\newpage

To characterize the set of 2D conflict-free trajectories, we examine the relative velocity vector $v_{ij}$ as a function of trajectory control bounds. For each aircraft $i \in \A$, we assume that the speed rate variable is lower bounded by $\qmin$ and upper bounded by $\qmax$, i.e.:
\begin{equation}\label{eq:qbound}
\qmin \leq \qi \leq \qmax, \qquad \forall i \in \A.
\end{equation}

We assume that the heading deviation is lower bounded by $\underline{\ti}$ and upper bounded by $\overline{\ti}$, i.e.:
\begin{equation}\label{eq:tbound}
\underline{\ti} \leq \ti \leq \overline{\ti},  \qquad \forall i \in \A.
\end{equation}

To derive lower and upper bounds on relative velocity components $\vijx$ and $\vijy$, we re-arrange Eq. \eqref{eq:v} using trigonometric identities:
\begin{subequations}\label{eq:v2}
\begin{align}
& \vijx = \qi \viz\cos(\tiz)\cos(\ti) - \qi \viz\sin(\tiz)\sin(\ti) - \qj \vjz\cos(\tjz)\cos(\tj) + \qj \vjz\sin(\tjz)\sin(\tj), \\
& \vijy = \qi \viz\sin(\tiz)\cos(\ti) + \qi \viz\cos(\tiz)\sin(\ti) - \qj \vjz\sin(\tjz)\cos(\tj) - \qj \vjz\cos(\tjz)\sin(\tj).
\end{align}
\end{subequations}

Let $\underline{v}_{ij,x},\overline{v}_{ij,x}$ and $\underline{v}_{ij,y},\overline{v}_{ij,y}$ be the lower and upper bounds for $\vijx$ and $\vijy$, respectively. These bounds can be determined using Eq. \eqref{eq:v2} and the bounds on speed and heading control provided in Eqs. \eqref{eq:qbound} and \eqref{eq:tbound}. The derived bounds on the relative velocity components can be used to define a box in the plane $\{(\vijx,\vijy) \in \mathbb{R}^2\}$.

\begin{definition}[Relative velocity box]
\label{box}
Consider a pair of aircraft $(i,j) \in \P$. Let $\B$ be the  subset of $\R^2$ defined as
\begin{equation}
\B = \left\{(\vijx,\vijy) \in \mathbb{R}^2 : \underline{v}_{ij,x}\leq \vijx \leq \overline{v}_{ij,x}, \underline{v}_{ij,y} \leq \vijy \leq \overline{v}_{ij,y}\right\}.
\end{equation}
$\B$ is the relative velocity box of $(i,j) \in \P$.
\end{definition}

The relative velocity box $\B$ characterize all possible trajectories for the pair $(i,j) \in \P$ based on the available 2D deconfliction resources, i.e. speed and heading controls. To characterize the set of conflict-free trajectories of a pair of aircraft $(i,j) \in \P$, we compare the relative position of the relative velocity box $\B$ with the conflict region of this pair of aircraft. Observe that the conflict region is convex and can be defined based by reversing the inequalities \eqref{eq:rs} and omitting the disjunction $\z \in \{0,1\}$.

\begin{definition}[Conflict region]
\label{cr}
Consider a pair of aircraft $(i,j) \in \P$. Let $\C$ be the subset of $\R^2$ defined as 
\begin{equation}
\C = \left\{(\vijx,\vijy) \in \R^2 : \vijy \gamma_{ij}^l - \vijx \phi_{ij}^l \geq 0 \land \vijy \gamma_{ij}^u - \vijx \phi_{ij}^u \leq 0 \right\}.
\end{equation}
$\C$ is the conflict region of $(i,j) \in \P$.
\end{definition}

\begin{figure}
	\centering
	{\includegraphics[width=0.3\textwidth]{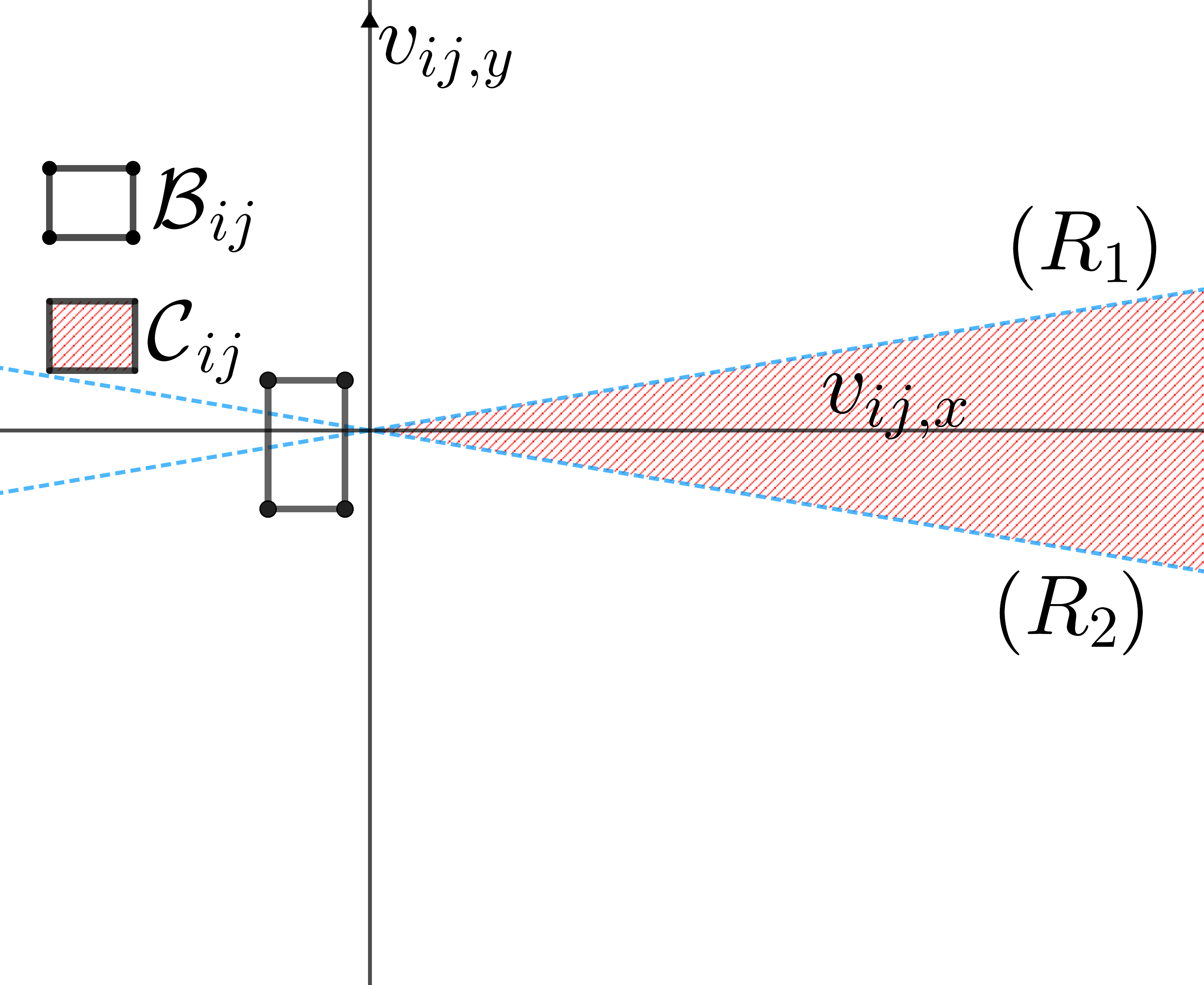}} 
	{\includegraphics[width=0.3\textwidth]{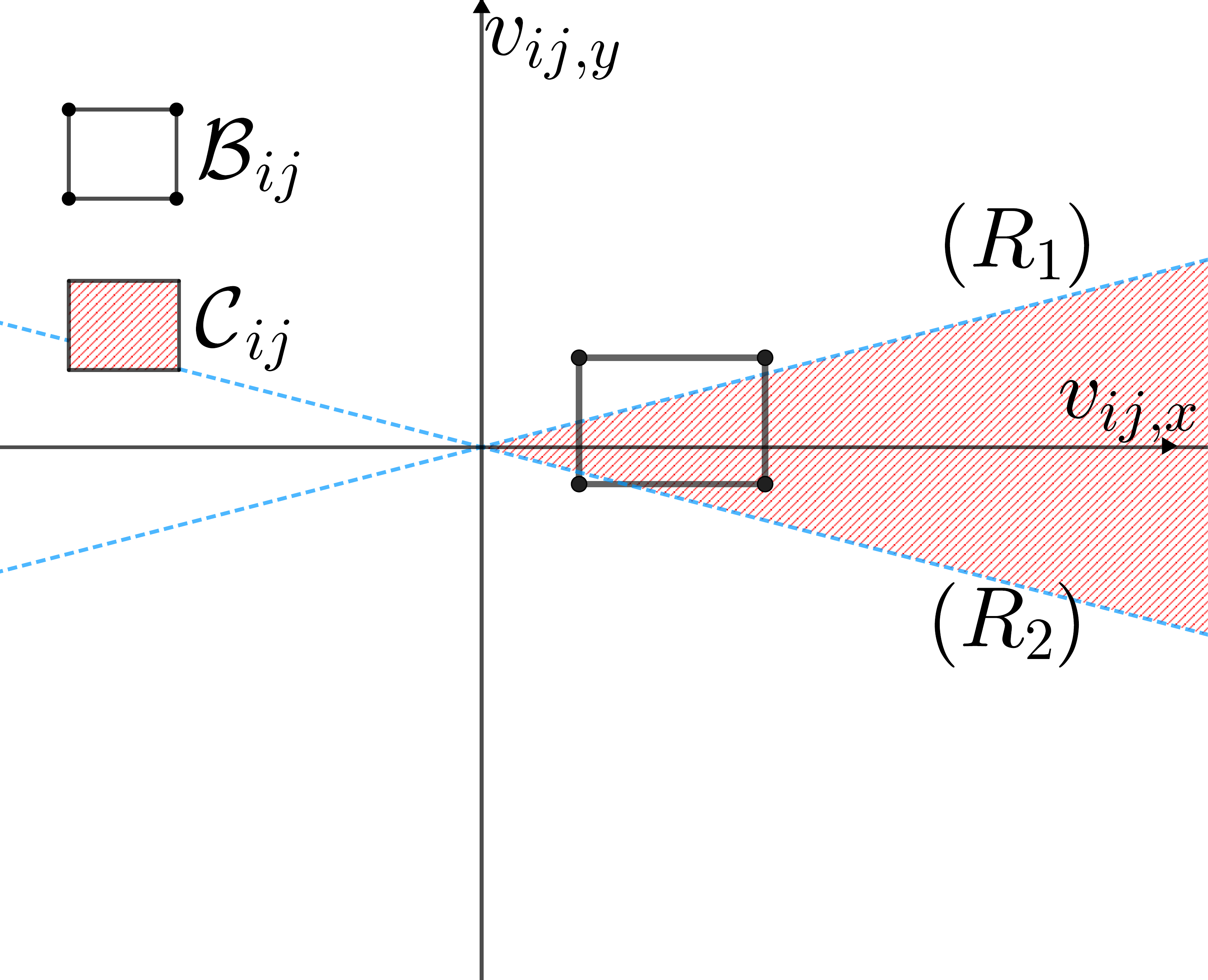}}
	{\includegraphics[width=0.3\textwidth]{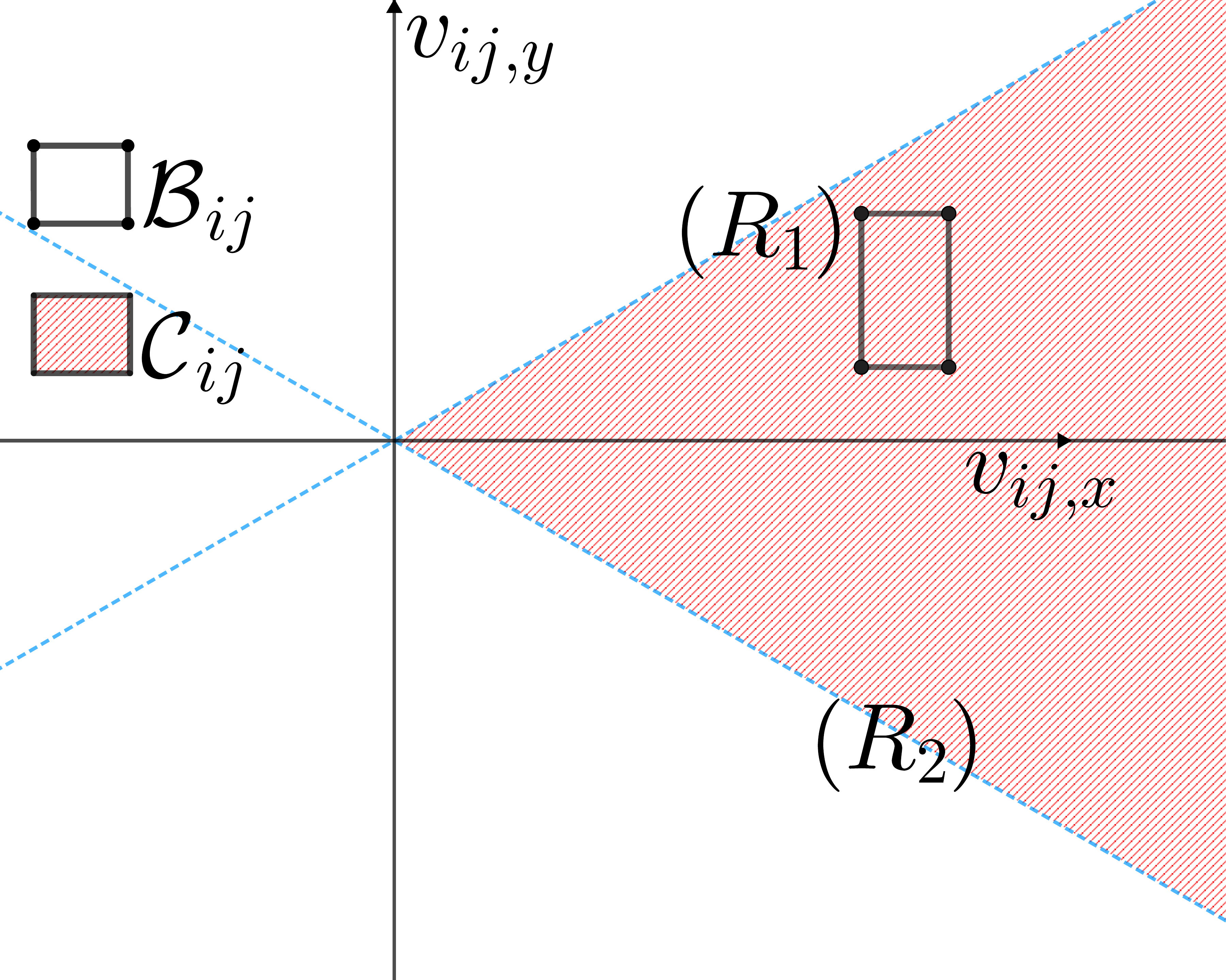}}\\
	\subfloat[Conflict-free aircraft pair, i.e. $\B \cap \C = \emptyset$.  Nominal headings are $\hat{\theta}_i$ = 2.01 and $\hat{\theta}_j$ = 1.30.] {\includegraphics[width=0.3\textwidth]{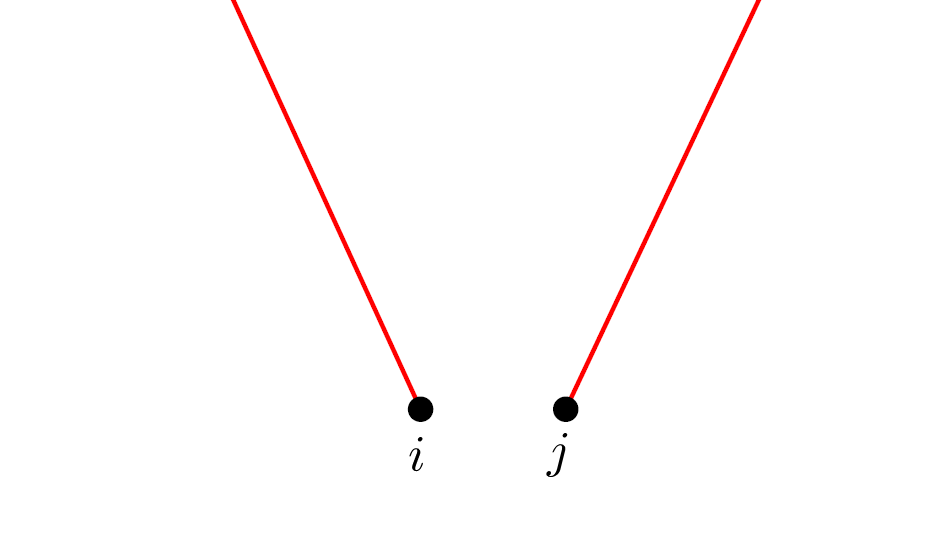}\label{alwayssep}} 
	\qquad 
	\subfloat[Separable aircraft pair, i.e. $\B \cap \C \neq \emptyset$ and $\B \not\subset \C$. Nominal headings are $\hat{\theta}_i$ = 1.25 and $\hat{\theta}_j$ = 1.88.] {\includegraphics[width=0.3\textwidth]{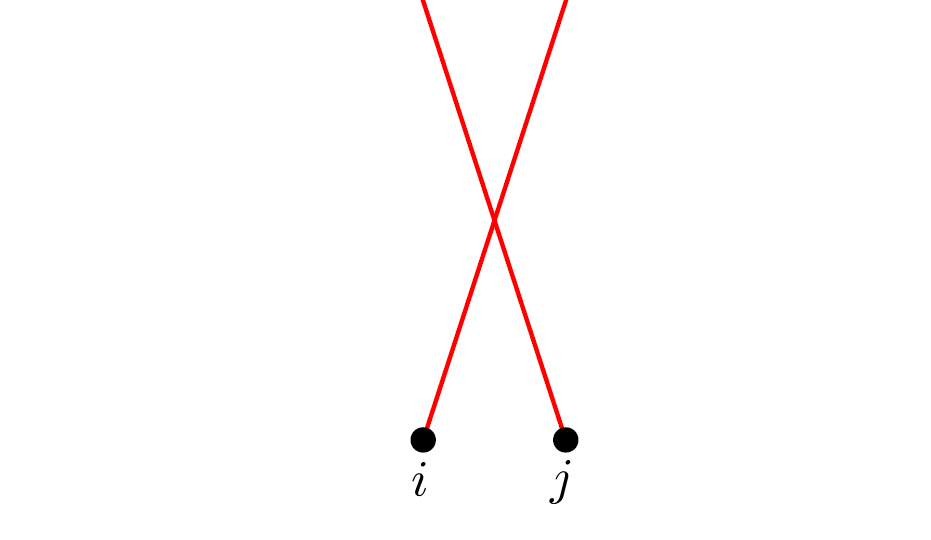}\label{P}}
	\qquad 
	\subfloat[Non-separable aircraft pair, i.e. $\B \subset \C$. Nominal headings are $\hat{\theta}_i$ = 1.04 and $\hat{\theta}_j$ = 2.09.] {\includegraphics[width=0.3\textwidth]{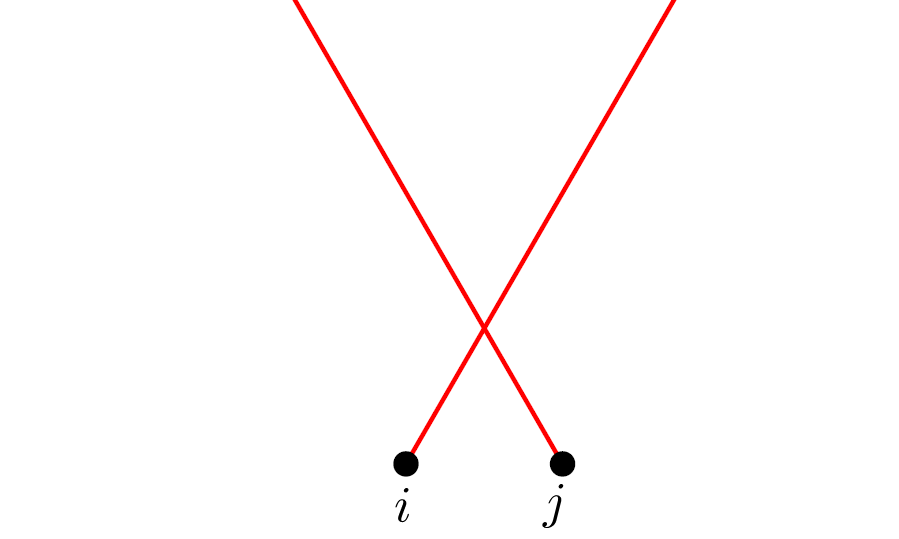}\label{infeas}}
	\caption{Three aircraft configurations illustrating conflict-free (\ref{alwayssep}), separable (\ref{P}) and non-separable (\ref{infeas}) pairs. Aircraft initial positions and trajectories are depicted at the bottom part of each sub-figure, where black dots correspond to initial positions and red lines represent initial trajectories. The top part of each sub-figure shows a graph of the relative velocity box $\B$ represented by a black rectangle, and the conflict region $\C$ represented by a red hashed area. Dashed blue lines represent the solutions of $g(\vijx,\vijy)=0$ ($R_1$ and $R_2$). In all cases, aircraft $i$ and $j$ have nominal speeds $\hat{v}_i = \hat{v}_j = 500$ NM/h and are initially separated by 30 NM in the x-axis direction.}
	\label{fig:allcases}
\end{figure}

The conflict region of pair of aircraft represents the set of relative velocity vectors which correspond to a conflict. The relative positions of the relative velocity box $\B$ and the conflict region $\C$ can be examined to determine the existence or not of a potential conflict. For any pair $(i,j) \in \P$, if $\B \cap \C = \emptyset$, then aircraft $i$ and $j$ are separated for any combination of controls; conversely if $\B \subset \C$ then $i$ and $j$ cannot be separated via speed or heading control within the assumed control bounds; otherwise, $\B$ and $\C$ intersect but do not completely overlap. This is illustrated in Figure \ref{fig:allcases} which illustrates the three possible cases. Figure \ref{alwayssep} illustrates the case where aircraft $i$ and $j$ are separated for any combination of speed and heading control---we say that such pairs are conflict-free. Figure \ref{P} depicts the case where $\B$ and $\C$ only partially intersect---we say that such pairs are separable. Last, Figure \ref{infeas} illustrates the case where $\B \subset \C$---we say that such pairs are non-separable. The following propositions provide methods to efficiently determine if, given controls bounds on speed and heading, a pair of aircraft is either conflict-free or non-separable.

\begin{prop}[Conflict-free aircraft pair]
Consider a pair of aircraft $(i,j) \in \P$, and let $LP(i,j)$ be the feasibility linear program defined as:
\begin{equation*}
LP(i,j): 
\begin{cases}
&\emph{Minimize } 1,  \\
&\emph{Subject to:}  \\
& \vijy \gamma_{ij}^l - \vijx \phi_{ij}^l \geq 0,\\
& \vijy \gamma_{ij}^u - \vijx \phi_{ij}^u \leq 0,\\
& (\vijx,\vijy) \in \B.
\end{cases}
\end{equation*}

The pair $(i,j)$ is conflict-free for any 2D control if and only if $LP(i,j)$ is infeasible. 
\label{prop:cf}
\end{prop}

\begin{prop}[Non-separable aircraft pair]
Consider a pair of aircraft $(i,j) \in \P$. The pair $(i,j)$ is non-separable if and only if each of the four extreme points of $\B$ corresponds to a conflict.
\label{prop:ns}
\end{prop}

Using Propositions \ref{prop:cf} and \ref{prop:ns} (the proofs are provided in  \ref{pp1} and \ref{pp2}), we can design an efficient pre-processing algorithm to partition the set of aircraft pairs $\P$ of a 2D ACRP instance into three categories: conflict-free pairs denoted $\PF$, separable pairs denoted $\PS$ and non-separable pairs $\PI$.\\

\begin{algorithm}[H]
\KwIn{$\A$, $\hat{\bm{\theta}}$, $\hat{\bm{v}}$, $\qlb$, $\qub$, $\tlb$, $\tub$}
\KwOut{$\P$, $\PF$, $\PS$, $\PI$}
$\P \gets \{i,j \in \A : i <j\}$ \\
$\PF, \PS, \PI \gets \emptyset$ \\
\For{$(i,j) \in \P$}{
Solve $LP(i,j)$\\
\If{$LP(i,j)$ \emph{ is infeasible}}{
    $\PF \gets \PF \cup \{(i,j)\}$\\
}
\Else{
$k \gets 0$ \\
\For{$(\vijx,\vijy) \in E(\B)$}{
    \If{$g_{ij}(\vijx,\vijy) < 0 \land \tm(\vijx,\vijy) > 0 $}{
    $k \gets k + 1$\\
    }
}
\If{$k = 4$}{
    $\PI \gets \PI \cup \{(i,j)\}$\\
}
}
\Else{
    $\PS \gets \PS \cup \{(i,j)\}$\\
}
}
\caption{Pre-processing of aircraft pairs}
\label{algo:prepro}
\end{algorithm}

\newpage
To identify conflict-free pairs, $LP(i,j)$ is solved and a pair of aircraft is conflict-free if and only if the LP is infeasible. Observe that the feasibility linear program $LP(i,j)$ can be solved by enumerating all extreme points of its feasible region and tests if this corresponds to a conflict or not. Since the $LP(i,j)$ consists of four bound constraints and two shared constraints, there is a total of 13 extreme points to test (the combinations of the  bound constraints of a variable can be excluded). To identify non-separable pairs, we denote $E(\B)$ the set of extreme points of the relative velocity box $\B$ for any pair $(i,j) \in \P$ and use the separation condition \eqref{eq:sep1} to determine if all extreme points are conflicts or not. This procedure is summarized in Algorithm \ref{algo:prepro} (we use boldface to denote vectors). Observe that pairwise variables and constraints need only to be indexed by the set of separable pairs $\PS$ since pairs in $\PF$ are always conflict-free. Further, any 2D conflict resolution problem such that $|\PI|>0$ is infeasible.

\subsection{Complex Formulation for the 2D Aircraft Conflict Resolution Problem}
\label{2dnonconvex}

An alternative way to represent aircraft motion equation is via the complex number formulation introduced by \cite{rey2017complex}. Let $V_i$ be the complex number defined as:
\begin{equation}
V_i = \qi(\cos(\ti) + \ii\sin(\ti)), \qquad \forall i \in \A.
\end{equation}

Representing $V_i$ in its polar form with $\dx = \qi\cos(\ti)$ and $\dy = \qi\sin(\ti)$, yield:
\begin{equation}
V_i = \dx + \ii\dy, \qquad \forall i \in \A.
\end{equation}

The magnitude of $V_i$ is $|V_i| = \sqrt{\dx^2 + \dy^2} = \qi$ and its argument is $\arg (V_i) = \arctan\frac{\dx}{\dy} = \ti$. This approach is inspired by complex number formulation for the optimal power flow problem in power systems, such as in \cite{hijazi2017convex} and \cite{coffrin2015qc}. Accordingly, the relative motion equations of a pair of aircraft can be rewritten as:
\begin{subequations}\label{eq:speedxdy}
\begin{align}
&\vijx =  \dx\hat{v}_i\cos(\hat{\ti}) - \dy\hat{v}_i\sin(\hat{\ti}) - \delta_{j,x}\hat{v}_j\cos(\hat{\ti}) + \delta_{j,y}\hat{v}_j\sin(\hat{\tj}), &&\forall (i,j) \in \PS, \\
&\vijy =  \dy\hat{v}_i\cos(\hat{\ti}) - \dx\hat{v}_i\sin(\hat{\ti}) - \delta_{j,y}\hat{v}_j\cos(\hat{\ti}) + \delta_{j,x}\hat{v}_j\sin(\hat{\tj}), &&\forall (i,j) \in \PS.
\end{align}
\end{subequations}

The variables $\dx$ and $\dy$ are used as the main control variables in this formulation and their bounds are derived from the bounds of the original control variables $\qi$ and $\ti$:
\begin{subequations}\label{eq:boundsdxdy}
\begin{align}
& \qmin\cos(\max\{|\overline{\ti}|,|\underline{\ti}|\}) \leq \dx \leq \qmax, && \forall i \in \A, \label{eq:bounddx}\\
& \qmax\sin(\underline{\ti}) \leq \dy \leq \qmax\sin(\overline{\ti}), && \forall i \in \A. \label{eq:bounddy}
\end{align}
\end{subequations}

The speed control constraint \eqref{eq:qbound} can be reformulated in quadratic form as:
\begin{subequations}\label{eq:boundq}
\begin{align}
& \qmin^2 \leq \dx^2 + \dy^2, && \forall i \in \A,\label{eq:lboundq}\\
& \qmax^2 \geq \dx^2 + \dy^2, && \forall i \in \A.\label{eq:uboundq}
\end{align}
\end{subequations}

The heading control constraint \eqref{eq:tbound} can be reformulated in linear form as:
\begin{subequations}\label{eq:boundtheta}
\begin{align}
& \dy \geq \dx\tan(\underline{\ti}), && \forall i \in \A,\\
& \dy \leq \dx\tan(\overline{\ti}), && \forall i \in \A.
\end{align}
\end{subequations}

To design the objective function, we introduce a preference weight $w \in \ ]0,1[$ to balance the trade-offs among velocity controls, i.e. speed and heading. We extend the objective function proposed by \citet{rey2017complex} as follows:
\begin{equation}\label{eq:obj2d}
\text{minimize} \sum_{i \in \A} w\dy^2 + (1 - w)(1 - \dx)^2.
\end{equation}

We next show that the objective function \eqref{eq:obj2d} attains its minimum value when aircraft have deviation-free trajectories, i.e. $\qi=1$ and $\ti=0$ for all $i\in \A$.

\begin{prop}
The 2D objective function \eqref{eq:obj2d} is convex with regards to decision variables $\dx$ and $\dy$ for any value $w \in \ ]0,1[$, and is minimal for $\qi=1$ and $\ti=0$ for all aircraft $i \in \A$.
\label{prop:conv}
\end{prop}

Proposition \ref{prop:conv} (the proof is provided in  \ref{pp3}) shows that the proposed 2D objective function achieves an optimal value for deviation-free aircraft trajectories for any preference weight $w \in \ ]0,1[$. We will show in Section \ref{sa} that increasing $w$ increases the penalization of heading deviations whereas decreasing $w$ increases the penalization of speed deviations. The 2D nonconvex formulation for the ACRP is summarized in Model \ref{mod:2dnonconvex}. 

\newpage



\begin{model}[Nonconvex 2D Formulation]
\label{mod:2dnonconvex}
\begin{subequations}
\begin{align}
&\emph{Minimize Total 2D Deviation} \quad \eqref{eq:obj2d}, \nonumber \\
&\emph{Subject to:}  && \nonumber \\
&\emph{Motion Equations} \quad \eqref{eq:speedxdy}, \nonumber\\
&\emph{Separation Conditions} \quad \eqref{eq:zd}, \eqref{eq:rs}, \nonumber \\
&\emph{Speed and Heading Control Constraints} \quad  \eqref{eq:boundq}, \eqref{eq:boundtheta}, \nonumber\\
&\emph{Variable Bounds} \quad \eqref{eq:boundsdxdy}, \nonumber \\
& \vijx,\vijy \in \B, &&\forall (i,j) \in \PS, \\
& \z \in \{0,1\}, &&\forall (i,j) \in \PS, \\
& \dx, \dy \in \R, && \forall i \in \A. 
\end{align}
\end{subequations}
\end{model}

Model \ref{mod:2dnonconvex} provides a compact formulation for the 2D ACRP with continuous speed and heading control variables which requires a single binary variable per pair of aircraft. This formulation is nonconvex due to the speed lower bound constraint \eqref{eq:lboundq} which is nonconvex quadratic. Note that we use the convex hull reformulation of the linear On/Off constraints \eqref{eq:zd} and \eqref{eq:rs} as derived in \cite{Hijazi14}. This reformulation does not require the introduction of auxiliary variables and is proved to provide the tightest continuous relaxation for each On/Off constraint. Coefficients $\gamma_{ij}^l$, $\phi_{ij}^l$ and $\gamma_{ij}^u$, $\phi_{ij}^u$ (present in \eqref{eq:zd} and \eqref{eq:rs}) can be pre-processed based on the sign of $\hat{x}_{ij}$ and $\hat{y}_{ij}$. For implementation details, a fully reproducible formulation can be found at: \small \url{https://github.com/acrp-lib/acrp-lib}.\normalsize 

\subsection{Complex Formulation Formulation for the 2D+FL Aircraft Conflict Resolution Problem}
\label{3dnonconvex}

To model FL changes, we assume that each aircraft $i \in \A$ is initially assigned to a base FL denoted $\hat{\rho}_i$. We assume that adjacent FLs are vertically separated (e.g. by 1000 ft). Thus, we only need to impose separation constraints on pairs of aircraft which share the same FL. Let $\Zi$ denote the set of available FLs for each aircraft $i \in \A$, and consider the binary variable $\rik$ defined as:
\begin{equation}\label{eq:flight}
\rik = \begin{cases}
1 & \text{ if aircraft $i \in \A$ is on FL $k \in \Zi$},\\
0 & \text{ otherwise}.
\end{cases}
\end{equation}

By design, $\hat{\rho}_i \in \Zi$ and we require that each aircraft $i \in \A$ be assigned to exactly one FL in its reachable set $\Zi$ via constraint \eqref{eq:flightassignment}:
\begin{align}\label{eq:flightassignment}
\sum_{k \in \Zi} \rik = 1, &&\forall i \in \A.
\end{align}

Let $\PPS$ be the set of aircraft pairs which are not conflict-free and may share the same FL, i.e. $\PPS = \{(i,j) \in \PS\cup\PI : \Zi \cap \Zj \neq \emptyset\}$. Let $\vpij$ be the binary variable defined as: 
\begin{equation}\label{eq:flightshare}
\vpij=\begin{cases}
1 & \text{ if aircraft $i$ and $j$ are assigned to the same FL},\\
0 & \text{ otherwise}.
\end{cases}
\end{equation}

Variable $\vpij$ can be linked to binary variables $\rik$ and $\rjk$ via the constraint:
\begin{equation}
\label{eq:flightshare2}
\rik +\rjk \leq \vpij + 1, \quad \forall (i,j) \in \PPS, \forall k \in \Zi \cap \Zj.
\end{equation}

The separation conditions \eqref{eq:zd}-\eqref{eq:rs} are  rewritten to account for altitude separation as follows: 
\begin{subequations}\label{eq:sepdis}
\begin{align}
\vijy\hat{x}_{ij} - \vijx\hat{y}_{ij} \leq 0, &\quad \text{ if } \z=1 \quad \text{and} \quad  \vpij =1, &&\forall (i,j) \in \PPS,\\
\vijy\hat{x}_{ij} - \vijx\hat{y}_{ij} \geq 0, &\quad \text{ if } \z=0   \quad \text{and} \quad  \vpij =1, &&\forall (i,j) \in \PPS,\\
\vijy \gamma_{ij}^l - \vijx \phi_{ij}^l \leq 0, &\quad \text{ if } \z=1 \quad \text{and} \quad  \vpij =1, &&\forall (i,j) \in \PPS,\\
\vijy \gamma_{ij}^u - \vijx \phi_{ij}^u \geq 0, &\quad \text{ if } \z=0 \quad \text{and} \quad \vpij =1, &&\forall (i,j) \in \PPS.
\end{align}
\end{subequations}

FL changes are typically less desirable compared to other deconfliction maneuvers such as speed or heading control \citep{bilimoria1996effects,hu2002optimal,alonso2011collision}. This is due to a number of practical considerations including an increase in fuel consumption, passenger discomfort due to climbing or descending and the need for extended monitoring. Under these considerations, we propose a lexicographic optimization approach for the 2D+FL ACRP. The number of FL changes is first minimized; before the total 2D deviation of flights is minimized. The proposed objective function for minimizing the number of FL changes is:
\begin{equation}\label{eq:obj3d}
\text{minimize} \sum_{i \in \A} \left|\sum\limits_{k \in \Zi} k\rik - \hat{\rho}_i\right|.
\end{equation}

The resulting nonconvex lexicographic 2D+FL formulation is summarized in Model \ref{mod:3dnonconvex}.

\begin{model}[Nonconvex Lexicographic 2D+FL Formulation]
\label{mod:3dnonconvex}
\begin{subequations}
\begin{align}
&\emph{1. Minimize Total FL Deviation } \eqref{eq:obj3d}, \nonumber \\
&\emph{2. Minimize Total 2D Deviation } \eqref{eq:obj2d}, \nonumber \\		
&\emph{Subject to:}  && \nonumber \\
&\emph{Motion Equations} \quad \eqref{eq:speedxdy}, \nonumber \\
&\emph{Separation Conditions} \quad \eqref{eq:sepdis}, \nonumber \\
&\emph{Speed and Heading Control Constraints} \quad  \eqref{eq:boundq}, \eqref{eq:boundtheta}, \nonumber\\
&\emph{Variable Bounds} \quad \eqref{eq:boundsdxdy}, \nonumber \\
&\emph{Flight Level Control} \quad \eqref{eq:flightassignment}, \eqref{eq:flightshare2}, \nonumber\\
& \vijx,\vijy \in \mathbb{R}, &&\forall (i,j) \in \PPS, \\
& \z, \vpij \in \{0,1\}, &&\forall (i,j) \in \PPS, \\
& \dx,\dy \in \R, &&\forall i \in \A, \\
& \rik \in \{0,1\}, &&\forall i \in \A, k \in \Zi.
\end{align}
\end{subequations}
\end{model}

Compared to Model \ref{mod:2dnonconvex}, Model \ref{mod:3dnonconvex} requires additional binary decision variables $\rik$ and $\vpij$. The former are used to assign aircraft to separated FLs and the latter ensures that aircraft sharing the same FL are separated via the 2D separation conditions.\\ 

We next propose exact solution methods for these 2D and 2D+FL nonconvex aircraft conflict formulations.

\section{Exact Solution Methods for the Aircraft Conflict Resolution Problem}
\label{methods}

We present exact solution methods for the 2D and the 2D+FL ACRPs that build on and extend the convex relaxations presented by \citet{rey2017complex}. We first present a convex relaxation of the 2D ACRP that fully relaxes the speed control constraint in the complex number formulation (Section \ref{miqp}). This relaxation yields a Mixed-Integer Quadratic Program (MIQP) which solution may violate the speed control bounds of the problem. To eliminate these potential violations, we incorporate convex quadratic constraints and piecewise linear outer approximations in a Mixed-Integer Quadratically Constrained Program (MIQCP). We then propose a constraint generation algorithm to iteratively refine the piece-wise linear approximation and show that our approach converges to optimal solutions of the 2D ACRP (Section \ref{miqcp}). 

To solve the lexicographic optimization formulation for the 2D+FL ACRP, we propose a two-step decomposition approach. We first solve a restricted flight assignment problem that only implicitly accounts for aircraft trajectories and yields an optimal solution with regards to the first objective function (total FL deviation). We then use the optimal solution of the flight assignment formulation to assign aircraft to FLs and solve a series of 2D problems (one per FL) to construct an optimal solution with regards to the second objective function (total 2D deviation). Both steps are iterated until a global solution is found (Section \ref{3dalgo}). We next present in detail the proposed exact solution methods.

\subsection{Mixed-Integer Quadratic Relaxation for the 2D Aircraft Conflict Resolution Problem}
\label{miqp}

An initial convex relaxation of Model \ref{mod:2dnonconvex} was proposed by \citet{rey2017complex} by relaxing the speed control constraint \eqref{eq:boundq}. The resulting formulation is a MIQP summarized in Model \ref{mod:2dmiqp}.

\begin{model}[MIQP 2D Formulation]
\label{mod:2dmiqp}
\begin{subequations}
\begin{align}
&\emph{Minimize Total 2D Deviation} \quad \eqref{eq:obj2d}, \nonumber \\
&\emph{Subject to:}  && \nonumber \\
&\emph{Motion Equations} \quad \eqref{eq:speedxdy}, \nonumber\\
&\emph{Separation Conditions} \quad \eqref{eq:zd}, \eqref{eq:rs}, \nonumber \\
&\emph{Heading Control Bounds} \quad  \eqref{eq:boundtheta}, \nonumber\\
&\emph{Variable Bounds} \quad \eqref{eq:boundsdxdy}, \nonumber \\
& \vijx,\vijy \in \B, &&\forall (i,j) \in \PS, \\
& \z \in \{0,1\}, &&\forall (i,j) \in \PS, \\
& \dx, \dy \in \R, && \forall i \in \A. 
\end{align}
\end{subequations}
\end{model}

Model \ref{mod:2dmiqp} yields a lower bound on the optimal objective value of Model \ref{mod:2dnonconvex} and a solution which is a global optimum if the relaxed constraint \eqref{eq:boundq} is not violated. 


\subsection{Mixed-Integer Quadratically Constrained Relaxation and Constraint Generation Algorithm}
\label{miqcp}

To tighten the MIQP relaxation given in Section \ref{miqp}, we build on and extend the MIQCP relaxation proposed by \citet{rey2017complex} by incorporating the speed control constraint \eqref{eq:boundq} using convex quadratic and piecewise linear constraints. Observe that the speed upper bound constraint \eqref{eq:uboundq} is convex quadratic, hence it can be incorporated directly in the MIQCP formulation. 


To incorporate the speed lower bound constraint \eqref{eq:lboundq}, for each aircraft $i \in \A$, we introduce real variables $\tdx$ and $\tdy$ to approximate $\dx^2$ and $\dy^2$, respectively, via convex quadratic constraints:
\begin{subequations}\label{eq:tdxtdy}
\begin{align}
& \tdx \geq \dx^2, && \forall i \in \A,\label{eq:dx} \\
& \tdy \geq \dy^2, && \forall i \in \A.\label{eq:dy}
\end{align}
\end{subequations}

Accordingly, we require:
\begin{align}
\qmin^2 \leq \tdx + \tdy, && \forall i \in \A. 
\label{eq:dxdylower}
\end{align}

To impose upper bounds on $\tdx$ and $\tdy$, we use their McCormick envelopes \citep{mccormick1976computability}:
\begin{subequations}
\begin{align}
& \tdx \leq (\qmax + \qmin\cos(\max\{|\underline\ti|,|\overline\ti|\}))\dx - \qmax\qmin\cos(\max\{|\underline\ti|,|\overline\ti|\}), && \forall i \in \A. \label{eq:tdx}\\
& \tdy \leq \qmax(\sin(\underline{\ti}) + \sin(\overline{\ti}))\dy - \qmax^2\sin(\overline{\ti})\sin(\underline{\ti}), && \forall i \in \A. \label{eq:tdy}
\end{align}
\label{eq:linearcuts}
\end{subequations}

\begin{figure}[t]
\centering
\subfloat[Graph of $\dx^2$ (in red) over the domain given by Eq. \eqref{eq:bounddx}. The mixed-integer cuts \eqref{eq:micdx} are illustrated for a partition of $|\gx|=2$ segments.]{\includegraphics[width=5.5cm,height=4cm]{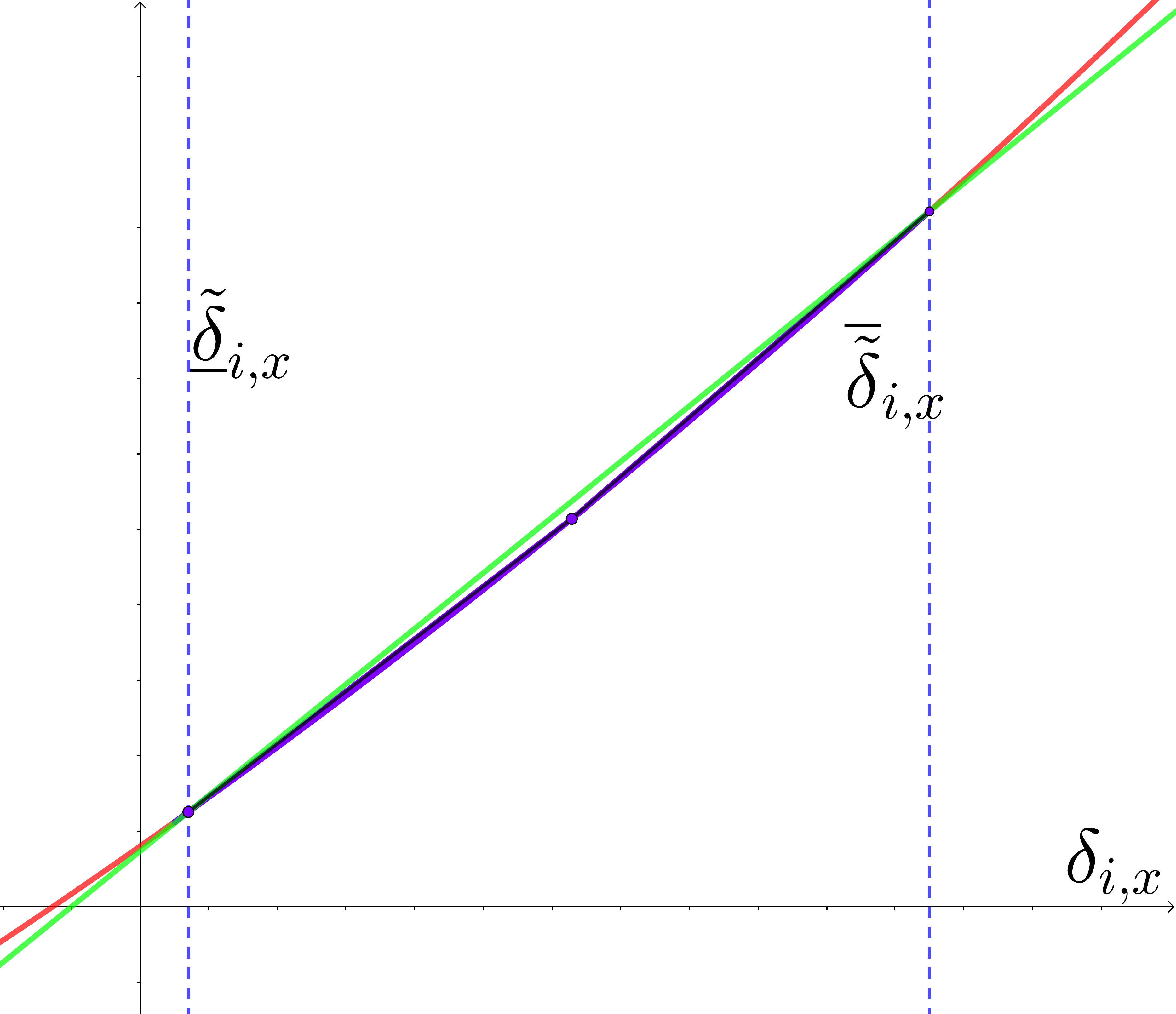}\label{fig:tdxcuts}}%
\qquad
\subfloat[Graph of $\dy^2$ (in red) over the domain given by Eq. \eqref{eq:bounddy}. The mixed-integer cuts \eqref{eq:micdy} are illustrated for a partition of $|\gy|=4$ segments.]{\includegraphics[width=7.5cm,height=4cm]{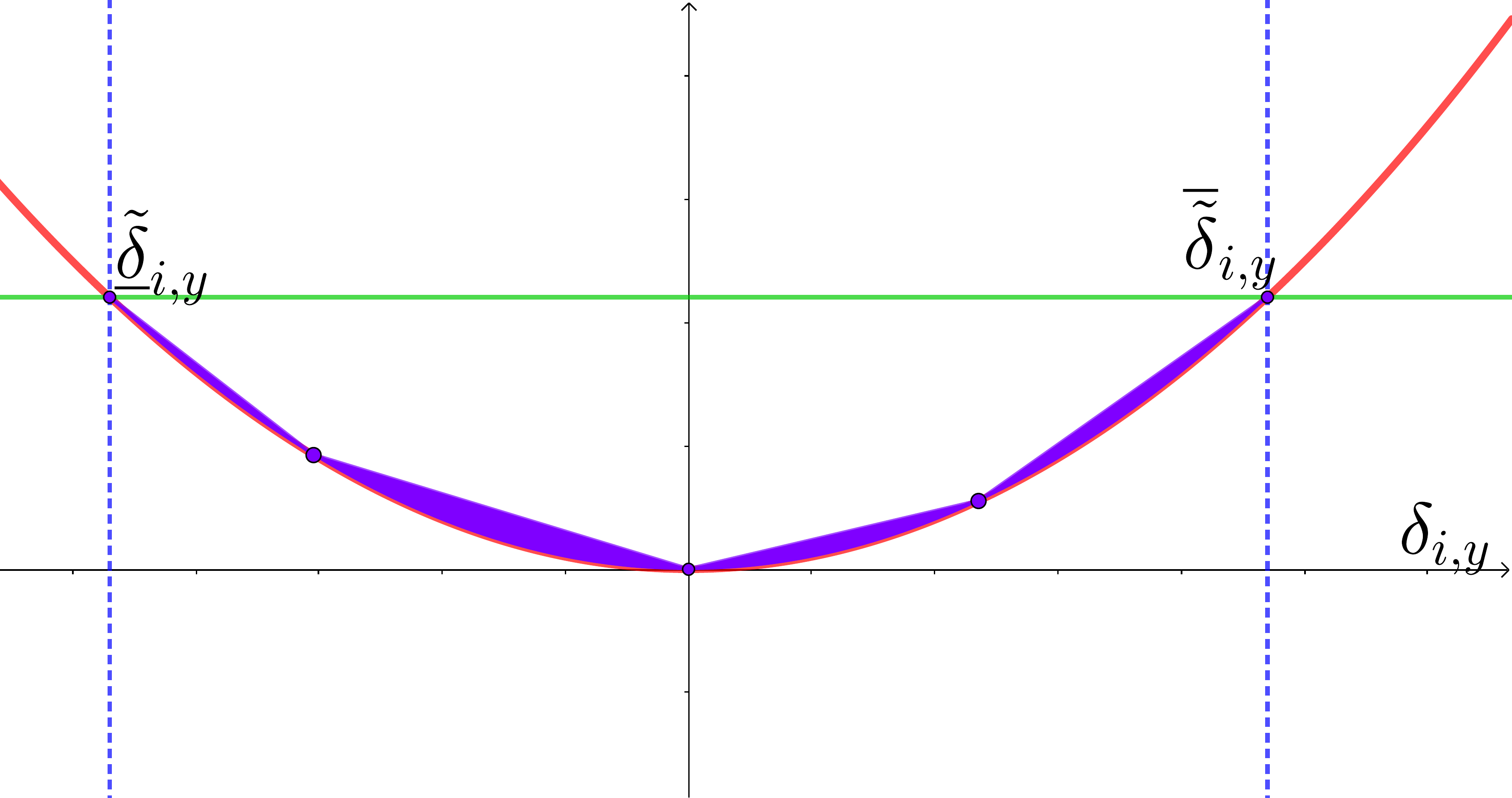}\label{fig:tdycuts}}%
\caption{Piecewise linear approximation of $\dx^2$ and $\dy^2$ (in red). The green lines represent the initial McCormick relaxation given by \eqref{eq:linearcuts}. The purple shaded regions represent the refined feasible region of $\tdx$ and $\tdy$ after imposing the mixed-integer cuts \eqref{eq:mic}. Control bounds are $\qmin = 0.94$, $\qmax = 1.03$, $\tmin = -\frac{\pi}{6}$ and $\tmax = +\frac{\pi}{6}$.}%
\label{fig:tdxtdy}%
\end{figure}

Constraints \eqref{eq:tdxtdy}-\eqref{eq:linearcuts} restrict variables $\tdx$ and $\tdy$ to convex regions thus providing an initial relaxation of the speed lower bound constraint \eqref{eq:lboundq}. This initial relaxation is illustrated in Figure \ref{fig:tdxtdy} which depicts the variation of $\dx^2$ and $\dy^2$ (in red) over the domain of $\dx$ and $\dy$ for realistic speed and heading control bounds. The green lines in Figures \ref{fig:tdxcuts} and \ref{fig:tdycuts} represent the initial McCormick envelopes. This relaxation may still yield infeasible aircraft speeds. To refine this convex relaxation, we introduce mixed-integer cuts that can be generated on-the-fly in a constraint generation algorithm.

We next present the general structure of the proposed mixed-integer cuts before discussing how these cuts are generated iteratively. Let $\gx$ and $\gy$ be partitions of the domain of variables $\dx$ and $\dy$, respectively. Since $\dx^2$ and $\dy^2$ are convex, any line joining two extremities of a segment in $\gx$ and $\gy$ refines the initial McCormick upper envelopes. Let $\axk$ (resp. $\ayk$) and $\bxk$ (resp. $\byk$) be the slope and the intercept corresponding to segment $k \in \gx$ (resp. $k \in \gy$). Further, let $\sxk$ (resp. $\syk$) be a binary variable taking value 1 if $\dx$ (resp. $\dy$) belongs to segment $k$ of partition $\gx$ (resp. $\gy$). The proposed mixed-integer cuts take the form of:
\begin{subequations}
\begin{align}
\tdx \leq \axk\dx + \bxk, \quad \text{ if } \quad \sxk = 1, \quad &\forall i \in \A, k \in \gx, \label{eq:micdx}\\
\tdy \leq \ayk\dy + \byk, \quad \text{ if } \quad \syk = 1, \quad &\forall i \in \A, k \in \gy. \label{eq:micdy}
\end{align}
\label{eq:mic}
\end{subequations}

Let $\dxlk$ and $\dxuk$ (resp. $\dylk$ and $\dyuk$) be the extremities of segment $k \in \gx$ (resp. $k \in \gy$). Binary variables $\sxk$ and $\syk$ are defined as:
\begin{subequations}
\begin{align}
\sxk =\begin{cases}
1 & \text{ if }\dx \in [\dxlk,\dxuk[,\\
0 & \text{ otherwise}, \\
\end{cases}\quad \forall i \in \A, k \in \gx,\label{eq:defdx}\\
\syk =\begin{cases}
1 & \text{ if }\dy \in [\dylk,\dyuk[,\\
0 & \text{ otherwise}.
\end{cases}\quad \forall i \in \A, k \in \gy.\label{eq:defdy}
\end{align}
\label{eq:def}
\end{subequations}

We require the following cut selection constraints:
\begin{subequations}
\begin{align}
\sum_{k \in \gx} \sxk = 1 \quad &\forall i \in \A, \label{eq:bindx}\\
\sum_{k \in \gy} \syk = 1 \quad &\forall i \in \A. \label{eq:bindy}
\end{align}
\label{eq:bin}
\end{subequations}

The proposed mixed-integer cuts are illustrated in Figures \ref{fig:tdxcuts} and \ref{fig:tdycuts} for the case of $|\gx|=2$ and $|\gy|=4$ segments, respectively. The shaded purple region represents the feasible region of $\tdx$ and $\tdy$ after imposing the mixed-integer cuts. The resulting formulation is an MIQCP with mixed-integer cuts summarized in Model \ref{mod:miqcp2d}. This formulation can be solved by off-the-shelf commercial optimization software and provides a relaxation of the Model \ref{mod:2dnonconvex} which can be tightened as desired by refining the partitions $\gx$ and $\gy$ of the domain of variables $\dx$ and $\dy$, respectively. 

\begin{model}[MIQCP 2D Formulation]
\label{mod:miqcp2d}
\begin{subequations}
\begin{align}
&\emph{Minimize Total 2D Deviation} \quad \eqref{eq:obj2d}, \nonumber \\
&\emph{Subject to:}  && \nonumber \\
&\emph{Motion Equations} \quad \eqref{eq:speedxdy}, \nonumber\\
&\emph{Separation Conditions} \quad \eqref{eq:zd}, \eqref{eq:rs}, \nonumber \\
&\emph{Heading Control Bounds} \quad  \eqref{eq:boundtheta}, \nonumber\\
& \emph{Speed Control Upper Bound} \quad  \eqref{eq:uboundq}, \nonumber\\
&\emph{Speed Control Relaxed Lower Bound} \quad \eqref{eq:tdxtdy},\eqref{eq:dxdylower},\eqref{eq:linearcuts}, \nonumber \\
&\emph{Mixed Integer Cuts} \quad \eqref{eq:mic},\eqref{eq:def},\eqref{eq:bin}, \nonumber \\
&\emph{Variable Bounds} \quad \eqref{eq:boundsdxdy}, \nonumber \\
& \vijx,\vijy \in \B, &&\forall (i,j) \in \PS, \\
& \z \in \{0,1\}, &&\forall (i,j) \in \PS, \\
& \dx, \dy \in \R, && \forall i \in \A,\\
& \tdx, \tdy \in \R, && \forall i \in \A, \\ 
& \sxk \in \{0,1\}, &&\forall i \in \A, k \in \gx, \\
& \syk \in \{0,1\}, &&\forall i \in \A, k \in \gy.
\end{align}
\end{subequations}
\end{model}

We propose a constraint generation algorithm which starts from a relaxed formulation and iteratively refines the piecewise linear outer approximation of the quadratic terms $\dx^2$ and $\dy^2$ via mixed-integer cuts \eqref{eq:mic}. At each iteration, the proposed constraint generation algorithm examines the solution $(\dx, \dy)$ of each aircraft $i \in \A$ for speed lower bound violations. If $\dx^2 + \dy^2 < \qmin^2$, then at least one of the relaxed auxiliary variables $\tdx$ or $\tdy$ must be such that $\tdx > \dx^2$ or $\tdy > \dy^2$. To eliminate the current infeasible solution, the partition(s) $\gx$ or $\gy$ of the violating variable(s) is augmented by dividing the segment(s) corresponding to $\dx$ or $\dy$ into two segments which meet at $\dx$ or $\dy$ and corresponding binary variable(s) $\sxk$ or $\syk$ are added to the formulation. The process is repeated until all aircraft have feasible speed profiles which corresponds to a global optimum of Model \ref{mod:2dnonconvex}. To further improve the convergence of the solution algorithm, the heuristic procedure outlined by \citet{rey2017complex} is used to attempt to find improving upper bounds. After each solve of the relaxed formulation \ref{mod:miqcp2d}, Model \ref{mod:2dnonconvex} is solved as a nonlinear program (NLP) by fixing variable $\z$. The pseudo-code of the resulting algorithm is summarized in Algorithm \ref{algo:2d}.

\begin{algorithm}
\KwIn{$\A$, $\hat{\bm{\theta}}$, $\hat{\bm{v}}$, $\qlb$, $\qub$, $\tlb$, $\tub$, $\epsilon$}
\KwOut{$\bm{q}^\star$, $\bm{\theta}^\star$, LB, UB}
$\P, \PI,\PF,\PS \gets$ Algorithm \ref{algo:prepro}\\
LB $\gets 0$ \\
UB $\gets +\infty$ \\
$\bm{q}, \bm{\theta}, \bm{\delta}_x, \bm{\delta}_y, \bm{z}$, LB $\gets$ Solve MIQP \ref{mod:2dmiqp} and calculate $\bm{q}$ and $\bm{\theta}$ from $\bm{\delta}_x$ and  $\bm{\delta}_y$\\
\If{\emph{MIQP \ref{mod:2dmiqp} is infeasible}}{
    Return \textsc{infeasible}\\
}
\If{$\bm{q}$ \emph{is feasible (no speed violation)}}{
    UB $\gets$ LB\\
	$\bm{q}^\star \gets \bm{q}$\\
	$\bm{\theta}^\star \gets \bm{\theta}$\\		
}
\Else{
    \textsc{converged} $\gets$ \texttt{False}\\
    \While{\textsc{converged} $=$ \emph{\texttt{False}}}{
    	$\bm{q}', \bm{\theta}'$, UB-NLP $\gets$ Solve Model \ref{mod:2dnonconvex} as NLP with fixed $\bm{z}$\\
    	\If{\emph{UB-NLP} $<$ \emph{UB}}{
    		UB $\gets$ UB-NLP\\
    		$\bm{q}^\star \gets \bm{q}'$\\
    		$\bm{\theta}^\star \gets \bm{\theta}'$\\
    	}
    	\For{$i \in \A$}{
    	    \If{$\dx^2 + \dy^2 > \qmax^2$}{
    	        Add constraint \eqref{eq:uboundq}
    	    }		
    	    \ElseIf{$\dx^2 + \dy^2 < \qmin^2$}{
    	        \If{$\tdx > \dx^2$}{
    	        Add segment to $\gx$ at $\dx$, variable $\sxk$, constraints \eqref{eq:micdx}, \eqref{eq:defdx}, and update \eqref{eq:bindx} \\
    	        }
    	        \If{$\tdy > \dy^2$}{
    	        Add segment to $\gy$ at $\dy$, variable $\syk$, constraints \eqref{eq:micdy}, \eqref{eq:defdy}, and update \eqref{eq:bindy}\\
    	        }
    	    }
    	}
    	$\bm{q}, \bm{\theta}, \bm{\delta}_x, \bm{\delta}_y, \tilde{\bm{\delta}}_x, \tilde{\bm{\delta}}_y, \bm{z}$, LB $\gets$ Solve MIQCP \ref{mod:miqcp2d} and calculate $\bm{q}$ and $\bm{\theta}$ from $\bm{\delta}_x$ and  $\bm{\delta}_y$\\
        \If{\emph{MIQCP \ref{mod:miqcp2d} is infeasible}}{
            Return \textsc{infeasible}\\
        }    	
        \If{$\bm{q}$ \emph{is feasible (no speed violation)}}{
            UB $\gets$ LB\\
        	$\bm{q}^\star \gets \bm{q}$\\
        	$\bm{\theta}^\star \gets \bm{\theta}$\\
        	\textsc{converged} $\gets$ \texttt{True}\\
        }
    	\If{\emph{(UB-LB)/UB} $\leq \epsilon$}{
    	    \textsc{converged} $\gets$ \texttt{True}\\
    	}        
    }
}
\caption{Solution algorithm for the 2D ACRP}
\label{algo:2d}
\end{algorithm}

\subsection{Decomposition Algorithm for the 2D+FL Aircraft Conflict Resolution Problem}
\label{3dalgo}

We next introduce a two-step decomposition approach for the nonconvex lexicographic 2D+FL conflict resolution problem represented by Model \ref{mod:3dnonconvex}. The first objective function \eqref{eq:obj3d} focuses on minimizing the number of FL changes. Observe that this objective function is null and minimal if all aircraft $i \in \A$ can remain at their initial FL $\hat{\rho}_i$. Since minimizing aircraft FL re-assignment is the highest priority objective function, we need only to identify combinations of aircraft which are non-separable in 2D and ensure that such combinations are not assigned to the same FL. To address this first step, we propose a compact aircraft FL assignment model that implicitly accounts for 2D non-separable aircraft combinations. Let $\Omega_{\textrm{I}} \subseteq 2^{\A}$ be the set of aircraft combinations which are 2D non-separable, i.e. for any $\omega \in \Omega_{\textrm{I}}$, the subset of aircraft $\omega$ cannot be separated in 2D. To ensure that all aircraft are assigned to 2D separable FLs, we require:
\begin{equation}
\sum_{i \in \omega} \rik \leq |\omega|-1, \quad \forall \omega \in \Omega_{\textrm{I}}, \forall k \in \bigcap_{i \in \omega} \Zi.
\label{eq:sepfl}
\end{equation}


Let $\Delta\rho_{i} \geq 0$ be a variable representing the absolute flight level deviation for aircraft $i \in \A$. Objective function \eqref{eq:obj3d} can be linearized as:
\begin{equation}\label{eq:obj3dlin}
\text{minimize} \sum_{i \in \A} \Delta\rho_{i},
\end{equation}

subject to the constraints: 
\begin{subequations}
\begin{align}
\Delta\rho_{i} \geq \sum\limits_{k \in \Zi} k\rik - \hat{\rho}_i, &\quad \forall i \in \A, \\
\Delta\rho_{i} \geq \hat{\rho}_i - \sum\limits_{k \in \Zi} k\rik, &\quad \forall i \in \A.
\end{align}
\label{eq:lin}
\end{subequations}

Combining the FL separation constraint \eqref{eq:sepfl} and the above linearized objective function, yields a compact FL assignment formulation summarized in Model \ref{mod:fa} which is a MILP with an exponential number of constraints. 

\begin{model}[FL Assignment Formulation]
\label{mod:fa}
\begin{subequations}
\begin{align}
&\emph{Minimize Total FL Deviation}\quad \eqref{eq:obj3dlin} \nonumber\\
&\emph{Subject to:}  && \nonumber \\
&\emph{FL Control} \quad \eqref{eq:flightassignment}, \nonumber\\ 
&\emph{FL Separation} \quad \eqref{eq:sepfl}, \nonumber\\
&\emph{Linearization} \quad \eqref{eq:lin}, \nonumber\\
& \rik \in \{0,1\}, &&\forall i \in \A, k \in \Zi,\\
& \Delta\rho_i \geq 0, &&\forall i \in \A.
\end{align}
\end{subequations}
\end{model}

We propose to solve Model \ref{mod:fa} by first restricting Constraint \eqref{eq:sepfl} to subsets of size two, i.e. 2D non-separable aircraft pairs, which can be efficiently identified using Algorithm \ref{algo:prepro} in pre-processing. We then decompose the 2D+FL problem into a series of 2D conflict resolution problems, one per FL, based on the optimal solution $\bm{\rho}^\star$ of the relaxed Model \ref{mod:fa}. For each FL, we solve the corresponding 2D problem using the exact constraint generation approach of Algorithm \ref{algo:2d} with the aircraft set $\A_{k} = \{i \in \A: \rho^\star_{ik} = 1\}$. If Algorithm \ref{algo:2d} returns \textsc{infeasible} for FL $k$, then the corresponding FL separation constraint \eqref{eq:sepfl} with $\omega = \A_k$ is generated and Model \ref{mod:fa} is re-solved with the additional constraint(s). The process is repeated until all 2D problems are feasible. Let $\Z$ be the set of all FLs, i.e. $\Z = \cup_{k \in \Zi : i \in \A} \Zi$, we denote $\bm{q}_k$ and $\bm{\theta}_k$ the vectors of speed and heading controls for aircraft assigned to FL $k \in \Z$, respectively. The proposed approach for the lexicographic 2D+FL conflict resolution problem is summarized in Algorithm \ref{algo:3d}. Note that this decomposition approach offers the possibility to solve all 2D problems in parallel via the \textbf{for} loop at Line 8.\\

\begin{algorithm}[H]
\KwIn{$\A$, $[\Zi]_{i\in\A}$, $\hat{\bm{\theta}}$, $\hat{\bm{v}}$, $\hat{\bm{\rho}}$, $\qlb$, $\qub$, $\tlb$, $\tub$}
\KwOut{$\bm{q}^\star$, $\bm{\theta}^\star$, $\bm{\rho}^\star$}
$\P, \PF,\PS,\PI \gets$ Algorithm \ref{algo:prepro}\\
$\Omega_{\textrm{I}} \gets \{(i,j) \in \PI : \Zi \cap \Zj \neq
\emptyset\}$\\
$\Z \gets \cup_{k \in \Zi : i \in \A} \Zi$\\
\textsc{converged} $\gets$ \texttt{False}\\
\While{\textsc{converged} $=$ \emph{\texttt{False}}}{
    $\bm{\rho}^\star \gets$ Solve MILP \ref{mod:fa} \\
    \textsc{converged} $\gets$ \texttt{True}\\
    \For{$k \in \Z$}{
        $\A_{k} \gets \{i \in \A: \rho^\star_{ik} = 1\}$ \\
        $\bm{q}_k^\star, \bm{\theta}_k^\star \gets$ Algorithm \ref{algo:2d} with $\A = \A_k$\\
        \If{\emph{Algorithm \ref{algo:2d} returns} \textsc{infeasible}}{
            $\Omega_{\textrm{I}} \gets \Omega_{\textrm{I}} \cup \A_k$\\
            \textsc{converged} $\gets$ \texttt{False}\\
        }
    }
}
$\bm{q}^\star, \bm{\theta}^\star \gets [\bm{q}_k^\star]_{k\in\Z}, [\bm{\theta}_k^\star]_{k \in \Z}$\\
\caption{Solution algorithm for the lexicographic 2D+FL conflict resolution problem}
\label{algo:3d}
\end{algorithm}

\section{Numerical Results}
\label{num}

We first introduce the experimental framework used to test the proposed mixed-integer formulations and algorithms in Section \ref{ed}. We then explore the behavior of the proposed 2D objective function in Section \ref{sa}. Numerical results for the 2D ACRP are presented in Section \ref{2dresults}, and results for the 2D+FL problem are presented in Section \ref{3dresults}.

\subsection{Experiments Design}
\label{ed}

\begin{figure}[t]	
	\centering	
		\subfloat[CP-7: circle problem with 7 aircraft]{\includegraphics[width=0.22\textwidth]{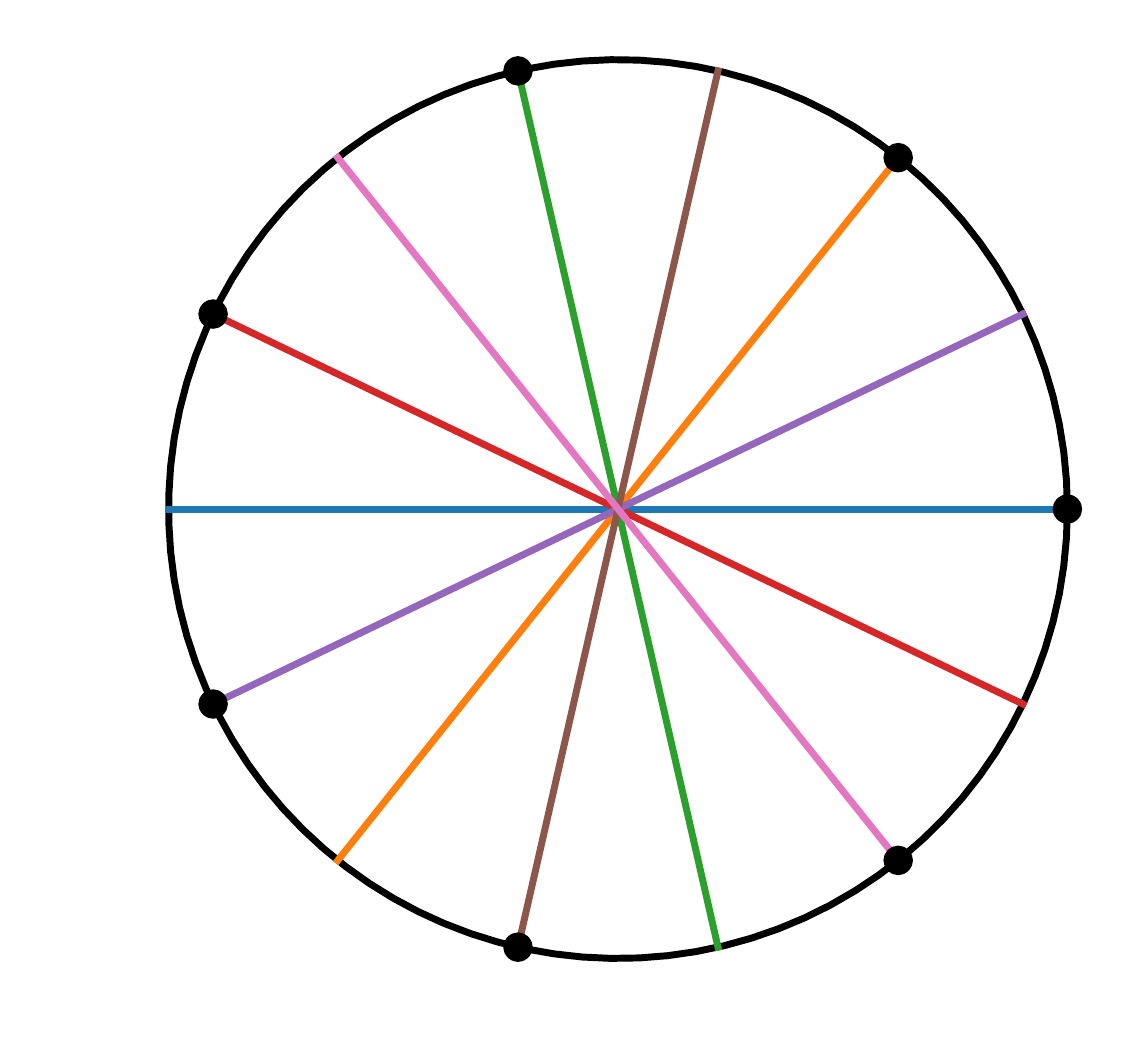}\label{fig:cp}} 
		\hfil
		\subfloat[RCP-10: random circle problem with 10 aircraft]{\includegraphics[width=0.22\textwidth]{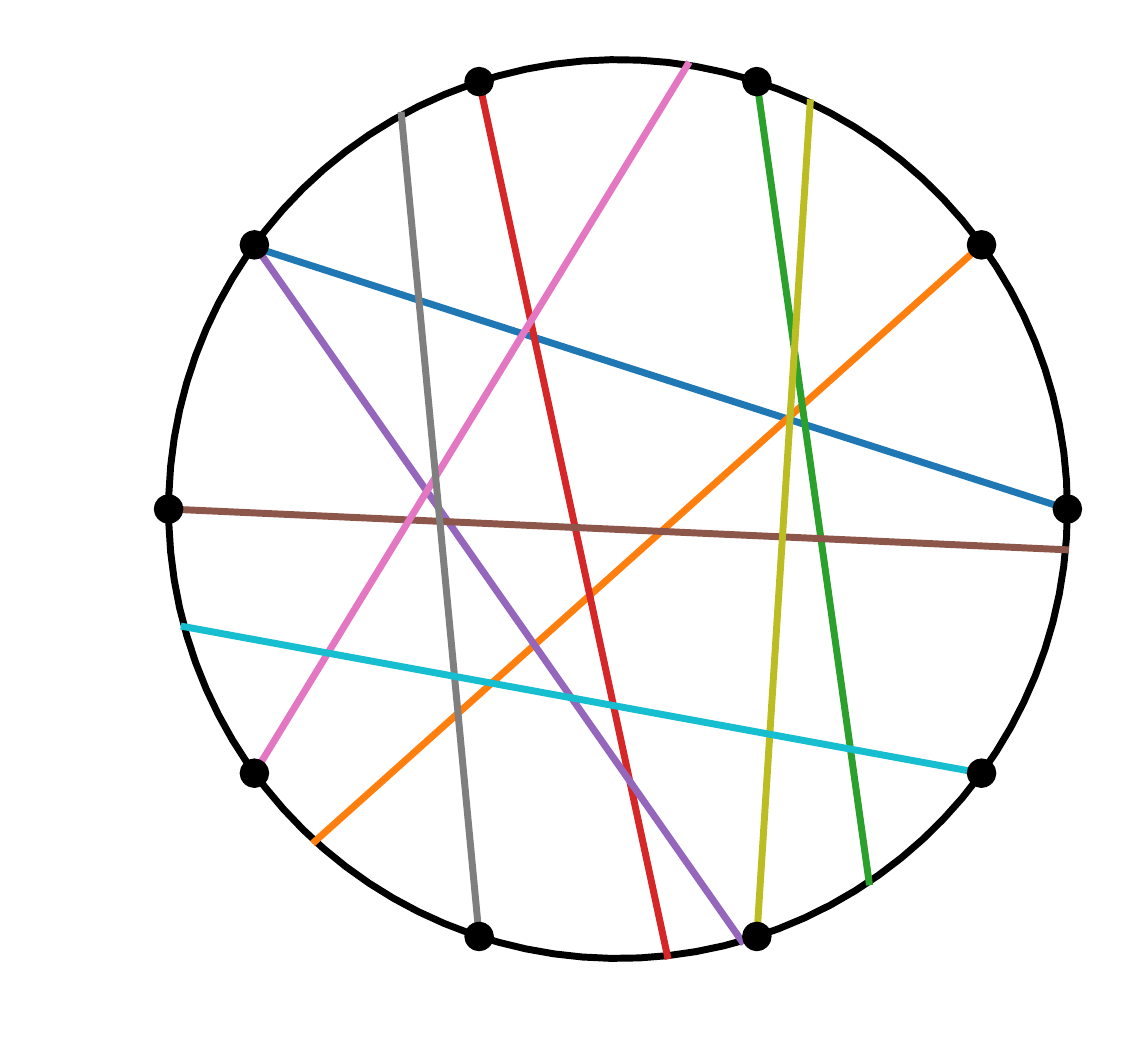} \label{fig:rcp10-5}}
		\hfil
		\subfloat[FP-5: flow problem with 10 aircraft and $\alpha=\frac{\pi}{6}$]{\includegraphics[width=0.22\textwidth]{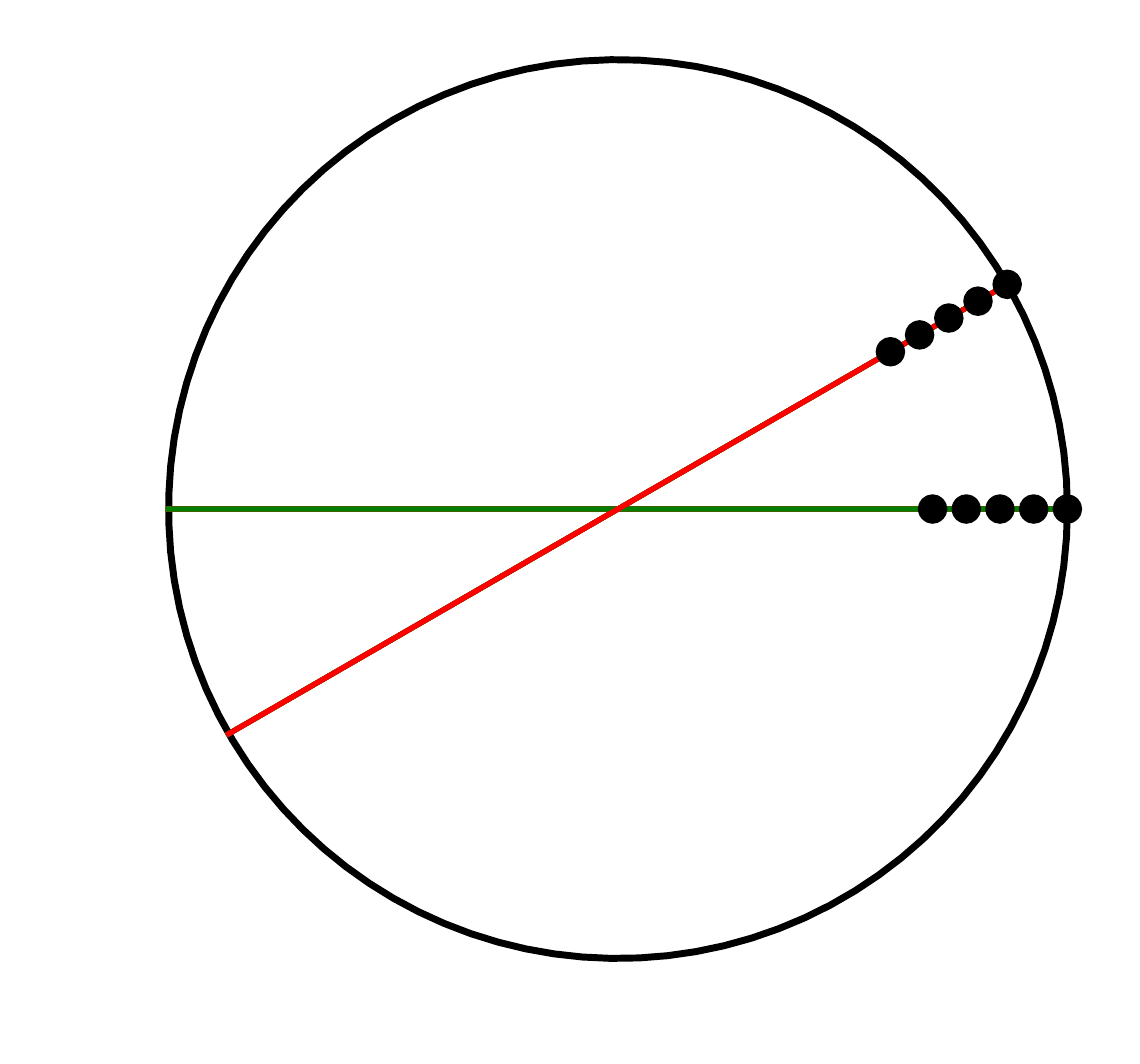} \label{fig:fp}}
		\hfil
		\subfloat[GP-8: grid problem with 32 aircraft and $\alpha=\frac{\pi}{2}$]{\includegraphics[width=0.22\textwidth]{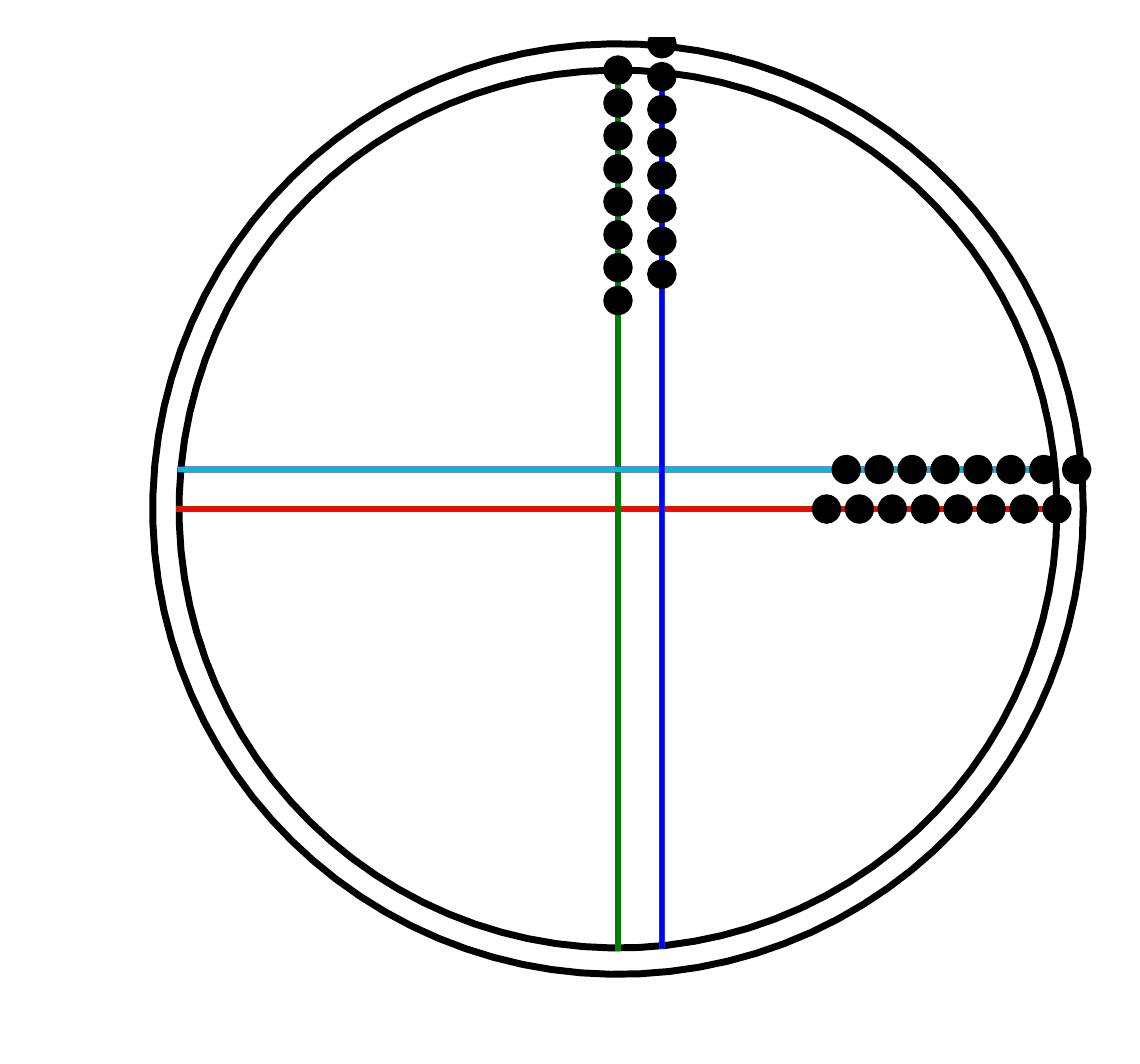} \label{fig:gp}}
	\caption{Example of 2D benchmarking instances for the Circle Problem (CP), Random Circle Problem (RCP), Flow Problem (FP) and Grid Problem (GP).}
	\label{case1}
\end{figure}

We test the performance of the proposed mixed-integer formulations and algorithms using four benchmarking problems from the literature: the Circle Problem (CP), the Flow Problem (FP), the Grid Problem (GP) and the Random Circle Problem (RCP). The four types of benchmarking instances are illustrated in Figure \ref{case1}. The CP consists of a set of aircraft uniformly positioned on the circle heading towards its centre. Aircraft speeds are assumed to be identical, hence the problem is highly symmetric (see Fig. \ref{fig:cp}). The CP is notoriously difficult due to the geometry of aircraft initial configuration has been widely used for benchmarking CD\&R algorithms in the literature \citep{durand2009ant,rey2015equity,cafieri2017mixed,cafieri2017maximizing,rey2017complex}. To break the symmetry of CP benchmarking instances, \citet{vanaret2012benchmarking} introduced the RCP which builds on the same framework, but aircraft initial speeds and headings are randomly deviated within specified ranges to create randomized instances with less structure (see Fig. \ref{fig:rcp10-5}). CP and RCP instances are named CP-N and RCP-N, respectively, where N is the total number of aircraft.  \citet{lehouillier2017two} formally introduced two additional structured problems which aim to represent more realistic air traffic configurations: the FP and the GP. The FP consists of two streams of aircraft separated by an angle $\alpha$ and anchored on the circumference of a circle. In each stream, aircraft are separated by at least 5 NM (we use 15 NM in our experiments) from each other (see Figure \ref{fig:fp}). The GP consists of two FP instances separated by 15 NM diagonally (see Figure \ref{fig:gp}). FP and GP instances are named FP-N and GP-N, respectively, where N denotes the number of aircraft per stream.
 
In all experiments, we use a circle of radius 200NM. For CP, FP and GP instances, all aircraft have an initial speed of 500 NM/h. For RCP instances, aircraft initial speeds are randomly chosen in the range 486-594 NM/h and their initial headings are deviated from a radial trajectory (i.e. towards the centre of the circle) by adding a randomly chosen angle between $-\frac{\pi}{6}$ and $+\frac{\pi}{6}$. For FP and GP instances, we use $\alpha = \frac{\pi}{6}$ and $\alpha = \frac{\pi}{2}$, respectively, as proposed by \citet{lehouillier2017two}.

We report numerical results for problems with a subliminal speed control range of. $[-6\%, +3\%]$ \citep{bonini2009erasmus}. We consider two heading control ranges, we first assume that aircraft can modify their heading within the range $[-30^\circ,+30^\circ]$ as commonly used in the literature \citep{cafieri2017mixed,rey2017complex}, and we also consider a reduced heading control range of $[-15^\circ,+15^\circ]$. For conflict resolution problems with altitude control, we randomly assign each aircraft $i \in \A$ to a FL $\hat{\rho}_i \in \Z$, and we assume that only adjacent flight levels are available for aircraft, i.e. $\Zi = \{\hat{\rho}_i-1,\hat{\rho}_i,\hat{\rho}_i+1\}$.\\

The proposed disjunctive formulations are compared against the so-called shadow model, which uses the separation conditions introduced initially by \cite{pallottino2002conflict} and more recently adapted for speed, heading control altitude control by \cite{alonso2011collision, alonso2016exact}. In our implementation of the shadow formulation, the control variables and constraints are identical to that of the disjunctive formulation, and the only difference between both formulations is the set of separation constraints used, i.e. Constraints \eqref{eq:zd} and \eqref{eq:rs} are replaced with the shadow separation conditions, and the number of binary variables required to express these On/Off constraints. Specifically, for the 2D ACRP, the disjunctive formulation only requires a single binary variable per aircraft pair ($z_{ij}$) whereas the shadow formulation requires four binary variables per aircraft pair \citep{pallottino2002conflict,alonso2011collision,alonso2016exact}. Both formulations are implemented using the same pre-processing procedure (Algorithm \ref{algo:prepro}) to eliminate conflict-free aircraft pairs. A detailed formulation of our implementation of the 2D shadow formulation is provided in  \ref{shadow}. 

All 2D ACRPs are solved by implementing Algorithm \ref{algo:2d} with an optimality gap $\epsilon = 0.01$ and a time limit of 10 minutes. All 2D+FL problems are solved by implementing Algorithm \ref{algo:3d} which calls Algorithm \ref{algo:2d} using the same optimality gap as for the 2D problems, and we allocate a time limit of 10 minutes per FL. All models are implemented using \Python on a personal computer with 16 GB of RAM and an Intel i7 processor at 2.9GHz. The MIQPs and MIQCPs are solved with \Cplex v12.10 \citep{cplex2009v12} API for \Python using default options. For reproducibility purposes, all formulations and instances used are available at \small \url{https://github.com/acrp-lib/acrp-lib}.\normalsize

We next conduct a sensitivity analysis on the preference weight $w$ in the 2D objective function to explore its impact on aircraft trajectories in Section \ref{sa}. We present results on 2D problems in Section \ref{2dresults} and results on 2D+FL problems in Section \ref{3dresults}.

\subsection{Sensitivity Analysis on the Preference Weight $w$}
\label{sa}

To quantify the impact of the preference weight $w$ in the proposed 2D objective function \eqref{eq:obj2d}, we conduct numerical experiments on one instance of each of the four types of benchmarking instances for varying values of $w$. For this experiment, we focus on the typical heading control range $[-30^\circ,+30^\circ]$. We report the total speed deviation $\Sigma_q = \sum_{i \in \A} (1 - \qi)^2$, and the total heading deviation $\Sigma_\theta = \sum_{i \in \A} \ti^2$. Our goal is to show that by varying the preference weight $w \in \ ]0,1[$, the decision-maker can control the desired level of trade-off among total speed deviation and total heading deviation. Recall that in objective function \eqref{eq:obj2d}, $w$ is the coefficient of $\dy^2 = (\qi \sin (\ti))^2$ which is minimal for $\ti=0$; while $(1-w)$ is the coefficient of $(1 - \dx)^2 = (1 - \qi\cos(\ti))^2$ which is minimal for $\qi = 1$ and $\ti = 0$. Hence, one can expect that increasing (resp. decreasing) $w$ will tend to penalize heading (resp. speed) deviations more than speed (resp. heading) deviations. 
 
\begin{figure}[t]
    \centering
    \subfloat[CP-8]{\includegraphics[width=0.4\textwidth]{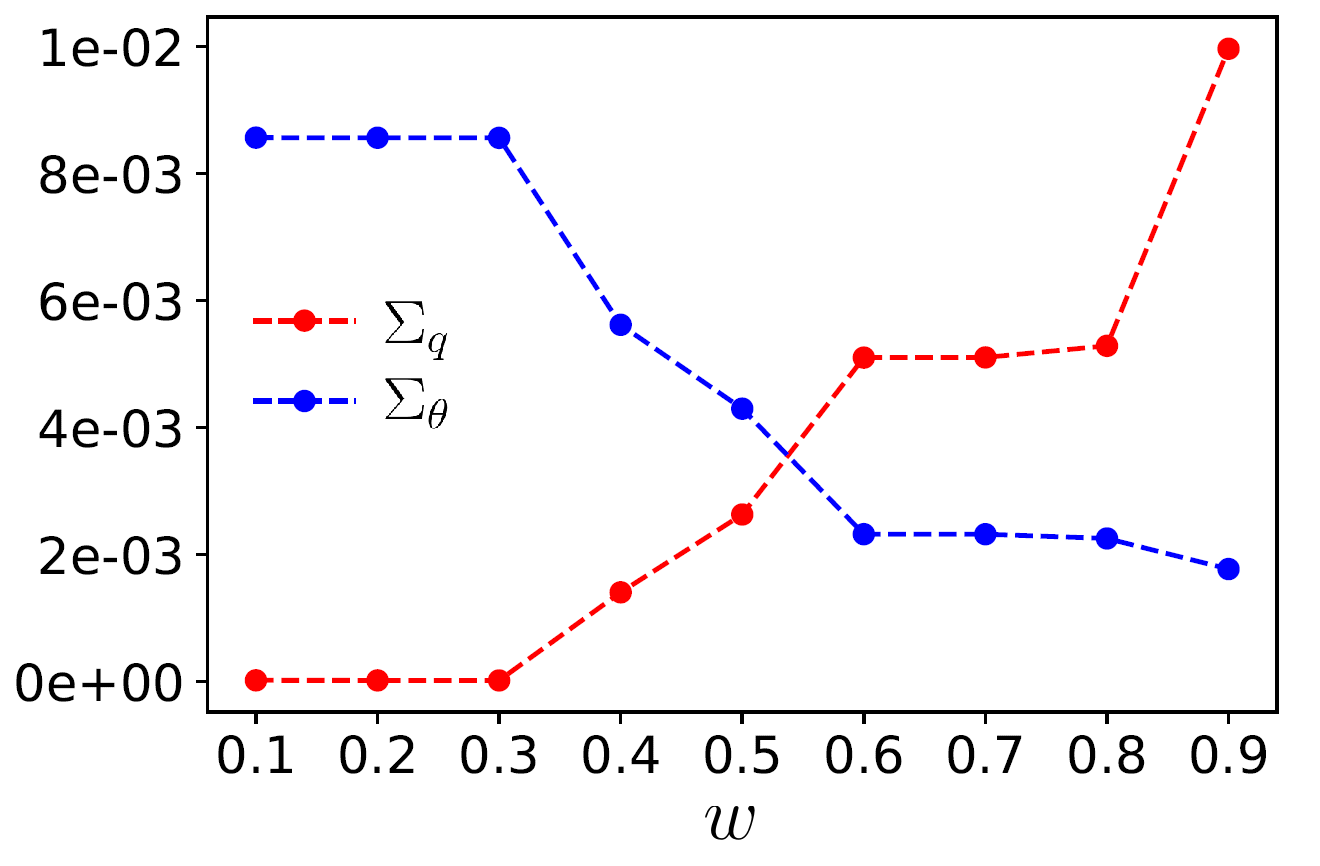}\label{cp8}}
    \hfil
    \subfloat[FP-10]{\includegraphics[width=0.4\textwidth]{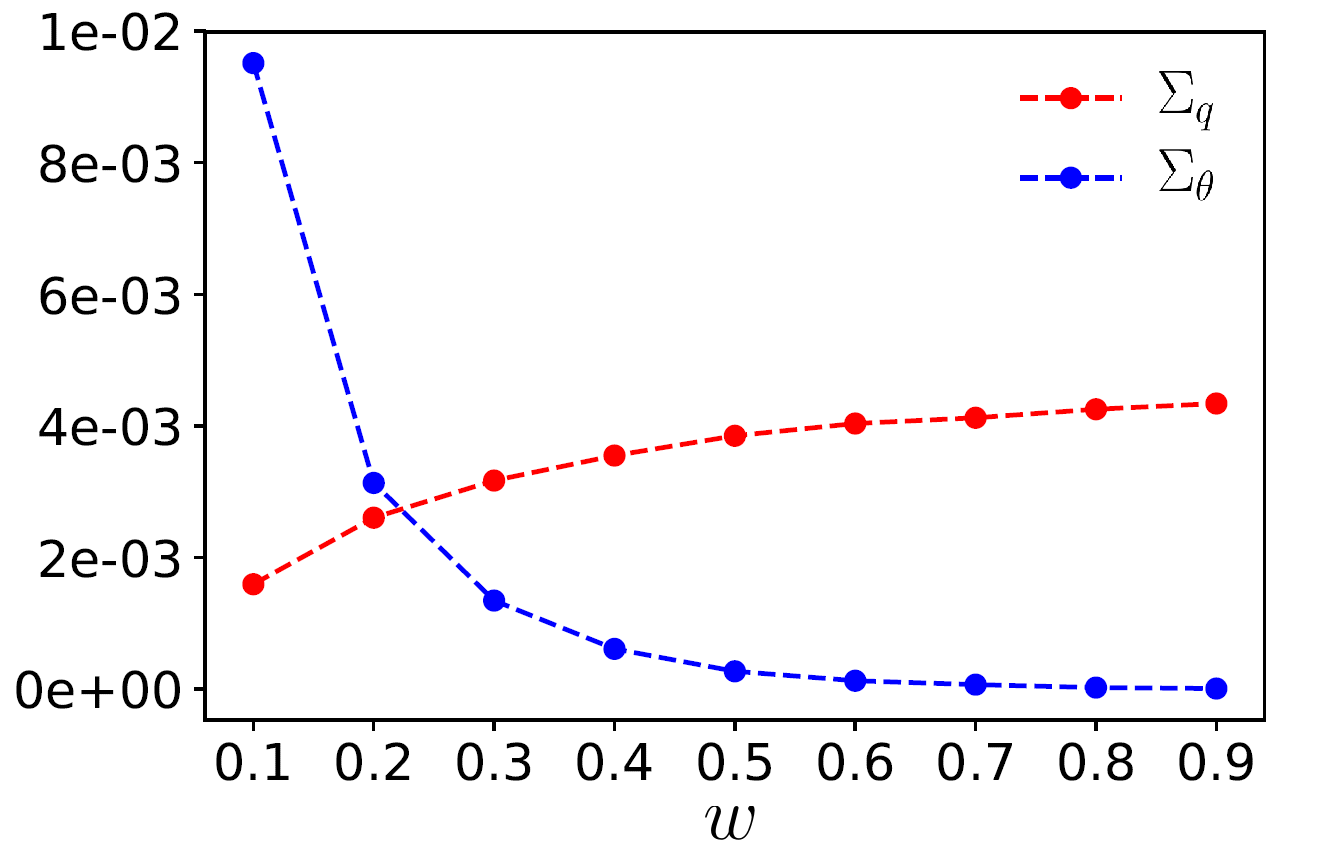}\label{fp10}}\\
    \subfloat[GP-10]{\includegraphics[width=0.4\textwidth]{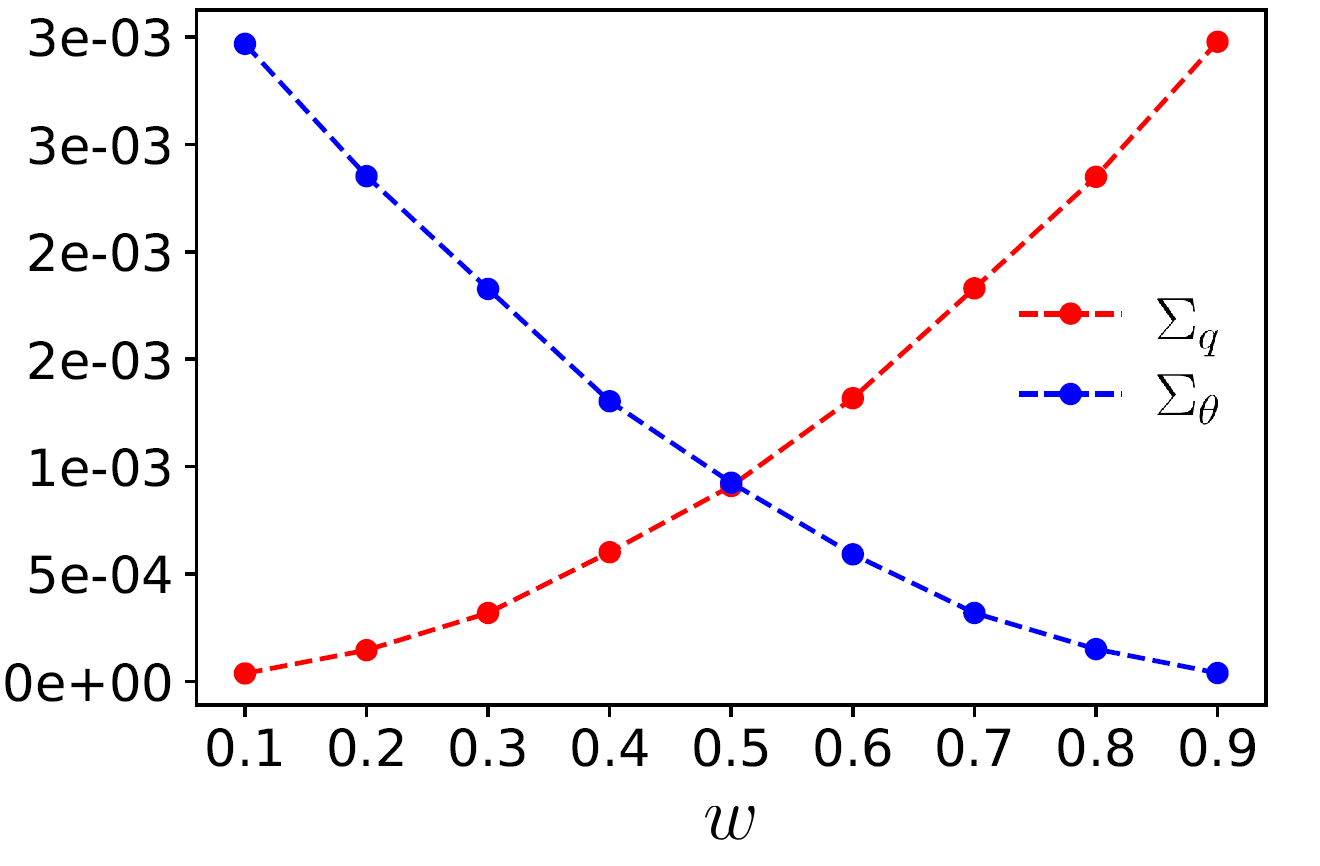}\label{gp10}}
    \hfil
	\subfloat[RCP-30]{\includegraphics[width=0.4\textwidth]{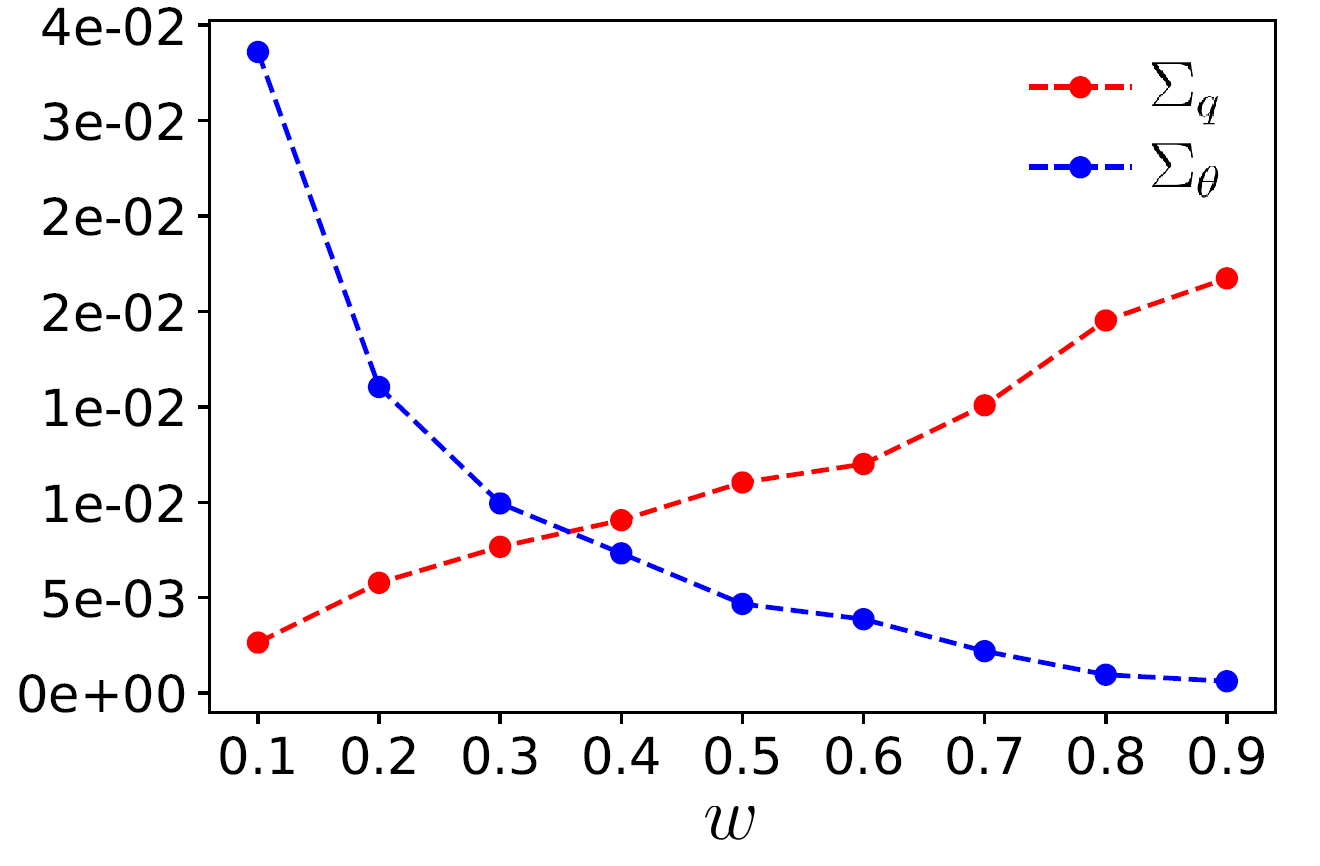}\label{rcp30}}
	\caption{Sensitivity analysis on the preference weight $w$ of the 2D objective function \eqref{eq:obj2d}. For all figures, $\Sigma_q$ represents the total speed deviation defined as $\sum_{i \in \mathcal{A}}(1-q_i)^2$ (in red) and $\Sigma_{\theta}$ represents the total heading deviation defined as $\sum_{i \in \mathcal{A}}\theta_i^2$ (in blue).}
	\label{fig:sa}
\end{figure}

This behavior is confirmed in our numerical experiments. Specifically, we solve the 2D instances CP-8, FP-10, GP-10 and one RCP-30 instance for $w = 0.1, \ldots, 0.9$ in steps of size $0.1$, i.e. for a total of 9 values of $w$ per instance. All instances are solved to optimality using Algorithm \ref{algo:2d} with no MIQCP iterations, i.e. the MIQP returned a global optimal solution for all tests. The change in the total speed deviation $\Sigma_q$ and in the total heading deviation $\Sigma_\theta$ are reported in Figure \ref{fig:sa}. For all four instances tested, we find that increasing $w$ monotonically decreases the total heading deviation and monotonically increases the total speed deviation. Further, we observe that in all cases both the total speed and heading deviations are of similar order of magnitudes. 

This sensitivity analysis shows that using the proposed 2D objective function, the decision-maker can control which manoeuvre is prioritized by scaling up or down the preference weight $w$ accordingly. Higher values of $w$ will minimize the total heading deviation while lower values of $w$ will minimize the total speed deviation. We use $w=0.5$ in our numerical experiments presented in the remaining of the paper. We next discuss numerical results on 2D instances in Section \ref{2dresults} before discussing model performance on 2D+FL instances in Section \ref{3dresults}. 

\subsection{Results on 2D instances}\label{2dresults}

To report the performance of the proposed formulations on 2D ACRPs, we conduct four groups of numerical experiments, one per instance type. We present the results for 12 CP instances ranging from 4 to 15 aircraft in Tables \ref{CP-30} and \ref{CP-15}. Results for FP and GP instances are reported in Tables \ref{FP-30}-\ref{GP-15}, respectively, for 12 instances each with 4 to 15 aircraft per stream. Results for RCP instances are reported in Tables \ref{RCP-30} and \ref{RCP-15} for four instance sizes with 10, 20, 30 and 40 aircraft per group. For each instance size, 100 RCP instances are randomly generated and we report the average performance along with the standard deviation in parenthesis. We compare the performance of the proposed formulations for both the standard and reduced heading control ranges.


Each row in the results tables represents an instance (CP, FP and GP) or a group of instances (RCP). The header of the results tables is presented from left to right. The left-most column, Instance, identifies the instance; $|A|$ is the number of aircraft and $n_c$ is the number of the conflicts. The next three columns summarizes the performance of the pre-processing algorithm \ref{algo:prepro}: $|\PF|/|\P|$ is the proportion of conflict-free aircraft pairs; $|\PI|/|\P|$ is the proportion of non-separable aircraft pairs; and Time is the runtime of Algorithm \ref{algo:prepro} in seconds. The next six columns summarize the performance of \ref{algo:2d} using the disjunctive formulation: LB and UB are the lower and upper bound; Gap is the optimality gap in percent calculated using LB and UB; Time is the total runtime in seconds; $n_i$ is the number of iterations of the \textbf{while} loop in Algorithm \ref{algo:2d}, and $n_t$ represents the proportion of instances that could not be solved within the time limit (10 minutes), i.e. the number of time-outs. The next four columns summarize the performance of the shadow formulation: $\Delta$ UB is the upper bound of the shadow formulation minus that of the disjunctive formulation; and we also report the optimality gap, total runtime and proportion of time-outs. The right-most column, Gain, reports the performance gain for instances solved within the time limit calculated as the runtime of the shadow formulation minus that of the disjunctive formulation in percentage: a positive value indicates that the disjunctive formulation is faster. 

\begin{table*}[t]
\resizebox{1\columnwidth}{!}{%
\begin{tabular}{llll lll lllllll llll l}
\toprule
&&& \multicolumn{3}{l}{Pre-processing} &\multicolumn{6}{l}{Disjunctive} & \multicolumn{4}{l}{Shadow}\\
\cmidrule(l) {4-6} \cmidrule(l){7-12} \cmidrule(l){13-16}
Instance & $|\A|$ &$n_c$ & $\frac{|\PF|}{|\P|}$ & $\frac{|\PI|}{|\P|}$ & Time & LB & UB & Gap & Time & $n_i$ & $n_t$ &$\Delta$UB & Gap & Time & $n_t$ & Gain \\
\midrule
CP-4	&	4	&	6	&	0 & 0	&	0.01	&	6.2E-4	&	6.2E-4	&	0.00	&	0.20	&	0	&	0	&	0.00 &	0.00		&	1.31	&	0	&	84.5	\\
CP-5	&	5	&	10	&	0 & 0	&	0.03	&	1.1E-3	&	1.1E-3	&	0.00	&	0.40	&	0	&	0	&	0.00 &	0.00		&	1.02	&	0	&	60.1	\\
CP-6	&	6	&	15	&	0 & 0	&	0.07	&	1.8E-3	&	1.8E-3	&	0.00	&	0.64	&	0	&	0	&	0.00 &	0.00		&	1.17	&	0	&	45.3	\\
CP-7	&	7	&	21	&	0 & 0	&	0.12	&	2.4E-3	&	2.4E-3	&	0.00	&	0.39	&	0	&	0	&	0.00 &	0.00		&	1.28	&	0	&	69.5	\\
CP-8	&	8	&	28	&	0 & 0	&	0.18	&	3.5E-3	&	3.5E-3	&	0.02	&	3.04	&	0	&	0	&	0.00 &	0.02		&	4.23	&	0	&	28.0	\\
CP-9	&	9	&	36	&	0 & 0	&	0.24	&	4.3E-3	&	4.3E-3	&	0.02	&	7.97	&	0	&	0	&	0.00 &	0.02		&	11.2	&	0	&	29.3	\\
CP-10	&	10	&	45	&	0 & 0	&	0.31	&	5.6E-3	&	5.6E-3	&	0.02	&	72.3	&	0	&	0	&	0.00 &	0.02		&	73.1	&	0	&	1.10	\\
CP-11	&	11	&	55	&	0 & 0	&	0.37	&	6.8E-3	&	6.9E-3	&	1.99	&	600 	&	0	&	100	&	0.00 &	15.6		&	600	    &	100	&	-	\\
CP-12	&	12	&	66	&	0 & 0	&	0.48	&	5.0E-3	&	8.4E-3	&	40.2	&	600	    &	0	&	100	&	0.00 &	37.6		&	600 	&	100	&	-	\\
CP-13	&	13	&	78	&	0 & 0	&	0.63	&	4.4E-3	&	9.9E-3	&	55.2	&	600 	&	0	&	100	&	0.00 &	53.1		&	600 	&	100	&	-	\\
CP-14	&	14	&	91	&	0 & 0	&	0.72	&	3.8E-3	&	1.1E-2	&	66.8	&	600 	&	0	&	100	&	0.00 &	63.0		&	600 	&	100	&	-	\\
CP-15	&	15	&	105	&	0 & 0	&	0.89	&	3.6E-3	&	1.5E-2	&	73.5	&	600 	&	0	&	100	&	0.00 &	74.1		&	600 	&	100	&	-	\\
\bottomrule
\end{tabular}
}
\caption{Results for 2D CP instances with a speed control range of $[-6\%,+3\%]$ and a heading control range of $[-30^\circ,+30^\circ]$. Times (Time) are reported in seconds. The proportions of conflict-free ($\frac{|\PF|}{|\P|}$) and non-separable ($\frac{|\PI|}{|\P|}$) pairs, optimality gaps (Gap), proportions of time-outs ($n_t$) and the performance gain (Gain) are reported in \%.}
\label{CP-30}
\end{table*}

\begin{table*}[t]
\resizebox{1\columnwidth}{!}{%
\begin{tabular}{llll lll lllllll llll l}
\toprule
&&& \multicolumn{3}{l}{Pre-processing} &\multicolumn{6}{l}{Disjunctive} & \multicolumn{4}{l}{Shadow}\\
\cmidrule(l) {4-6} \cmidrule(l){7-12} \cmidrule(l){13-16}
Instance & $|\A|$ &$n_c$ & $\frac{|\PF|}{|\P|}$ & $\frac{|\PI|}{|\P|}$ & Time & LB & UB & Gap & Time & $n_i$ & $n_t$ &$\Delta$UB & Gap & Time & $n_t$ & Gain \\
\midrule
CP-4	&	4	&	6	&	0   &	0   & 0.01 	&	6.2E-4	&	6.2E-4	&	0.00	&	0.25		&	0	&	0	&   0.00	&	0.00	&	1.27	&	0	&	80.2	\\
CP-5	&	5	&	10	&	0   &	0  & 0.03	&	1.1E-3	&	1.1E-3	&	0.00	&	0.03		&	0	&	0	&	0.00	&	0.00	&	1.16	&	0	&	97.3	\\
CP-6	&	6	&	15	&	0   &	0  & 0.11	&	1.8E-3	&	1.8E-3	&	0.00	&	0.23		&	0	&	0	&	0.00	&	0.00	&	1.24	&	0	&	81.4	\\
CP-7	&	7	&	21	&	0   &	0  & 0.18	&	2.4E-3	&	2.4E-3	&	0.00	&	0.19		&	0	&	0	&	0.00	&	0.01	&	1.20	&	0	&	84.4	\\
CP-8	&	8	&	28	&	0   &	0  & 0.25	&	3.5E-3	&	3.5E-3	&	0.02	&	3.17		&	0	&	0	&	0.00	&	0.00	&	3.72	&	0	&	14.7	\\
CP-9	&	9	&	36	&	0   &	0  & 0.31	&	4.3E-3	&	4.3E-3	&	0.02	&	8.11		&	0	&	0	&	0.00	&	0.02	&	8.80	&	0	&	7.80	\\
CP-10	&	10	&	45	&	0   &	0  & 0.42	&	5.6E-3	&	5.6E-3	&	0.02	&	75.1		&	0	&	0	&	0.00	&	0.02	&	85.1	&	0	&	11.7	\\
CP-11	&	11	&	55	&	0   &	0  & 0.52	&	6.6E-3	&	6.8E-3	&	5.25     &	600 		&	0	&	100	&	0.00	&	10.2	&	600 	&	100	&	-	\\
CP-12	&	12	&	66	&	0   &	0  & 0.62	&	4.9E-3	&	8.4E-3	&	40.7	&	600	    	&   0	&	100	&	0.00	&	39.9	&	600	    &	100	&	-	\\
CP-13	&	13	&	78	&	0   &	0  & 0.72	&	4.4E-3	&	9.9E-3	&	55.5	&	600	    	&	0	&	100	&	0.00	&	50.1	&	600	    &	100	&	-	\\
CP-14	&	14	&	91	&	0   &	0  & 0.80	&	3.8E-3	&	1.1E-2	&	66.5	&	600 		&	0	&	100	&	0.00	&	65.7	&	600	    &	100	&	-	\\
CP-15	&	15	&	105	&	0   &	0  & 0.88	&	3.7E-3	&	1.4E-2	&	72.9	&	600 		&	0	&	100	&	0.00	&	70.6	&	600 	&	100	&	-	\\
\bottomrule
\end{tabular}
}
\caption{Results for 2D CP instances with a speed control range of $[-6\%,+3\%]$ and a heading control range of $[-15^\circ,+15^\circ]$. Times (Time) are reported in seconds. The proportions of conflict-free ($\frac{|\PF|}{|\P|}$) and non-separable ($\frac{|\PI|}{|\P|}$) pairs, optimality gaps (Gap), proportions of time-outs ($n_t$) and the performance gain (Gain) are reported in \%.}
\label{CP-15}
\end{table*}

The implementation of the pre-processing procedure (Algorithm \ref{algo:prepro}) across CP, FP and GP instances with both the standard heading control range of $[-30^\circ,+30^\circ]$ and the reduced range of $[-15^\circ,+15^\circ]$ reveal that no aircraft pair can be eliminated (i.e. is conflict-free) or is non-separable. We also note that the runtime of Algorithm \ref{algo:prepro} on all these instances grows quadratically with the number of aircraft and evolve under from 0.01 s to 1.23 s.  

The experiments on the CP instances (Tables \ref{CP-30} and \ref{CP-15}) show that, as expected, the upper bound (UB) and the lower bound (LB) increase with the number of aircraft. Using a heading control range of $[-30^\circ,+30^\circ]$ (Table \ref{CP-30}), both the disjunctive and shadow formulations are able to solve CP instances with 4 to 10 aircraft within the available time limit. In both cases, runtime increases exponentially with the number of aircraft. For CP instances with 11 to 15 aircraft, the MIQP iteration of Algorithm \ref{algo:2d} is unable to converge. 
In terms of performance gain, the disjunctive formulation outperforms the shadow formulation for instances solved but both formulations return comparable optimality gaps for instances that timed-out. Notably, the optimality gap of CP-15 is 73.5\% and 74.1\% using the disjunctive and shadow formulations, respectively, which highlights the difficulty of CP instances. Using a reduced heading control range of $[-15^\circ,+15^\circ]$ (Table \ref{CP-15}), yields a comparable performance, although the proposed formulations are able to achieve slightly better UB values for CP-11 and CP-15, suggesting that the reduction of the solution space may help performance. We also find that the same UB values are obtained for both heading control ranges for instances solved to optimality.


The numerical results obtained for FP instances (Tables \ref{FP-30} and \ref{FP-15}) for both heading control ranges reveal that all corresponding 24 ACRPs can be solved to optimality within the available time limit by the disjunctive formulation. The shadow fails to solve FP-15 with a reduced heading control range (Table \ref{FP-15}). Both formulations find the same optimal objective value, i.e. $\Delta$ UB = 0 for all FP instances tested. The disjunctive formulation is able to solve all FP instances in less than a minute, which contrasts with the performance of the shadow formulation which requires significantly more time to solve FP instances with 9 or more aircraft per stream. This corresponds to an average performance gain of 90.9\% across all FP instances solved to optimality by both formulations. 

The outcome of the numerical experiments for GP instances (Tables \ref{GP-30} and \ref{GP-15}) reveal a similar trend. The disjunctive formulation is able to solve all GP instances with runtimes of 125 s and 132 s for GP-15, which corresponds to a total of 60 aircraft, with a standard and a reduced heading control range, respectively. The shadow formulation is able to solve all but three GP instances across all 24 GPs tested. The average performance gain for GP instances solved to optimality by both formulations is 68.0\%. We also observe that reducing the heading control range from the standard $[-30^\circ,+30^\circ]$ to $[-15^\circ,+15^\circ]$ does not affect the optimal solution of these problems and the shadow formulation is consistently able to find the same best UB as the disjunctive formulation. 

The experiments performed on RCP instances reveal that the pre-processing procedure (Algorithm \ref{algo:prepro}) is able to eliminate approximately 8\% of aircraft pairs when using a reduced heading control range of $[-15^\circ,+15^\circ]$ (Table \ref{RCP-15}), whereas no aircraft pairs are conflict-free using a standard heading control range of $[-30^\circ,+30^\circ]$ (Table \ref{RCP-30}). 

The implementation of Algorithm \ref{algo:2d} on RCP instances (Tables \ref{RCP-30} and \ref{RCP-15}) reveals that while all 10- and 20-aircraft RCP instances can be solved via the MIQP iteration, 30- and 40-aircraft RCP instances may require additional MIQCP iterations. Using the standard heading control range (Table \ref{RCP-30}), on average, RCP-10 instances can be solved in 0.05 s and 0.14 s using the disjunctive and shadow formulations, respectively. RCP-20 instances require 0.26 s and 2.52 s, on average, using the disjunctive and shadow formulations, respectively. The performance of the proposed formulations on RCP-10 and RCP-20 using a reduced heading control range (Table \ref{RCP-15}) is of similar order of magnitude due to the relatively low average number of conflicts, 3.1 and 13.1, respectively, in these problems. 

RCP-30 and RCP-40 instances have on average 32.9 and 59.3 conflicts, respectively, and present considerable computing challenges, notably the latter. The disjunctive formulation requires an average of 0.4 MIQCP iterations for RCP-30 instances using a standard heading control range of $[-30^\circ,+30^\circ]$ (Table \ref{RCP-30}) compared to 1.4 when using a reduced heading control range of $[-15^\circ,+15^\circ]$ (Table \ref{RCP-15}). These figures increase to 0.8 and 1.6 for RCP-40 instances. Overall, we observe that reducing the heading control range tends to improve the performance of the proposed formulations while retaining comparable optimal solutions, as indicated by the similar UB values obtained. Using a standard heading control range (Table \ref{RCP-30}), the disjunctive is able to solve all but 3\% of the RCP-30 instances whereas the shadow formulation times-out on 20\% of these instances. The performance gain is 41.7\% on RCP-30 solved to optimality by both formulations. For RCP-40 instances, the disjunctive formulation is able to 28\% of the problems (72\% of time-outs) and the shadow formulation is unable to solve any of the 100 instances tested within the available time limit. Using a reduced heading control range (Table \ref{RCP-15}), the disjunctive formulation times-out on 3\% of the RCP-30 instances and on 59\% of the RCP-40 instances. In contrast, the shadow formulation times-out on 10\% of the RCP-30 instances and is unable to solve any of the RCP-40 instances. We also observe that the shadow formulation may fail to find an UB as competitive as that found by the disjunctive formulation as noted by the average $\Delta$UB value of 0.02 in Tables \ref{RCP-30} and \ref{RCP-15}.

\begin{table*}[t]
\resizebox{1\columnwidth}{!}{%
\begin{tabular}{llll lll lllllll llll l}
\toprule
&&& \multicolumn{3}{l}{Pre-processing} &\multicolumn{6}{l}{Disjunctive} & \multicolumn{4}{l}{Shadow}\\
\cmidrule(l) {4-6} \cmidrule(l){7-12} \cmidrule(l){13-16}
Instance & $|\A|$ &$n_c$ & $\frac{|\PF|}{|\P|}$ & $\frac{|\PI|}{|\P|}$ & Time & LB & UB & Gap & Time & $n_i$ & $n_t$ &$\Delta$UB & Gap & Time & $n_t$ & Gain \\
\midrule
FP-4	&	8	&	4	& 0 & 0 & 0.17 & 	8.1E-4	&	8.2E-4	&	0.10	&	0.05	&	0	&	0	&	0.00	&	0.10	&	0.26	&	0	&	80.7	\\
FP-5	&	10	&	5	&0 & 0 & 0.28 & 	1.1E-3	&	1.1E-3	&	0.03	&	0.05	&	0	&	0	&	0.00	&	0.08	&	4.17	&	0	&	98.8	\\
FP-6	&	12	&	6	&0 & 0 & 0.36 & 	1.5E-3	&	1.5E-3	&	0.05	&	0.16	&	0	&	0	&	0.00	&	0.72	&	6.00	&	0	&	97.3	\\
FP-7	&	14	&	7	&0 & 0 & 0.40 & 	2.1E-3  &	2.1E-3	&	0.05	&	0.44	&	0	&	0	&	0.00	&	0.32	&	6.80	&	0	&	93.5	\\
FP-8	&	16	&	8	&0 & 0 & 0.45 & 	2.8E-3	&	2.7E-3	&	0.04	&	2.31	&	0	&	0	&	0.00	&	0.04	&	7.66	&	0	&	69.8	\\
FP-9	&	18	&	9	&0 & 0 & 0.52 & 	3.7E-3	&	3.7E-3	&	0.03	&	2.28	&	0	&	0	&	0.00	&	0.15	&	126	&	0	&	98.1	\\
FP-10	&	20	&	10	&0 & 0 & 0.61 & 	5.2E-3	&	5.2E-3	&	0.02	&	8.41	&	0	&	0	&	0.00	&	0.23	&	138	&	0	&	93.9	\\
FP-11	&	22	&	12	&0 & 0 & 0.71 &	5.8E-3	&	5.8E-3	&	0.14	&	12.3	&	0	&	0	&	0.00	&	0.05	&	126	&	0	&	90.2	\\
FP-12	&	24	&	14	&0 & 0 & 0.82 &	6.6E-3	&	6.6E-3	&	0.03	&	16.2	&	0	&	0	&	0.00	&	0.25	&	261	&	0	&	93.7	\\
FP-13	&	26	&	16	&0 & 0 & 0.94 &	7.4E-3	&	4.5E-3	&	0.12	&	21.4	&	0	&	0	&	0.00	&	0.12	&	342	&	0	&	93.7	\\
FP-14	&	25	&	18	&0 & 0 & 1.05 &	2.7E-2	&	2.7E-2	&	0.14	&	25.3	&	0	&	0	&	0.00	&	0.16	&	465	&	0	&	94.5	\\
FP-15	&	30	&	22	&0 & 0 & 1.15 &	3.5E-2	&	3.5E-2	&	0.13	&	42.2	&	0	&	0	&	0.00	&	0.24	&	576	&	0	&	92.6	\\
\bottomrule
\end{tabular}
}
\caption{Results for 2D FP instances with a speed control range of $[-6\%,+3\%]$ and a heading control range of $[-30^\circ,+30^\circ]$. Times (Time) are reported in seconds. The proportions of conflict-free ($\frac{|\PF|}{|\P|}$) and non-separable ($\frac{|\PI|}{|\P|}$) pairs, optimality gaps (Gap), proportions of time-outs ($n_t$) and the performance gain (Gain) are reported in \%.}
\label{FP-30}
\end{table*}

\begin{table*}[t]
\resizebox{1\columnwidth}{!}{%
\begin{tabular}{llll lll lllllll llll l}
\toprule
&&& \multicolumn{3}{l}{Pre-processing} &\multicolumn{6}{l}{Disjunctive} & \multicolumn{4}{l}{Shadow}\\
\cmidrule(l) {4-6} \cmidrule(l){7-12} \cmidrule(l){13-16}
Instance & $|\A|$ &$n_c$ & $\frac{|\PF|}{|\P|}$ & $\frac{|\PI|}{|\P|}$ & Time & LB & UB & Gap & Time & $n_i$ & $n_t$ &$\Delta$UB & Gap & Time & $n_t$ & Gain \\
\midrule
FP-4	&	8	&	4	& 0 & 0 & 0.09 & 	8.1E-4	&	8.2E-4	&	0.11	&	0.06	&	0	&	0	&	0.00	&	0.09	&	0.20	&	0	&	70.0	\\
FP-5	&	10	&	5	&0 & 0 & 0.22 & 	1.1E-3	&	1.1E-3	&	0.03	&	0.07	&	0	&	0	&	0.00	&	0.10	&	4.57	&	0	&	98.4	\\
FP-6	&	12	&	6	&0 & 0 & 0.29 & 	1.5E-3	&	1.5E-3	&	0.03	&	0.21	&	0	&	0	&	0.00	&	0.67	&	6.23	&	0	&	96.6	\\
FP-7	&	14	&	7	&0 & 0 & 0.41 & 	2.1E-3  &	2.1E-3	&	0.06	&	0.47	&	0	&	0	&	0.00	&	0.30	&	6.81	&	0	&	93.1	\\
FP-8	&	16	&	8	&0 & 0 & 0.45 & 	2.8E-3	&	2.7E-3	&	0.04	&	2.56	&	0	&	0	&	0.00	&	0.41	&	7.66	&	0	&	66.5	\\
FP-9	&	18	&	9	&0 & 0 & 0.62 & 	3.7E-3	&	3.7E-3	&	0.03	&	2.76	&	0	&	0	&	0.00	&	0.15	&	156	&	0	&	98.2	\\
FP-10	&	20	&	10	&0 & 0 & 0.72 & 	5.2E-3	&	5.2E-3	&	0.02	&	9.52	&	0	&	0	&	0.00	&	0.24	&	143	&	0	&	93.3	\\
FP-11	&	22	&	12	&0 & 0 & 0.85 &	5.8E-3	&	5.8E-3	&	0.04	&	10.3	&	0	&	0	&	0.00	&	0.04	&	176	&	0	&	94.4	\\
FP-12	&	24	&	14	&0 & 0 & 0.92 &	6.6E-3	&	6.6E-3	&	0.13	&	12.2	&	0	&	0	&	0.00	&	0.14	&	226	&	0	&	94.6	\\
FP-13	&	26	&	16	&0 & 0 & 0.99 &	7.4E-3	&	4.5E-3	&	0.12	&	18.4	&	0	&	0	&	0.00	&	0.01	&	338	&	0	&	94.5	\\
FP-14	&	25	&	18	&0 & 0 & 1.02 &	2.7E-2	&	2.7E-2	&	0.04	&	22.3	&	0	&	0	&	0.00	&	0.04	&	476	&	0	&	95.3	\\
FP-15	&	30	&	22	&0 & 0 & 1.23 &	3.5E-2	&	3.5E-2	&	0.23	&	32.2	&	0	&	0	&	0.00	&	0.14	&	600	&	100	&	-	\\
\bottomrule
\end{tabular}
}
\caption{Results for 2D FP instances with a speed control range of $[-6\%,+3\%]$ and a heading control range of $[-15^\circ,+15^\circ]$. Times (Time) are reported in seconds. The proportions of conflict-free ($\frac{|\PF|}{|\P|}$) and non-separable ($\frac{|\PI|}{|\P|}$) pairs, optimality gaps (Gap), proportions of time-outs ($n_t$) and the performance gain (Gain) are reported in \%.}
\label{FP-15}
\end{table*}

\begin{table*}[t]
\resizebox{1\columnwidth}{!}{%
\begin{tabular}{llll lll lllllll llll l}
\toprule
&&& \multicolumn{3}{l}{Pre-processing} &\multicolumn{6}{l}{Disjunctive} & \multicolumn{4}{l}{Shadow}\\
\cmidrule(l) {4-6} \cmidrule(l){7-12} \cmidrule(l){13-16}
Instance & $|\A|$ &$n_c$ & $\frac{|\PF|}{|\P|}$ & $\frac{|\PI|}{|\P|}$ & Time & LB & UB & Gap & Time & $n_i$ & $n_t$ &$\Delta$UB & Gap & Time & $n_t$ & Gain \\
\midrule
GP-4	&	16	&	12	& 0 & 0 & 0.03 &	9.2E-4	&	9.2E-4	&	0.01	&	0.16	&	0	&	0	&	0.00	&	0.11	&	1.22	&	0	&	86.8	\\
GP-5	&	20	&	16	&0 & 0 & 0.10 &	1.3E-3	&	1.3E-3	&	0.02	&	0.98	&	0	&	0	&	0.00	&	0.25	&	2.34	&	0	&	58.1	\\
GP-6	&	24	&	20	&0 & 0 & 0.21 &	1.8E-3	&	1.8E-3	&	0.04	&	2.85	&	0	&	0	&	0.00	&	0.26	&	12.1	&	0	&	76.4	\\
GP-7	&	28	&	24	&0 & 0 & 0.32&	2.4E-3	&	2.4E-3	&	0.03	&	11.3	&	0	&	0	&	0.00	&	0.33	&	14.4	&	0	&	21.5	\\
GP-8	&	32	&	28	&0 & 0 & 0.44 &	3.2E-3	&	3.2E-3	&	0.02	&	20.7	&	0	&	0	&	0.00	&	0.12	&	59.5	&	0	&	65.2	\\
GP-9	&	36	&	32	&0 & 0 & 0.58 &	4.3E-3	&	4.3E-3	&	0.04	&	46.2	&	0	&	0	&	0.00	&	0.24	&	165	&	0	&  72.0	\\
GP-10	&	40	&	36	&0 & 0 & 0.69 &	6.1E-3	&	6.1E-3	&	0.02	&	47.5	&	0	&	0	&	0.00	&	0.82	&	194	&	0	&	75.5	\\
GP-11	&	44	&	40	&0 & 0 & 0.72 &	7.9E-3	&	7.9E-3	&	0.24	&	52.3	&	0	&	0	&	0.00	&	0.12	&	225	&	0	&	77.6	\\
GP-12	&	48	&	44	&0 & 0 & 0.89 &	8.3E-3	&	8.3E-3	&	0.31	&	61.2	&	0	&	0	&	0.00	&	0.13	&	383	&	0	&	83.7	\\
GP-13	&	52	&	48	&0 & 0 & 0.98 &	9.5E-3	&	9.5E-3	&	0.12	&	75.3	&	0	&	0	&	0.00	&	0.21	&	476	&	0	&   83.5	 \\
GP-14	&	56	&	52	&0 & 0 & 1.09 &	3.4E-2	&	3.4E-2	&	0.56	&	92.8	&	0	&	0	&	0.00	&	0.24	&	600	&	100	&	-	\\
GP-15	&	60	&	54	&0 & 0 & 1.24 &	3.8E-2	&	3.8E-2	&	0.45	&	125	&	0	&	0	&	0.00	&	0.25	&	600	&	100	&	-	\\
\bottomrule
\end{tabular}
}
\caption{Results for 2D GP instances with a speed control range of $[-6\%,+3\%]$ and a heading control range of $[-30^\circ,+30^\circ]$. Times (Time) are reported in seconds. The proportions of conflict-free ($\frac{|\PF|}{|\P|}$) and non-separable ($\frac{|\PI|}{|\P|}$) pairs, optimality gaps (Gap), proportions of time-outs ($n_t$) and the performance gain (Gain) are reported in \%.}
\label{GP-30}
\end{table*}

\begin{table*}[t]
\resizebox{1\columnwidth}{!}{%
\begin{tabular}{llll lll lllllll llll l}
\toprule
&&& \multicolumn{3}{l}{Pre-processing} &\multicolumn{6}{l}{Disjunctive} & \multicolumn{4}{l}{Shadow}\\
\cmidrule(l) {4-6} \cmidrule(l){7-12} \cmidrule(l){13-16}
Instance & $|\A|$ &$n_c$ & $\frac{|\PF|}{|\P|}$ & $\frac{|\PI|}{|\P|}$ & Time & LB & UB & Gap & Time & $n_i$ & $n_t$ &$\Delta$UB & Gap & Time & $n_t$ & Gain \\
\midrule
GP-4	&	16	&	12	& 0 & 0 & 0.05 &	9.2E-4	&	9.2E-4	&	0.03	&	0.24	&	0	&	0	&	0.00	&	0.18	&	1.58	&	0	&	84.6	\\
GP-5	&	20	&	16	&0 & 0 & 0.05 &	1.3E-3	&	1.3E-3	&	0.02	&	1.81	&	0	&	0	&	0.00	&	0.21	&	2.16	&	0	&	16.2	\\
GP-6	&	24	&	20	&0 & 0 & 0.15 &	1.8E-3	&	1.8E-3	&	0.05	&	3.91	&	0	&	0	&	0.00	&	0.24	&	11.3	&	0	&	65.3	\\
GP-7	&	28	&	24	&0 & 0 & 0.24 &	2.4E-3	&	2.4E-3	&	0.03	&	11.2	&	0	&	0	&	0.00	&	0.65	&	18.6	&	0	&	39.7	\\
GP-8	&	32	&	28	&0 & 0 & 0.34 &	3.2E-3	&	3.2E-3	&	0.02	&	21.5	&	0	&	0	&	0.00	&	0.12	&	59.5	&	0	&	63.8	\\
GP-9	&	36	&	32	&0 & 0 & 0.48 &	4.3E-3	&	4.3E-3	&	0.05	&	46.5	&	0	&	0	&	0.00	&	0.24	&	136	&	0	&	65.8	\\
GP-10	&	40	&	36	&0 & 0 & 0.56 &	6.1E-3	&	6.1E-3	&	0.04	&	49.2	&	0	&	0	&	0.00	&	0.89	&	188	&	0	&	73.8	\\
GP-11	&	44	&	40	&0 & 0 & 0.67 &	7.9E-3	&	7.9E-3	&	0.14	&	50.3	&	0	&	0	&	0.00	&	0.32	&	205	&	0	&	75.4	\\
GP-12	&	48	&	44	&0 & 0 & 0.78 &	8.3E-3	&	8.3E-3	&	0.23	&	62.2	&	0	&	0	&	0.00	&	0.43	&	352	&	0	&	82.2	\\
GP-13	&	52	&	48	&0 & 0 & 0.81 &	9.5E-3	&	9.5E-3	&	0.22	&	78.3	&	0	&	0	&	0.00	&	0.31	&	421	&	0	&   81.4	 \\
GP-14	&	56	&	52	&0 & 0 & 0.95 &	3.4E-2	&	3.4E-2	&	0.24	&	102	&	0	&	0	&	0.00	&	0.42	&	525	&	0	&	80.5	\\
GP-15	&	60	&	54	&0 & 0 & 1.03 &	3.8E-2	&	3.8E-2	&	0.33	&	132	&	0	&	0	&	0.00	&	0.51	&	600	&	100	&	-	\\
\bottomrule
\end{tabular}
}
\caption{Results for 2D GP instances with a speed control range of $[-6\%,+3\%]$ and a heading control range of $[-15^\circ,+15^\circ]$. Times (Time) are reported in seconds. The proportions of conflict-free ($\frac{|\PF|}{|\P|}$) and non-separable ($\frac{|\PI|}{|\P|}$) pairs, optimality gaps (Gap), proportions of time-outs ($n_t$) and the performance gain (Gain) are reported in \%.}
\label{GP-15}
\end{table*}

\newcommand{\flrcpones}{
\resizebox{1\columnwidth}{!}{%
\begin{tabular}{lll lll lllllll}
\toprule
&&& \multicolumn{3}{l}{Pre-processing} & \multicolumn{6}{l}{Disjunctive}\\
\cmidrule(l){4-6}\cmidrule(l){7-12} 
Instance & $|\A|$ &$n_c$ & $\frac{|\PF|}{|\P|}$ & $\frac{|\PI|}{|\P|}$ & Time & LB & UB & Gap & Time & $n_i$ & $n_t$ \\
\midrule
RCP-10 & 10 & 3.10 (1.6) & 0 & 0 & 0.12  & 2.2E-4 (2E-4) & 2.2E-4 (2E-4) & 0.00 (0.0) & 0.05 (0.01) & 0.0 (0.0) & 0 \\
RCP-20 & 20 & 13.1 (3.5) & 0 & 0& 0.23  & 1.7E-3 (9E-4) & 1.7E-3 (9E-4) & 0.00 (0.0) & 0.26 (0.10) & 0.0 (0.0) & 0 \\
RCP-30 & 30 & 32.9 (5.6) & 0 & 0& 0.45  & 7.1E-3 (2E-3) & 7.1E-3 (2E-3) & 0.17 (0.7) & 135 (239) &  0.4 (0.6) & 3 \\
RCP-40 & 40 & 59.3 (7.1) & 0 & 0& 0.60  & 1.8E-2 (5E-3) & 2.4E-2 (1E-2) & 15.0 (25) & 516 (194) &  0.8 (0.5)  &  72 \\
\bottomrule
\end{tabular}
}
}
\newcommand{\flrcptwos}{
\resizebox{0.5\columnwidth}{!}{%
\begin{tabular}{l llll l}
\toprule
& \multicolumn{4}{l}{Shadow}\\
\cmidrule(l){2-5}
Instance & $\Delta$UB & Gap & Time & $n_t$ & Gain \\
\midrule
RCP-10 & 0.00 (0.0) & 0.02 (0.1) & 0.14 (0.1) & 0 & 71.8 \\
RCP-20 & 0.00 (0.0) & 0.01 (0.0) & 2.52 (0.6) & 0 &  89.6 \\
RCP-30 & 0.00 (0.0) & 2.36 (4.9) & 231 (218) & 20 & 41.7 \\
RCP-40 & 0.02 (0.1) & 29.8 (25) & 600 (0.0) &  100 & - \\			
\bottomrule
\end{tabular}
}
}
\begin{figure}[t]
\centering
\subfloat{\flrcpones}\\%
\subfloat{\flrcptwos}
\captionof{table}{Results for 2D RCP instances with a speed control range of $[-6\%,+3\%]$ and a heading control range of $[-30^\circ,+30^\circ]$. Times (Time) are reported in seconds. The proportions of conflict-free ($\frac{|\PF|}{|\P|}$) and non-separable ($\frac{|\PI|}{|\P|}$) pairs, optimality gaps (Gap), proportions of time-outs ($n_t$) and the performance gain (Gain) are reported in \%.}%
\label{RCP-30}%
\end{figure}

\newcommand{\flrcpone}{
\resizebox{1\columnwidth}{!}{%
\begin{tabular}{lll lll lllllll}
\toprule
&&& \multicolumn{3}{l}{Pre-processing} & \multicolumn{6}{l}{Disjunctive}\\
\cmidrule(l){4-6}\cmidrule(l){7-12} 
Instance & $|\A|$ &$n_c$ & $\frac{|\PF|}{|\P|}$ & $\frac{|\PI|}{|\P|}$ & Time & LB & UB & Gap & Time & $n_i$ & $n_t$ \\
\midrule
RCP-10 & 10 & 3.10 (1.6)  & 8.2 &0 & 0.33  & 2.2E-4 (2E-4) & 2.2E-4 (2E-4) & 0.03 (0.2) & 0.04 (0.0) & 0.0 (0.0) & 0 \\
RCP-20 & 20 & 13.1 (3.5) & 7.7 &0 & 0.51  & 1.7E-3 (9E-4) & 1.7E-3 (9E-4) & 0.01 (0.0) & 0.24 (0.1) & 0.0 (0.0) & 0\\
RCP-30 & 30 & 32.9 (5.6) & 7.7 &0 & 0.62  & 7.2E-3 (2E-3) & 7.2E-3 (2E-3) & 0.01 (0.1) & 66.3 (135) &  1.4 (0.8) & 3 \\
RCP-40 & 40 & 59.3 (7.1) & 7.9 & 0 & 0.75  & 1.8E-2 (5E-3) & 2.2E-2 (1E-2) & 13.1 (20.7) & 389 (261) &  1.6 (0.5) &  59 \\
\bottomrule
\end{tabular}
}
}

\newcommand{\flrcptwo}{
\resizebox{0.5\columnwidth}{!}{%
\begin{tabular}{l llll l}
\toprule
& \multicolumn{4}{l}{Shadow}\\
\cmidrule(l){2-5}
Instance & $\Delta$UB & Gap & Time & $n_t$ & Gain \\
\midrule
RCP-10 & 0.00 (0.0) & 0.03 (0.1) & 0.11 (0.1) & 0 & 62.6 \\
RCP-20 & 0.00 (0.0) & 0.01 (0.0) & 1.39 (1.1) & 0 &  81.3 \\
RCP-30 & 0.00 (0.0) & 0.89 (2.9) & 167 (181) & 10 & 60.5 \\
RCP-40 & 0.02 (0.1) & 25.8 (24) & 600 (0.0) &  100 & - \\			
\bottomrule
\end{tabular}
}
}

\begin{figure}[t]
\centering
\subfloat{\flrcpone}\\%
\subfloat{\flrcptwo}
\captionof{table}{Results for 2D RCP instances with a speed control range of $[-6\%,+3\%]$ and a heading control range of $[-15^\circ,+15^\circ]$. Times (Time) are reported in seconds. The proportions of conflict-free ($\frac{|\PF|}{|\P|}$) and non-separable ($\frac{|\PI|}{|\P|}$) pairs, optimality gaps (Gap), proportions of time-outs ($n_t$) and the performance gain (Gain) are reported in \%.}%
\label{RCP-15}%
\end{figure}

\subsection{Results on 2D+FL instances}\label{3dresults}

Results on 2D+FL instances are reported using similar tables as in Section \ref{2dresults}. Instead of reporting the number of aircraft and the number of conflicts, we report the average number of aircraft and conflicts per FL, i.e. $\frac{|\mathcal{A}|}{|\Z|}$, $\frac{n_c}{|\Z|}$ respectively. In addition, we include a section corresponding to the FL assignment formulation. In all our numerical experiments on 2D+FL instances, only a single pass through the \textbf{while} loop of Algorithm \ref{algo:3d} is required. Hence, we only report Obj which is the objective function value of MILP \ref{mod:fa}; and Time which is the corresponding computing runtime in seconds. For numerical experiment on 2D+FL problems, we focus on RCP instances named RCP-N-Z where N is the number of aircraft and Z is the number of FLs. We consider three number of aircraft: 50, 100 and 150; and two numbers of FLs: 3 and 5. In terms of conflict density, for the same number of aircraft, instances with 3 FLs have a greater number of conflicts compared to instances with 5 FLs and are more computationally challenging. The results are reported in Table \ref{FL-RCP-30} for the standard heading control range, and in Table \ref{FL-RCP-15} for the reduced heading control range.\\

Using the standard heading control range $[-30^\circ,+30^\circ]$ (Table \ref{FL-RCP-30}), the numerical experiments reveal that among all six groups of 2D+FL instances, only some RCP-150-3 and RCP-150-5 instances may require aircraft to change FLs. That is, for all other 2D+FL instances, all conflicts can be resolved using 2D trajectory control only and thus aircraft do not require to perform any FL change. This can be explained by observing that for instances RCP-150-3, the average number of aircraft per FL is 50 which corresponds to a denser aircraft configurations compared to the 2D RCP-40 instances, which all admit feasible solutions. For RCP-150-5, even though the average density per FL is 30, some instances may have denser FL requiring FL separation. Among all 100 RCP-150-3 instances, 35\% of the instances require a FL change. The maximum number of non-separable pairs ($|\PI|$) is 2 and the maximum objective value of Model \ref{mod:fa} is 1, indicating that only a single aircraft deviated from its initial FL. From a computational standpoint the average runtime of the Model \ref{mod:fa} on RCP-150-3 instances is 4.56 s. In comparison, only a single RCP-150-5 instance required a FL change for a single aircraft. For RCP-50-3, the disjunctive and shadow formulations are able to solve all instances in an average time of 0.88 s and 10.2 s, respectively. For RCP-100-3, which corresponds to an average number of 33.3 aircraft per FL, the disjunctive formulation solves all instances in an average time of 3.3 minutes, whereas the shadow formulation times-out on 46\% of the instances and requires an average time of 9.5 minutes for the instances solved. The denser RCP-150-3 instances, with an average of 50 aircraft per FL, present substantial computational challenges. The disjunctive formulation is able to solve only 40\% of these instances (60\% of time-outs) and an average of 2.8 MIQCP iterations are required with a standard deviation of 1.5. The average optimality gap is 10.2\% with a standard deviation of 23\%. In comparison, the shadow formulation is unable to solve any of these instances and returns an average optimality gap of 56.6\% with a standard deviation of 24\%. Increasing the number of FL from 3 to 5 reduces the density of aircraft per FL which translates into better computational performance for both formulations. All RCP-50-5 and RCP-100-5 instances are solved to optimality by both formulations. For RCP-150-5 instances, the disjunctive formulation is able to solve all instances in an average runtime of 142 s. On the same instances, the shadow formulation times-out on all the instances and returns and average optimality gap of 1.56\%.

Reducing the heading control range to $[-15^\circ,+15^\circ]$ (Table \ref{FL-RCP-15}), yields an overall similar performance. The main changes are observed in the pre-processing procedure which, as in the 2D RCP instances, eliminates an average of 8\% of conflict-free aircraft pairs. We also find that reducing the heading control range increases the number of non-separable pairs to 1.45\% for RCP-150-3 instances and 0.03\% for RCP-150-5 instances. This increase in non-separable aircraft pairs is reflected in the solution of the flight assignment formulation (MILP \ref{mod:fa}) which has an average value of 1.62 on RCP-150-3 instances and 0.02 on RCP-150-5 instances, respectively. For RCP-150-3, 67\% of the instances required a FL change, the maximum number of non-separable pairs ($|\PI|$) is 6 and the maximum objective value of Model \ref{mod:fa} is 3, indicating that three aircraft deviated from their initial FL. In terms of performance, the proposed formulations are able to solve slightly more instances compared to the configuration using the standard heading control range. All RCP-50-3, RCP-50-5 and RCP-100-5 instances are solved to optimality by both formulations. For RCP-100-3 instances, the disjunctive formulation is able to solve all instances in an average time of 168 s, whereas the shadow formulation times-out on 39\% of the instances and requires an average time of 453 s for instances solved to optimality. For RCP-150-3 instances, the disjunctive formulation times-out on 46\% of the instances and for RCP-150-5 the disjunctive formulation is able to solve all instances. In contrast, the shadow formulation is unable to solve any of the RCP-150-3 and RCP-150-5 instances within the available time limit.

\newcommand{\threeflrcpone}{
	\resizebox{0.6\columnwidth}{!}{%
            \begin{tabular}{lll lll ll}
			\toprule
			&&& \multicolumn{3}{l}{Pre-processing} & \multicolumn{2}{l}{FL assignment}\\
			\cmidrule(l){4-6}\cmidrule(l){7-8}
			Instance & $\frac{|\A|}{|\Z|}$ & $\frac{n_c}{|\Z|}$ & $\frac{|\PF|}{|\P|}$ & $\frac{|\PI|}{|\P|}$ & Time & Obj. & Time  \\
			\midrule
RCP-50-3  & 15 & 3.15 (1.5) & 0 & 0.00 & 0.42 & 0.00 (0.0) & 0.00 (0.0)  \\
RCP-100-3 & 33 & 35.6 (1.12) & 0 & 0.00 & 0.86 & 0.00 (0.0) & 0.00 (0.0) \\
RCP-150-3 & 50 & 64.1 (5.26) & 0 & 0.55 & 1.02 & 0.62 (0.08) & 4.56 (0.8)   \\
			\midrule
RCP-50-5  & 10 & 3.11 (1.52) & 0 & 0.00 & 0.62 & 0.00 (0.0) & 0.00 (0.0)  \\
RCP-100-5 & 20 & 14.1 (3.52) & 0 & 0.00 & 1.25 & 0.00 (0.0) & 0.00 (0.0) \\
RCP-150-5 & 30 & 34.2 (2.56)  & 0 & 0.01 & 1.92 & 0.01 (0.1) & 0.11 (1.1) \\
			\bottomrule
		\end{tabular}
	}
}
\newcommand{\threeflrcptwo}{
	\resizebox{1\columnwidth}{!}{%
            \begin{tabular}{l llllll llll l}
			\toprule
			& \multicolumn{6}{l}{Disjunctive}& \multicolumn{4}{l}{Shadow}\\
			\cmidrule(l){2-7}\cmidrule(l){8-11}
			Instance & LB & UB & Gap & Time   & $n_i$ & $n_t$ & $\Delta$UB & Gap  & Time  & $n_t$  & Gain  \\
			\midrule
RCP-50-3 & 1.2E-2 (6E-3) & 1.3E-2 (6E-3) & 0.02 (0.0) & 0.88 (0.2) & 0.0 (0.0) & 0 & 0.00 (0.0) & 0.01 (0.0) & 10.2 (14.8) & 0 & 91.3\\
RCP-100-3 & 1.4E-1 (3E-2) & 1.4E-1 (3E-2) & 0.12 (0.3) & 198 (123) & 1.2 (2.0) & 0 &0.05 (0.1) & 53.2 (0.7) & 572 (17.5) & 46 & 61.9\\
RCP-150-3 & 4.5E-2 (5E-3) & 5.2E-2 (1E-2) & 10.2 (23) & 520 (136) & 2.8 (1.5) & 60&0.03 (0.0) & 56.6 (24) & 600 (0.0) & 100 & - \\
\midrule
RCP-50-5 & 9.4E-3 (5E-3) & 9.4E-3 (5E-3) & 0.00 (0.0) & 0.26 (0.1) & 0.0 (0.0) & 0 & 0.00 (0.0) & 0.00 (0.0) & 1.43 (1.0) & 0 & 81.8\\
RCP-100-5 & 8.3E-3 (3E-2) & 8.3E-2 (3E-2) & 0.00 (0.0) & 2.51 (2.4) & 0.0 (0.0) & 0 & 0.00 (0.0)& 0.13 (0.1) & 126 (163) & 0 & 98.0\\
RCP-150-5 & 3.3E-1 (8E-2) & 3.3E-1 (8E-2) & 0.05 (0.2) & 142 (185) & 2.3 (0.6) & 0 & 0.08 (0.0)& 1.56 (0.6) & 600 (0.0) & 100 & -\\
			\bottomrule
		\end{tabular}}
}

\begin{figure}%
  \centering
  \subfloat{\threeflrcpone}\\%
  \subfloat{\threeflrcptwo}
  \captionof{table}{Results for 2D+FL RCP instances with a speed control range of $[-6\%,+3\%]$ and a heading control range of $[-30^\circ,+30^\circ]$. Times (Time) are reported in seconds. The proportions of conflict-free ($\frac{|\PF|}{|\P|}$) and non-separable ($\frac{|\PI|}{|\P|}$) pairs, optimality gaps (Gap), proportions of time-outs ($n_t$) and the performance gain (Gain) are reported in \%.}%
  \label{FL-RCP-30}%
\end{figure}

\newcommand{\fifteenrcpone}{
	\resizebox{0.6\columnwidth}{!}{%
            \begin{tabular}{lll lll ll}
			\toprule
			&&& \multicolumn{3}{l}{Pre-processing} & \multicolumn{2}{l}{FL assignment}\\
			\cmidrule(l){4-6}\cmidrule(l){7-8}
			Instance & $\frac{|\A|}{|\Z|}$ & $\frac{n_c}{|\Z|}$ & $\frac{|\PF|}{|\P|}$ & $\frac{|\PI|}{|\P|}$ & Time & Obj. & Time  \\
			\midrule
RCP-50-3  & 15 & 3.15 (1.5) & 7.0 & 0.00 & 0.52 & 0.00 (0.0) & 0.00 (0.0)  \\
RCP-100-3 & 33 & 35.6 (1.12) & 7.2 & 0.00 & 0.62 & 0.00 (0.0) & 0.00 (0.0) \\
RCP-150-3 & 50 & 64.1 (5.26) & 8.1 & 1.45 & 1.22 & 1.62 (0.1) & 5.62 (0.8)\\
			\midrule
RCP-50-5  & 10 & 3.11 (1.52) & 8.6 & 0.00 & 0.52 & 0.00 (0.0) & 0.00 (0.0)  \\
RCP-100-5 & 20 & 14.1 (3.52) & 7.8 & 0.00 & 1.22 & 0.00 (0.0) & 0.00 (0.0) \\
RCP-150-5 & 30 & 34.2 (2.56)  & 8.0 & 0.03 & 1.86 & 0.02 (0.2) & 0.21 (1.1)\\
			\bottomrule
		\end{tabular}
	}
}
\newcommand{\fifteenrcptwo}{
	\resizebox{1\columnwidth}{!}{%
            \begin{tabular}{l llllll llll l}
			\toprule
			& \multicolumn{6}{l}{Disjunctive}& \multicolumn{4}{l}{Shadow}\\
			\cmidrule(l){2-7}\cmidrule(l){8-11}
			Instance & LB & UB & Gap & Time  & $n_i$ & $n_t$ & $\Delta$UB & Gap  & Time  & $n_t$  & Gain  \\
			\midrule
RCP-50-3 & 1.2E-2 (6E-3) & 1.2E-2 (6E-3) & 0.02 (0.0) & 0.88 (0.2) & 0.0 (0.0) & 0 & 0.00 (0.0) & 0.01 (0.0) & 10.2 (15) & 0 & 91.3\\
RCP-100-3 & 1.4E-1 (3E-2) & 1.4E-1 (3E-2) & 0.15 (0.3) & 168 (103) & 2.2 (2.0) & 0 &0.08 (0.1) & 43.5 (0.7) & 453 (125) & 39 & 62.4\\
RCP-150-3 & 4.5E-2 (5E-3) & 5.2E-2 (1E-2) & 8.60 (12) & 420 (102) & 4.8 (1.5) & 46 &0.02 (0.0) & 42.1 (20) & 600 (0.0) & 100 & - \\
\midrule
RCP-50-5 & 9.4E-3 (5E-3) & 9.4E-3 (5E-3) & 0.00 (0.0) & 0.26 (0.1) & 0.0 (0.0) & 0 & 0.00 (0.0) & 0.00 (0.0) & 1.43 (1.0) & 0 & 81.8\\
RCP-100-5 & 8.3E-3 (3E-2) & 8.3E-2 (3E-2) & 0.00 (0.0) & 35.1 (2.4) & 0.0 (0.0) & 0 & 0.00 (0.0)& 0.13 (0.1) & 186 (163) & 0 & 81.2\\
RCP-150-5 & 3.3E-1 (8E-2) & 3.3E-1 (8E-2) & 0.51 (0.3) & 152 (78.2) & 4.3 (1.0) & 0 & 0.09 (0.0)& 12.3 (1.6) & 600 (0.0) & 100 & -\\
			\bottomrule
		\end{tabular}}
}

\begin{figure}%
  \centering
  \subfloat{\fifteenrcpone}\\%
  \subfloat{\fifteenrcptwo}
  \captionof{table}{Results for 2D+FL RCP instances with a speed control range of $[-6\%,+3\%]$ and a heading control range of $[-15^\circ,+15^\circ]$. Times (Time) are reported in seconds. The proportions of conflict-free ($\frac{|\PF|}{|\P|}$) and non-separable ($\frac{|\PI|}{|\P|}$) pairs, optimality gaps (Gap), proportions of time-outs ($n_t$) and the performance gain (Gain) are reported in \%.}%
  \label{FL-RCP-15}%
\end{figure}


\section{Conclusion}
\label{con}

We summarize our findings in Section \ref{find} and  discuss future research directions in Section \ref{future}.

\subsection{Summary of Findings}
\label{find}

We proposed new mixed-integer formulations and exact solution methods for  aircraft conflict resolution problems (ACRP). We first considered the 2D ACRP with continuous speed and heading control maneuvers and proposed compact disjunctive separation conditions. The proposed disjunctive separation conditions are linear with regards to aircraft relative velocity variables and only require a single binary variable per pair of aircraft. We have formally shown that the proposed disjunctive linear separation conditions fully characterize the set of conflict-free aircraft trajectories. We introduced a simple pre-processing algorithm to identify aircraft pairs which are conflict-free or non-separable for any combination of controls, which may help in reducing the size of ACRPs by omitting conflict-free pairs. We built on and extended the complex number formulation for the ACRP introduced by \citet{rey2017complex} by augmenting its objective function with a preference weight to balance the trade-off between speed and heading deviations. The resulting formulation is a nonconvex mixed-integer program. This 2D formulation is extended to the context of altitude control by flight level (FL) change and we proposed a lexicographic optimization to solve the 2D+FL ACRP which aims to minimize the number of FL changes in priority and resolve outstanding conflicts by 2D trajectory control.

We presented novel exact solution algorithms to solve the nonconvex 2D and 2D+FL ACRPs. The proposed algorithms refine the convex relaxations introduced by \citet{rey2017complex}. For the 2D problem, the nonconvex formulations are first relaxed to mixed-integer quadratic programs (MIQP) which solution may violate the speed control constraint. If such violations occur, convex quadratic constraints are added together with a constraint generation algorithm that iteratively refines an outer piecewise linear approximation of the speed control constraint by solving a sequence of mixed-integer quadratically constrained programs (MIQCP). The 2D+FL lexicographic optimization problem is solved by decomposing the nonconvex problem into a flight assignment problem and a series of FL-based 2D problems. The proposed flight assignment formulation is based on a reformulation of the FL separation constraint which requires an exponential number of constraints and is embedded into an iterative approach to generate altitude separation constraints as needed.

The performance of the proposed mixed-integer formulations and algorithms was tested on a total of $2072$ benchmarking instances. These instances include four types of ACRPs found in the literature with up to 60 aircraft per instance for 2D problems and 150 aircraft per instance for 2D+FL problems. The performance of the proposed solution algorithms highlights the scalability of the approach compared to existing methods in the literature. Further, we find that the combination of the pre-processing algorithm with the MIQP convex relaxation is sufficient to solve FP and GP instances with up to 60 aircraft, and most RCP instances with up to 30 aircraft. We also find that the number of MIQCP iterations remains low on average when solving larger problems. For 2D+FL lexicographic optimization problems, we find that the pre-processing procedure generates enough altitude separation constraints to solve dense instances with an average of 50 aircraft per FL. 

In our numerical experiments, we also compared the performance of the proposed disjunctive separation conditions with the existing shadow separation conditions. Our results reveal that the proposed disjunctive formulation is able to solve 87.8\% of the benchmarking instances compared to 69.7\% using the shadow formulation. Upon examining all instances that were solved by both formulations, we find that the disjunctive formulation is systematically faster than the shadow formulation, and we report an average decrease in computing time of 86.6\%. This is due to its compact form which only requires one binary variable per aircraft pair while the shadow formulation requires four binary per aircraft pair. For reproducibility purposes, all formulations and instances are made available at the public repository \small \url{https://github.com/acrp-lib/acrp-lib}.\normalsize

\subsection{Future Research and Perspectives}
\label{future}

The proposed approaches for ACRPs rely on several assumptions which may be limiting in practice. One of the modeling assumption at the core of the proposed mixed-integer formulations is the assumption of uniform motion laws, which translates into infinite acceleration and deceleration rates. While such an assumption may be plausible for constrained aircraft speed control, further research is needed to assess the practicality of this assumption when considering varying types of aircraft or different airspace environments, e.g. urban air mobility.

The usage of heading deviations for conflict resolution also raises concerns regarding aircraft trajectory recovery. It is well-acknowledged that resolving conflicts does not guarantee conflict-free recovery trajectories. \citet{alonso2014exact} have proposed an efficient method to calculate the earliest time of recovery for aircraft but this approach does not ensure the existence of conflict-free recovery trajectories. \citet{dias2020two} presented a two-step algorithm for the aircraft conflict resolution with trajectory recovery which decomposes the collision avoidance procedure from the trajectory recovery component but the approach remains sub-optimal, and further research is needed to find global optimal solutions.

Future research may also explore asynchronous conflict resolution problems where aircraft do not require to be coordinated. Further, more efforts to incorporate trajectory prediction uncertainty in ACRPs, which may be caused by adverse weather or disrupted navigation systems, is needed to develop more robust conflict resolution formulations.


\bibliographystyle{spr-chicago}
\bibliography{biblio}

\begin{appendix}

\section{Proofs of the Lemmas, Theorem and Propositions}
\label{landt}

\subsection{Proof of Lemma \ref{lem:normal}:}
\label{pl1}
\begin{proof}
We first show that \eqref{normalplane} is the bisector of one of the two angles formed between $R_1$ and $R_2$ (see Figure \ref{fig:qznp}).
Note that the slope of the plane defined by \eqref{normalplane} is $\frac{\pzy}{\pzx}$. Without any loss of generality, assume that $|\pzx| > d$ (the case $|\pzx| \le d$ can be treated similarly).

Let $r_1 = \frac{\pzx \pzy + d \sqrt{\pzx^2 + \pzy^2 - d^2}}{\pzx^2 - d^2}$ and $r_2 = \frac{\pzx \pzy - d \sqrt{\pzx^2 + \pzy^2 - d^2}}{\pzx^2 - d^2}$ be the slopes of the lines defined by \eqref{l3} and \eqref{l4}. The angle of the bisector of these lines is $\mu = \frac{1}{2}(\arctan{r_1} + \arctan{r_2})$ and with slope:
\begin{align*}
\tan(\mu) &= \tan\left(\frac{1}{2}\left(\arctan(r_1) + \arctan(r_2\right)\right) = \tan\left(\frac{1}{2}\left(\arctan\left(\frac{r_1 + r_2}{1 - r_1 r_2}\right)\right)\right). 
\end{align*}
If $r_1 r_2 = 1$, recall that $\lim_{X \rightarrow \pm\infty} \arctan(X) = \pm\pi/2$, thus $\tan(\mu) = \tan(\frac{1}{2}\frac{\pm\pi}{2}) = \tan(\pm \pi/4) = \pm 1$. In addition, $r_1 r_2 = \frac{\pzy^2 - d^2}{\pzx^2 - d^2}$, hence $r_1 r_2 = 1 \Leftrightarrow \pzy^2 = \pzx^2$ and the slope of \eqref{normalplane} is $\pm 1$.

Assume now $r_1 r_2 \ne 1$. Combining the half-angle formula $\tan(\frac{X}{2}) = \frac{1 - \cos(X)}{\sin(X)}$ with the inverse formulas: $\cos(\arctan(X))=\frac{1}{\sqrt{1+X^2}}$ and $\sin(\arctan(X))=\frac{X}{\sqrt{1+X^2}}$, after simplification $\tan(\mu)$ can be written as:
\[
\tan(\mu) = \frac{\sqrt{1 + \left(\frac{r_1 + r_2}{1 - r_1 r_2}\right)^2} - 1}{\left(\frac{r_1 + r_2}{1 - r_1 r_2}\right)}.
\]
Since $\frac{r_1 + r_2}{1 - r_1 r_2} = \frac{2\pzx\pzy}{\pzx^2 - \pzy^2}$, this gives:
\[
\tan(\mu) = \frac{\sqrt{(\pzx^2 - \pzy^2)^2 + 4\pzx\pzy} - (\pzx^2 - \pzy^2)}{2\pzx\pzy} = \frac{\pzy}{\pzx}.
\]
Since \eqref{plane} is orthogonal to \eqref{normalplane}, the line \eqref{plane} is the bisector of the other angle between the two linear equations represented by $R_1$ and $R_2$. 
\end{proof}

\subsection{Proof of Lemma \ref{lem:negative}:}
\label{pl2}

\begin{proof}
Any point $(\vijx,\vijy)$ of \eqref{normalplane} verifies	$\vijy = \frac{\hat{y}_{ij} }{\hat{x}_{ij} }\vijx$. Replacing $\vijy$ in \eqref{gfunction} yields:
\begin{equation*}
\begin{split}
g\Big(\vijx,\frac{\hat{y}_{ij}}{\hat{x}_{ij}}\vijx\Big) &= \vijx^2(\hat{y}_{ij}^2 - d^2) + \Big(\frac{\hat{y}_{ij} }{\hat{x}_{ij} }\vijx\Big)^2(\hat{x}_{ij}^2  - d^2)\\
&- \vijx\frac{\hat{y}_{ij}}{\hat{x}_{ij}}\vijx(2\hat{x}_{ij}\hat{y}_{ij}).
\end{split}
\end{equation*}

Simplifying the previous expression yields:
\begin{equation*}
g\Big(\vijx,\frac{\hat{y}_{ij} }{\hat{x}_{ij} }\vijx\Big) = -d^2v^2_{ij,x}\Big(1 + \Big(\frac{\hat{y}_{ij} }{\hat{x}_{ij}}\Big)^2\Big) \leq 0.
\end{equation*}
\end{proof}

\subsection{Proof of Theorem \ref{theo:conditions}:}
\label{pt}
\begin{proof}
We prove this statement by showing that the conditions \eqref{eq:zd}-\eqref{eq:rs} are equivalent to the nonlinear conditions \eqref{eq:sepconditions} which are well-known to be equivalent to Eq. \eqref{eq:sepcond}. As shown by Lemmas \ref{lem:normal} and \ref{lem:negative}, the line \eqref{normalplane} splits the $(\vijx,\vijy)$-plane into two disjunctive regions and the region in which lies the normal line \eqref{normalplane} consists of conflicting trajectories. Thus, both half-planes induced by the normal plane \eqref{normalplane} contain sub-regions corresponding to conflict-free and conflicting trajectories, and Eq. \eqref{eq:zd} defines variables $\z$ accordingly. Consider the half-plane corresponding to $\z = 1$ (see Figure \ref{fig:qznp}). This half-plane can be further split into two convex sub-regions and inequality \eqref{eq:r1} characterizes convex conflict-free region in which all pairs of aircraft trajectories verify $g(\vijx,\vijy) \geq 0$ or $\tm \leq 0$. The same reasoning applies to the half-plane corresponding to $\z = 0$ if we substitute \eqref{eq:r1} by \eqref{eq:r2}. Hence, all pairs of future conflict-free trajectories for aircraft $i$ and $j$ are characterized by the separation conditions \eqref{eq:zd}-\eqref{eq:rs}.
\end{proof}

\subsection{Proof of Proposition \ref{prop:cf}:}
\label{pp1}
\begin{proof}
For any pair $(i,j) \in \P$, the relative velocity box $\B$ and the conflict region $\C$ are both convex sets, and $\B \cap \C$ can be fully characterized by linear inequalities. The feasible region of the proposed feasibility linear program $LP(i,j)$ represents the set $\left\{(\vijx,\vijy) \in \B \cap \C \right\}$. If this set is empty, then there exists no pair of feasible trajectories for aircraft $i$ and $j$ which correspond to a conflict, hence the pair of aircraft $(i,j) \in \P$ is conflict-free. 
\end{proof}

\subsection{Proof of Proposition \ref{prop:ns}:}
\label{pp2}
\begin{proof}
Since the conflict region $\C$ is the convex, if all four extreme points of $\B$ are inside $\C$, i.e. each extreme points of $E(\B)$ corresponds to a conflict, then the aircraft pair $(i,j) \in \P$ is non-separable for any controls.
\end{proof}

\subsection{Proof of Proposition \ref{prop:conv}:}
\label{pp3}
\begin{proof}
Let $c_i (\dx,\dy) = w\dy^2 + (1 - w)(1 - \dx)^2$ be the cost function of aircraft $i \in \A$. The first-order optimality conditions of $c_i(\dx,\dy)$ are:
\begin{align*}
\frac{\partial c_i(\dx,\dy)}{\partial \dx} &= (2\dx - 2)(1 - w) = 0 \\
\frac{\partial c_i(\dx,\dy)}{\partial \dy} &= 2w\dy =0
\end{align*}

Recall that $w \in \ ]0,1[$, hence the first-order optimality conditions of this objective function yield $\dx=1$ and $\dy=0$ for all aircraft $i \in \A$. Since $\dy=\qi\sin(\ti)$, and $\qi > 0$, this implies $\dy=0$ which is equivalent to $\ti=0$ if $\ti \in \ ]-\pi,\pi[$. Further, since $\dx=\qi\cos(\ti)$ and $\ti=0$, $\dx=1$ implies $\qi=1$. It is trivial to show that $c_i (\dx,\dy)$ is convex, thus the 2D objective function \eqref{eq:obj2d} is also convex.


\end{proof}
\section{Details of the implementation of the shadow formulation}
\label{shadow}

Our implementation of the shadow formulation as proposed by \cite{bilimoria2000geometric} and \cite{pallottino2002conflict} and further developed by \cite{alonso2011collision,alonso2016exact} is summarized below. Let $\sig^1, \sig^2, \sig^3, \sig^4 \in \{0,1\}$ be binary decision variables for all aircraft pair $(i,j) \in \PS$. The separation constraints are defined as:
\begin{subequations}
\begin{align}
& -v_{ij,x} \leq M^s_{ij}(1-\sig^1), &&\forall (i,j) \in \PS, \\
& v_{ij,x}\tan(l_{ij}) - v_{ij,y} \leq (1 - \sig^1)M^{tanl}_{ij}, &&\forall (i,j) \in \PS, \\
& -v_{ij,x} \leq M^s_{ij}(1-\sig^2), &&\forall (i,j) \in \PS, \\
& -v_{ij,x}\tan(r_{ij}) + v_{ij,y} \leq (1 - \sig^2)M^{tanr}_{ij}, &&\forall (i,j) \in \PS, \\
& v_{ij,x} \leq M^s_{ij}(1 - \sig^3), &&\forall (i,j) \in \PS, \\
& -v_{ij,x}\tan(l_{ij}) + v_{ij,y} \leq (1 - \sig^3)M^{tanl}_{ij}, &&\forall (i,j) \in \PS, \\
& v_{ij,x} \leq M^s_{ij}(1 - \sig^4), &&\forall (i,j) \in \PS, \\
& v_{ij,x}\tan(r_{ij}) - v_{ij,y} \leq (1 - \sig^4)M^{tanr}_{ij}, &&\forall (i,j) \in \PS, \\
& \sig^1 + \sig^2 + \sig^3 + \sig^4 \geq 1, &&\forall (i,j) \in \PS. 
\end{align}\label{eq:sep_shadow}
\end{subequations}

In our experiments, we embed the shadow separation constraints \eqref{eq:sep_shadow} within the proposed formulations by substituting the disjunctive separation constraints \eqref{eq:zd} and \eqref{eq:rs} with \eqref{eq:sep_shadow}.  We summarize the nonconvex complex number formulation for the 2D ACRP based on the shadow separation conditions in Model \ref{mod:2dshadow}

\begin{model}[Nonconvex 2D Formulation using Shadow Separation Conditions]
\label{mod:2dshadow}
\begin{subequations}
\begin{align}
&\emph{Minimize Total 2D Deviation} \quad \eqref{eq:obj2d}, \nonumber \\
&\emph{Subject to:}  && \nonumber \\
&\emph{Motion Equations} \quad \eqref{eq:speedxdy}, \nonumber\\
&\emph{Separation Conditions} \quad \eqref{eq:sep_shadow}, \nonumber \\
&\emph{Speed and Heading Control Constraints} \quad  \eqref{eq:boundq}, \eqref{eq:boundtheta}, \nonumber\\
&\emph{Variable Bounds} \quad \eqref{eq:boundsdxdy}, \nonumber \\
& \vijx,\vijy \in \B, &&\forall (i,j) \in \PS, \\
& \sig^1, \sig^2, \sig^3, \sig^4 \in \{0,1\}, &&\forall (i,j) \in \PS, \\
& \dx, \dy \in \R, && \forall i \in \A. 
\end{align}
\end{subequations}
\end{model}










\end{appendix}

\end{document}